%% file: least-squares-convection-diffusion.tex
\documentclass[12pt]{article}

\usepackage{a4wide,epsfig}
\usepackage{amsfonts,amsmath}
\usepackage{caption}
\usepackage[labelformat=simple]{subcaption}

\usepackage{enumerate}
\usepackage{pgfplotstable}
\usepackage{booktabs}
\usepackage{stmaryrd}
\usepackage{adjustbox}

\newtheorem{theorem}{Theorem}[section]

\newtheorem{lemma}[theorem]{Lemma}

\newtheorem{remark}{Remark}[section]
\newtheorem{example}{Example}[section]
\newcommand{\proof} [1]
   { \noindent {\bf Proof.} #1 \hfill\rule{0.5em}{1.2ex} \par\medskip}

\newcommand{\norm}[1]{{\left\lVert{#1}\right\rVert}}

\DeclareMathOperator{\osc}{osc}
\DeclareMathOperator{\divv}{div}

\numberwithin{equation}{section} 

\begin{document}

\setcounter{page}{1}

\title{Adaptive least-squares space-time finite element methods for 
  convection-diffusion problems}
\author{Christian~K\"othe, Olaf~Steinbach}
\date{Institut f\"ur Angewandte Mathematik, TU Graz, \\[1mm]
  Steyrergasse 30, 8010 Graz, Austria}

\maketitle

\centerline{\em In memoriam Raytcho D. Lazarov (1943--2024)}

\bigskip

\begin{abstract}
  In this paper we formulate and analyse adaptive (space-time)
  least-squares finite element methods for the solution of
  convection-diffusion equations. The convective derivative
  $\mbox{\boldmath $v$} \cdot \nabla u$ is
  considered as part of the total time derivative
  $\frac{d}{dt}u = \partial_t u + \mbox{\boldmath $v$} \cdot \nabla u$,
  \linebreak
  and therefore we can use a rather standard stability and error analysis
  for related space-time finite element methods. For stationary problems we
  restrict the ansatz space $H^1_0(\Omega)$ such that the convective
  derivative is considered as an element of the dual $H^{-1}(\Omega)$ of
  the test space $H^1_0(\Omega)$, which also allows unbounded
  velocities $\mbox{\boldmath $v$}$. While the discrete finite element
  schemes are always unique solvable, the numerical solutions may
  suffer from a bad approximation property of the finite element space
  when considering convection dominated problems, i.e., small diffusion
  coefficients. Instead of adding suitable stabilization terms, we 
  aim to resolve the solutions by using adaptive (space-time) finite
  element methods. For this we introduce a least-squares approach where
  the discrete adjoint defines local a posteriori error indicators to drive
  an adaptive scheme. Numerical examples illustrate the theoretical
  considerations.
\end{abstract}



\section{Introduction}
As documented by Bochev and Gunzburger \cite{BochevGunzburger:2009},
least-squares finite element methods are a well established approach for the
numerical solution of second order partial differential equations. In most
cases, the partial differential equation is rewritten as first order system,
and the residuals of both equations are minimized, using appropriate norms
and weights. For early contributions of R.~Lazarov using the $L^2(\Omega)$
norm for both residuals, see
\cite{CaiLazarovManteuffelMcCormick:1994,PehlivanovCareyLazarov:1994}, and
\cite{BrambleLazarovPasciak:1997,BrambleLazarovPasciak:1998} when considering
the residual of the equilibrium equation in $H^{-1}(\Omega)$. Later on,
see \cite{LazarovTobiskaVassilevski:1997}, a least-squares formulation
in $L^2(\Omega)$ for a first-order system for a convection dominated
convection-diffusion equation was analyzed, for a related minimization
approach in $H^{-1}(\Omega)$, see \cite{LazarovVassilevski:1999}. For more
recent contributions on first order least-squares systems we refer to, e.g.,
\cite{Fuehrer:2024,FuehrerHeuerKarkulik:2022,SchafelnerVassilevski:2021}.

For a stable and accurate numerical solution of 
diffusion-convection-reaction equations there exists a huges amount of
literature, here we mention, e.g.,
\cite{AugustinJohn:2011,Bacuta:2024,BacutaHayes:2024,BrooksHughes:1982,
  CohenDahmenWelper:2012,DemkowiczHeuer:2013,
  FrancaHaukeMasud:2006,RoosStynesTobiska:2008}, just to name a few.
In the particular context of space-time variational formulations we
refer to, e.g., \cite{Bank:2017,BurmanErn:2012,Toulopoulos:2025}.

In our recent work \cite{KLS:2023} we have formulated and analysed a 
least-squares approach for the numerical solution of rather general
operator equations, including elliptic, parabolic, and hyperbolic
partial differential equations, for related work on parabolic
evolution equations, see also
\cite{Andreev:2013,Dahmen:2023,Monsuur:2024,Stevenson:2021}.
Instead of rewriting a second order partial differential equation as
first order system, which is then solved by a least-squares approach,
we consider a least-squares approach for the original equation minimizing
the residual in appropriate norms which are induced by elliptic and
self-adjoint operators. This requires the solution of a saddle point
variational formulation, where we can use standard arguments as known
in mixed finite element methods \cite{Boffi:2013}. While the adjoint
variable turns out to be zero in the continuous setting, its discrete
approximation can be used to define local a posteriori error indicators
to drive an adaptive scheme.

In this work we extend the approach of \cite{KLS:2023} to solve stationary
and instationary convection-diffusion equations by using adaptive
(space-time) finite element methods. For parabolic evolution equations
we combine the convective derivative with the partial time derivative
and introduce the total time derivative. The resulting and standard space-time
finite element method can be analyzed as in
\cite{GanglGobrialSteinbach:2025,Steinbach:2015}, but here we present
a more general proof of surjectivity when a convective derivative appears.
For an adaptive space-time finite element scheme we then apply the
least-squares approach as in \cite{KLS:2023}. In the case of an elliptic
convection-diffusion equation, and motivated by the parabolic case,
we consider the first order convective term in the dual $H^{-1}(\Omega)$
of the variational test space $H^1_0(\Omega)$. This allows to consider
even unbounded velocities. For a direct finite element discretization
we present a related numerical analysis showing the relations between
the size of the diffusion coefficient, and the finite element mesh size
which has to be sufficiently small in the case of convection dominated
problems. Instead of using appropriate stabilization techniques, here we
aim to construct accurate solutions by using an adaptive least-squares
approach.

\section{Stationary convection-diffusion problems}
For a bounded Lipschitz domain $\Omega \subset {\mathbb{R}}^n$, $n=1,2,3$,
we consider, as a model problem, the Dirichlet boundary value problem
for the convection-diffusion equation,
\begin{equation}\label{DBVP Advection Diffusion}
  - \mbox{div} [\alpha(x) \nabla u(x) ] +
  \mbox{\boldmath $v$}(x) \cdot \nabla u(x) = f(x) \quad
  \mbox{for} \; x \in \Omega, \quad
  u(x) = 0 \quad \mbox{for} \; x \in \partial \Omega ,
\end{equation}
where $f(x)$, $x \in \Omega$, is given.
We assume that the diffusion coefficient $\alpha(x)$ is bounded and strictly
positive, i.e.,
\begin{equation}\label{Assumption alpha}
  0 < \underline{\alpha} \leq \alpha(x) \leq \overline{\alpha} \quad
  \mbox{for all} \; x \in \Omega .
\end{equation}
Assumptions on the given velocity field $\mbox{\boldmath $v$}$ will be
specified later on. Note that in \eqref{DBVP Advection Diffusion} we may
add some reaction term $c(x)u(x)$ with some non-negative function
$c(x)$, $x \in \Omega$, and we may consider inhomogeneous Dirichlet
boundary conditions, or boundary conditions of mixed type, including
Neumann and Robin type boundary conditions.

\subsection{Variational formulation}
The variational formulation of the Dirichlet boundary value problem
\eqref{DBVP Advection Diffusion} is to find $u \in X$ such that
\begin{equation}\label{VF DBVP Advection Diffusion}
  b(u,q) := \int_\Omega \alpha(x) \,  \nabla u(x) \cdot
  \nabla q(x) \, dx + \int_\Omega
  \mbox{\boldmath $v$}(x) \cdot \nabla u(x) \, q(x) \, dx =
  \int_\Omega f(x) q(x) \, dx
\end{equation}
is satisfied for all $ q \in Y := H^1_0(\Omega)$, where the related
energy norm of the test space $Y$ is given by
\[
  \| q \|^2_Y := \int_\Omega \alpha(x) \, |\nabla q(x)|^2 \, dx ,
\]
and for which we have the norm equivalence inequalities
\begin{equation}\label{Norm equivalence Y}
  \sqrt{\underline{\alpha}} \, \| \nabla q \|_{L^2(\Omega)} \leq
  \| q \|_Y \leq
  \sqrt{\overline{\alpha}} \, \| \nabla q \|_{L^2(\Omega)}
  \quad \mbox{for all} \; q \in Y = H^1_0(\Omega) .
\end{equation}
According to the variational formulation \eqref{VF DBVP Advection Diffusion}
we assume $f \in H^{-1}(\Omega) := [H^1_0(\Omega)]^*$ with the norm
\begin{equation}\label{dual norm}
  \| f \|_{Y^*} := \sup\limits_{0 \neq q \in Y}
  \frac{\langle f , q \rangle_\Omega}{\| q \|_Y} ,
\end{equation}
using the duality pairing $\langle f , q \rangle_\Omega$
for $f \in Y^*$ and $q \in Y$ as extension of the $L^2(\Omega)$ inner product.
While the most standard choice for the ansatz space is
$X = H^1_0(\Omega)$, here we will use a slightly different approach.
We introduce the velocity dependent ansatz space
$X_v := \{ u \in Y : \mbox{\boldmath $v$} \cdot \nabla u \in Y^* \} \subset Y$,
with the graph norm
\[
  \| u \|_{X_v} := \sqrt{\| u \|_Y^2 +
    \| \mbox{\boldmath $v$} \cdot \nabla u \|_{Y^*}^2} =
  \sqrt{ \| u \|^2_Y + \| w_u \|^2_Y},
\]
where the Riesz representant $w_u \in Y$ is the unique solution of
the variational problem
\begin{equation}\label{Definition wu}
  a(w_u,q) :=
  \int_\Omega \alpha(x) \, \nabla w_u(x) \cdot \nabla q(x) \, dx =
  \int_\Omega \mbox{\boldmath $v$}(x) \cdot \nabla u(x) \, q(x) \, dx
  \quad \mbox{for all} \; q \in Y .
\end{equation}
Depending on the regularity of the given velocity $\mbox{\boldmath $v$}$,
the condition $\mbox{\boldmath $v$} \cdot \nabla u \in Y^*$ may result in the
fact that $X_v \subset Y$ is a real subspace, i.e., $X_v \neq Y$. At this
time we only assume that $X_v$ is not empty.
By the definition of the underlying function spaces we then conclude
\begin{eqnarray}
  |b(u,q)| \nonumber
  & = & \left| \int_\Omega \alpha(x) \,  \nabla u(x) \cdot
        \nabla q(x) \, dx + \int_\Omega
        \mbox{\boldmath $v$}(x) \cdot \nabla u(x) \, q(x) \, dx \right| \\
  & = & \nonumber
        \left| \int_\Omega \alpha(x) \,  \nabla u(x) \cdot
        \nabla q(x) \, dx + \int_\Omega
        \alpha(x) \, \nabla w_u(x) \cdot \nabla q(x) \, dx \right| \\
  & \leq & \nonumber \Big[ \| u \|_Y + \| w_u \|_Y \Big] \, \| q \|_Y \\
  & \leq & \sqrt{2} \, \sqrt{ \| u \|_Y^2 + \| w_u \|_Y^2} \, \| q \|_Y
           \, = \, \sqrt{2} \, \| u \|_{X_v} \| q \|_Y \quad
           \label{boundedness b}
\end{eqnarray}
for all $(u,q) \in X_v \times Y$, i.e., boundedness of the bilinear form
$b(u,q)$. A first rather standard assumption on the given velocity
$\mbox{\boldmath $v$}$ is as follows.

\begin{lemma}\label{Lemma convective positiv}
  Assume
  $\mbox{\rm div} \, \mbox{\boldmath $v$}(x) \leq 0$ for almost all
  $x \in \Omega$. Then,
  \begin{equation}\label{vuu positiv}
    \int_\Omega \mbox{\boldmath $v$}(x) \cdot \nabla u(x) \, u(x) \, dx
    \geq 0 \quad \mbox{for all} \; u \in X_v \subset Y = H^1_0(\Omega) .
  \end{equation}
\end{lemma}
\proof{For $ u \in X_v \subset Y$ we first note that the left hand side
  in \eqref{vuu positiv} is bounded.
  Using integration by parts for $ u \in X_v \subset Y = H^1_0(\Omega)$,
  \begin{eqnarray*}
    \int_\Omega \mbox{\boldmath $v$}(x) \cdot \nabla u(x) \, u(x) \, dx
    & = & \sum\limits_{k=1}^n
          \int_\Omega \partial_{x_k} u(x) \, v_k(x) u(x) \, dx \\
    & = & - \sum\limits_{k=1}^n
          \int_\Omega u(x) \, \partial_{x_k} [v_k(x)u(x)] \, dx \\
    & = & - \sum\limits_{k=1}^n \int_\Omega \Big( [u(x)]^2
          \partial_{x_k} v_k(x) +
          u(x) \, v_k(x) \, \partial_{x_k} u(x) \Big) \, dx,
  \end{eqnarray*}
  i.e.,
  \[
    2 \, \int_\Omega \mbox{\boldmath $v$}(x) \cdot \nabla u(x) \, u(x) \, dx =
    - \int_\Omega [u(x)]^2 \, \mbox{div} \,
    \mbox{\boldmath $v$}(x) \, dx \, \geq \, 0 .
  \]
}

\noindent
In order to ensure uniqueness for the solution of the
variational formulation \eqref{VF DBVP Advection Diffusion} we need to
have the following result.

\begin{lemma}
  Assume $\mbox{\rm div} \, \mbox{\boldmath $v$}(x) \leq 0$ for almost all
  $x \in \Omega$. Then the bilinear form $b(\cdot,\cdot)$ as defined in
  \eqref{VF DBVP Advection Diffusion} satisfies the inf-sup
  stability condition
  \begin{equation}\label{inf-sup advection-diffusion}
    \| u \|_{X_v} \leq \sup\limits_{0 \neq q \in Y} \frac{b(u,q)}{\| q \|_Y}
    \quad \mbox{for all} \; u \in X_v.
  \end{equation}
\end{lemma}

\proof{For an arbitrary but fixed $ u \in X_v \subset Y$ and the solution
  $w_u \in Y$ of the variational formulation \eqref{Definition wu}
  we define $q_u := u + w_u \in Y$, and we obtain
  \begin{eqnarray*}
    b(u,q_u)
    & = & \int_\Omega \alpha(x) \,  \nabla u(x) \cdot
          \nabla q_u(x) \, dx + \int_\Omega
          \mbox{\boldmath $v$}(x) \cdot \nabla u(x) \, q_u(x) \, dx \\
    & = & \int_\Omega \alpha(x) \,  \nabla u(x) \cdot
          \nabla q_u(x) \, dx + \int_\Omega
          \alpha(x) \, \nabla w_u(x) \cdot \nabla q_u(x) \, dx \\
    & = & \int_\Omega \alpha(x) \, \nabla [u(x)+w_u(x)]
          \cdot \nabla q_u \, dx \\
    & = & \int_\Omega \alpha(x) \, |\nabla q_u(x)|^2 \, dx \, = \,
          \| q_u \|^2_Y \, .
  \end{eqnarray*}
  On the other hand we have, using \eqref{vuu positiv},
  \begin{eqnarray*}
    \| q_u \|_Y^2
    & = & \| u + w_u \|_Y^2 \\
    && \hspace*{-15mm}
       = \, \int_\Omega \alpha(x) \, |\nabla_x[u(x)+w_u(x)]|^2 \, dx \\
    && \hspace*{-15mm}
       = \, \int_\Omega \alpha(x) \, |\nabla_x u(x)|^2 \, dx
          + \int_\Omega \alpha(x) \, |\nabla_x w_u(x)|^2 \, dx
          + 2 \,
          \int_\Omega \alpha(x) \, \nabla_x w_u(x) \cdot \nabla u(x) \, dx \\
    && \hspace*{-15mm}
       = \, \| u \|_Y^2 + \| w_u \|_Y^2 + 2 \,
          \int_\Omega \mbox{\boldmath $v$}(x) \cdot \nabla u(x) \, u(x) \, dx \\
    && \hspace*{-15mm} \geq \, \| u \|_{X_v}^2 .
  \end{eqnarray*}
  Hence, we have shown the inf-sup stability condition
  \[
    \| u \|_{X_v} \leq \frac{b(u,q_u)}{\| q_u \|_Y} \leq
    \sup\limits_{0 \neq q \in Y} \frac{b(u,q)}{\| q \|_Y} \quad
    \mbox{for all} \; u \in X_v.
  \]
}

\noindent
It remains to prove solvability of the variational formulation
\eqref{VF DBVP Advection Diffusion}, i.e., surjectivity of the
related operator $B : X_v \to Y^*$ which is defined by
\[
  \langle B u , q \rangle_\Omega := b(u,q) \quad
  \mbox{for all} \; (u,q) \in X_v \times Y.
\]
Recall that the definition of $X_v$ depends
on the regularity of the velocity field $\mbox{\boldmath $v$}$.
If the velocity is bounded, which is a standard assumption in many
applications, i.e., $\|\mbox{\boldmath $v$}(x)\|_2 \leq c_v$
for almost all $x \in \Omega$, we then obtain
\[
  \|
  \mbox{\boldmath $v$} \cdot \nabla u \|_{Y^*}
  = \sup\limits_{0 \neq q \in Y}
  \frac{\langle \mbox{\boldmath $v$} \cdot \nabla u,
    q \rangle_\Omega}{\| q \|_Y} \, \leq \, c_v \, \| \nabla u \|_{L^2(\Omega)}
  \sup\limits_{0 \neq q \in Y}
  \frac{\| q \|_{L^2(\Omega)}}{\| q \|_Y} \, \leq \, c \, \| u \|_Y,
\]
where we have used \eqref{Norm equivalence Y}, and Friedrich's inequality.
In this case we conclude that $\| u \|_Y$ defines an equivalent norm
in $X_v$, and $X_v = Y = H^1_0(\Omega)$ follows for all bounded
velocities $\mbox{\boldmath $v$}$. In any case, the assumption on the
boundedness of $\mbox{\boldmath $v$}$ can be relaxed as follows.

\begin{lemma}
  Assume
  \begin{equation}\label{Assumption v}
    c_v := \sup\limits_{0 \neq q \in H^1_0(\Omega)}
    \frac{\| q \, \mbox{\boldmath $v$}\|_{L^2(\Omega)}}
    {\| \nabla q \|_{L^2(\Omega)}} < \infty .
  \end{equation}
  Then,
   \begin{equation}\label{bound v p}
    \| \mbox{\boldmath $v$} \cdot \nabla u \|_{Y^*} \, \leq \, 
    \frac{c_v}{\underline{\alpha}} \;
    \| u \|_Y \quad \mbox{for all} \; u \in X_v .
  \end{equation}
\end{lemma}

\proof{When using the definition \eqref{dual norm} we have
  \begin{eqnarray*}
    \| \mbox{\boldmath $v$} \cdot \nabla u \|_{Y^*}
    & = & \sup\limits_{0 \neq q \in Y}
          \frac{\langle \mbox{\boldmath $v$} \cdot \nabla u,
          q \rangle_\Omega}{\| q \|_Y} \, = \,
          \sup\limits_{0 \neq q \in Y}
          \frac{\langle \nabla u,
          q \, \mbox{\boldmath $v$} \rangle_\Omega}{\| q \|_Y}  \\
    & \leq & \| \nabla u \|_{L^2(\Omega)}
             \sup\limits_{0 \neq q \in Y}
             \frac{\| q \, \mbox{\boldmath $v$} \|_{L^2(\Omega)}}{\| q \|_Y} \\
    & \leq & \frac{1}{\underline{\alpha}} \, \| u \|_Y
             \sup\limits_{0 \neq q \in Y}
             \frac{\| q \, \mbox{\boldmath $v$} \|_{L^2(\Omega)}}
             {\| \nabla q \|_{L^2(\Omega)}} \, = \,
             \frac{c_v}{\underline{\alpha}} \, \| u \|_Y \, .
  \end{eqnarray*}
}

\noindent
As a consequence of \eqref{bound v p} we conclude the norm equivalence
inequalities
\[
  \| u \|_Y^2 \leq \| u \|_{X_v}^2 =
  \| u \|_Y^2 + \| \mbox{\boldmath $v$} \cdot \nabla u \|^2_{Y^*}
  \leq \left( 1 + \frac{c_v^2}{\underline{\alpha}^2}
  \right) \, \| u \|^2_Y
\]
for all $u \in X_v$,
i.e., $\| u \|_Y$ defines an equivalent norm in $X_v$, and hence,
$X_v=Y=H^1_0(\Omega)$ follows when \eqref{Assumption v} is satisfied.
With this we are now in the position to state the surjectivity of $B$.

\begin{lemma}
  Assume $\mbox{\rm div} \, \mbox{\boldmath $v$}(x) \leq 0$ for almost all
  $x \in \Omega$ and \eqref{Assumption v} to be satisfied. For any
  $p \in Y \backslash \{ 0 \}$ there exists a $ u_p \in X_v$ such that
  \[
    b(u_p,p) > 0 \, .
  \]
\end{lemma}

\proof{Due to \eqref{Assumption v} we have $0 \neq p \in Y = X_v$, and
  hence we can choose $u_p = p \in X_v$ to conclude
  \[
    b(u_p,p) = \int_\Omega \alpha(x) \, \nabla p(x) \cdot \nabla p(x) \, dx
    + \int_\Omega \mbox{\boldmath $v$}(x) \cdot \nabla p(x) \, p(x) \, dx
    \, \geq \, \| p \|^2_Y > 0 .
  \]
}

\noindent
Now we are in a position to state the unique solvability result for
the solution of the variational formulation
\eqref{VF DBVP Advection Diffusion}.

\begin{lemma}\label{Lemma solvability}
  Assume \eqref{Assumption v} and
  $\mbox{\rm div} \, \mbox{\boldmath $v$}(x) \leq 0$ for almost all
  $x \in \Omega$. For any $f \in Y^*$ there exists a unique solution
  $u \in X_v$ of the variational formulation
  \eqref{VF DBVP Advection Diffusion} satisfying
  \[
    \| u \|_{X_v} \leq \| f \|_{Y^*} \, .
  \]
\end{lemma}

\proof{Since the bilinear $b(u,q)$ is bounded for all
  $(u,q) \in X_v \times Y$, satisfies the inf-sup stability
  condition \eqref{inf-sup advection-diffusion}, and is surjective,
  all assumptions of the Babu\v{s}ka--Ne\v{c}as theory, e.g.,
  \cite{Babuska:1971,Braess:2007,ErnGuermond:2004,Necas:1962}, are satisfied,
  and unique solvability of \eqref{VF DBVP Advection Diffusion} follows.
  Moreover, using \eqref{inf-sup advection-diffusion} we have
  \[
    \| u \|_{X_v} \leq \sup\limits_{0 \neq q \in Y}
    \frac{b(u,q)}{\| q \|_Y} =
    \sup\limits_{0 \neq q \in Y} \frac{\langle f , q \rangle_\Omega}{\| q \|_Y}
    = \| f \|_{Y^*} \, .
  \]
}

\begin{example}\label{rem:choice of v for n=1}
  For $n=1$, $\Omega = (0,1)$ and $\mbox{\boldmath $v$}(x) = v(x)$
  we consider the Dirichlet boundary value problem
  \[
    - u''(x) + v(x) u'(x) = f(x) \quad \mbox{for} \; x \in (0,1),
    \quad u(0)=u(1)=0,
  \]
  and the variational formulation to find $u \in X_v$ such that
  \begin{equation}\label{DBVP Example 1D}
    b(u,q) := \int_0^1 u'(x) q'(x) \, dx +
    \int_0^1 v(x) u'(x) q(x) \, dx =
    \int_0^1 f(x) q(x) \, dx
  \end{equation}
  is satisfied for all $q \in Y = H^1_0(0,1)$. For $q \in H^1_0(0,1)$
  we can write
  \[
    q(x) = \int_0^x q'(s) \, ds,
  \]
  and
  \begin{eqnarray*}
    |v(x) \, q(x)|^2
    & = & \left| v(x) \int_0^x q'(s) \, ds \right|^2 \\
    & \leq & [v(x)]^2 \int_0^x ds \int_0^x [q'(s)]^2 \, ds \, \leq \,
             x \, [v(x)]^2 \int_0^1 [q'(s)]^2 \, ds,
  \end{eqnarray*}
  and integration over $x \in (0,1)$ gives
  \[
    \| v \, q \|_{L^2(0,1)}^2 \leq
    \int_0^1 x \, [v(x)]^2 dx \; \| q' \|_{L^2(0,1)}^2,
  \]
  i.e.,
  \begin{equation}\label{Bound cv}
    c_v \leq \left( \int_0^1 x \, [v(x)]^2 \, dx \right)^{1/2} \, .
  \end{equation}
  Hence we conclude unique solvability of the variational
  formulation \eqref{DBVP Example 1D} for velocities $v(x)$ satisfying
  \eqref{Bound cv}. For example, we can consider
  \[
    v(x) = \frac{1}{\sqrt{x}}, \quad
    c_v \leq 1, \quad \mbox{\rm div} \, \mbox{\boldmath $v$}(x) =
    \frac{d}{dx} \frac{1}{\sqrt{x}} = - \frac{1}{2} \, x^{-3/2} < 0
    \quad \mbox{for} \; x \in (0,1).
  \]
  But for $v(x)=1/x$ we can not apply the above estimates, although we have
  \[
    \mbox{\rm div} \,
    \mbox{\boldmath $v$}(x) = \frac{d}{dx} \frac{1}{x} = - \frac{1}{x^2} < 0
    \quad \mbox{for} \; x \in (0,1).
  \]
  While the inf-sup condition \eqref{inf-sup advection-diffusion} remains
  valid for $v(x)=1/x$, we have to prove surjectivity in a different way:
  For $0 \neq p \in H^1_0(0,1)$ we define $u_p(x) := x \, p(x)$ with
  $u_p \in Y$. It remains to consider
  \[
    \| v \, u_p' \|_{Y^*} = \sup\limits_{0 \neq q \in Y}
    \frac{\langle v \, u_p' , q \rangle_{(0,1)}}{\| q \|_Y} \, ,
  \]
  where we have
  \begin{eqnarray*}
    \langle v \, u_p' , q \rangle_{(0,1)}
    & = & \int_0^1 \frac{1}{x} [p(x) + x p'(x)] q(x) \, dx \\
    & = & \int_0^1 \frac{1}{x} \, p(x) \, q(x) \, dx +
          \int_0^1 p'(x) q(x) \, dx \, .
  \end{eqnarray*}
  Note that the second term can be bounded as
  \begin{eqnarray*}
    \left| \int_0^1 p'(x) q(x) \, dx \right|
    & \leq & \left( \int_0^1 [p'(x)]^2 dx \right)^{1/2}
             \left( \int_0^1 [q(x)]^2 dx \right)^{1/2} \\
    && \hspace*{-3cm} = \, \| p \|_Y \left( \int_0^1
       \left[ \int_0^x q'(s) \, ds \right]^2 dx \right)^{1/2} \, \leq \,
       \| p \|_Y \left( \int_0^1 \int_0^x 1^2 ds \int_0^x [q'(s)]^2 ds \,
       dx \right)^{1/2} \\
    && \hspace*{-3cm} \leq \, \frac{1}{2} \, \| p \|_Y \| q \|_Y \, .
  \end{eqnarray*}
  In a similar way we have, recall $p,q \in H^1_0(0,1)$,
  \begin{eqnarray*}
    \left| \int_0^1 \frac{1}{x} \, p(x) \, q(x) \, dx \right|
    & = & \left| \int_0^1 \frac{1}{\sqrt{x}} \int_0^x p'(s) \, ds
          \, \frac{1}{\sqrt{x}} \int_0^x q'(s) \, ds \, dx \right| \\
    && \hspace*{-3cm} \leq
       \, \left( \int_0^1 \frac{1}{x} \left[ \int_0^x p'(s) \, ds \right]^2
             dx \right)^{1/2}
             \left( \int_0^1 \frac{1}{x} \left[ \int_0^x q'(s) \, ds \right]^2
             dx \right)^{1/2} \\
    & & \hspace*{-3cm}
        \leq \, \left( \int_0^1 \frac{1}{x} \int_0^x 1^2 ds \int_0^x [p'(s)]^2
             \, ds \, dx \right)^{1/2}
             \left( \int_0^1 \frac{1}{x} \int_0^x 1^2 ds \int_0^x [q'(s)]^2
        \, ds \, dx \right)^{1/2} \\[3mm]
    & & \hspace*{-3cm} \leq \, \| p \|_Y \| q \|_Y \, .
  \end{eqnarray*}
  Hence we conclude
  \[
    \| v \, u_p' \|_{Y^*} \leq \frac{3}{2} \, \| p \|_Y,
  \]
  i.e., $u_p \in X_v$. Then we obtain
  \begin{eqnarray*}
    b(u_p,p)
    & = & \int_0^1 u_p'(x) p'(x) \, dx + \int_0^1 \frac{1}{x} \, u_p'(x) \,
          p(x) \, dx \\
    & = & \int_0^1 [p(x) + x p'(x)] p'(x) \, dx +
          \int_0^1 \frac{1}{x} \, [p(x)+x p'(x)] p(x) \, dx \\
    & = & 2 \int_0^1 p(x) p'(x) \, dx +
          \int_0^1 x \, [p'(x)]^2 \, dx + \int_0^1 \frac{1}{x} \, [p(x)]^2 dx \\
    & = & \int_0^1 x \, [p'(x)]^2 \, dx + \int_0^1 \frac{1}{x} \, [p(x)]^2 dx
          \, > \, 0,
  \end{eqnarray*}
  i.e., surjectivity and therefore unique solvability of
  \eqref{DBVP Example 1D} follows also for $v(x)=1/x$.
\end{example}

\begin{remark}
  One possibility to prove the surjectivity result for the bilinear form
  $b(u,q)$ is to assume \eqref{Assumption v} to be able to conclude
  $X_v = Y$, which corresponds to the standard setting. As we have seen
  in the previous example, \eqref{Assumption v} is not necessary to
  establish surjectivity. At this time it remains open to generalize
  this result to more general situations. In any case, for the
  numerical realisation, surjectivity follows from injectivity, and
  assumption \eqref{Assumption v} is not needed.
\end{remark}

\subsection{Finite element discretization}
Let $X_H = Y_H = S_H^1(\Omega) \cap H^1_0(\Omega) = \mbox{span} \{
\varphi_k \}_{k=1}^M$ be a standard finite element space of, e.g.,
piecewise linear continuous basis functions $\varphi_k$ which are defined
with respect to some admissible decomposition of $\Omega$ into
shape regular simplicial finite elements $\tau_\ell$ of local mesh
size $H_\ell$, $\ell=1,\ldots,N$.
For any given $u \in X_v$, and similar as in \eqref{Definition wu},
we define $w_{u,H} \in Y_H$ as the unique solution
of the variational problem
\begin{equation}\label{Def wuh}
  \int_\Omega \alpha(x) \, \nabla w_{u,H}(x) \cdot \nabla q_H(x) \, dx =
  \int_\Omega \mbox{\boldmath $v$}(x) \cdot \nabla u(x) \, q_H(x) \, dx
  \quad \mbox{for all} \; q_H \in Y_H .
\end{equation}
Hence, we can introduce the discrete norm
\[
  \| u \|_{X_v,H} := \sqrt{\| u \|_Y^2 + \| w_{u,H} \|_Y^2} \leq
  \sqrt{\| u \|_Y^2 + \| w_u \|_Y^2} = \| u \|_{X_v} \quad
  \mbox{for all} \; u \in X_v .
\]
With respect to this discrete norm and as
in the continuous case, see Lemma \ref{Lemma solvability}, we
can prove a discrete inf-sup stability condition:

\begin{lemma}
  Assume 
  $\mbox{\rm div} \, \mbox{\boldmath $v$}(x) \leq 0$ for almost all
  $x \in \Omega$. Let $X_H \subset X_v$ and $Y_H \subset Y$ be conforming
  finite element spaces satisfying $X_H \subseteq Y_H$. Then there holds
  the discrete inf-sup stability condition
\begin{equation}\label{discrete inf-sup}
  \| u_H \|_{X_v,H} 
  \leq \sup\limits_{0 \neq q_H \in Y_H} \frac{b(u_H,q_H)}{\| q_H \|_Y}
  \quad \mbox{for all} \; u_H \in X_H .
\end{equation}

\end{lemma}
\proof{For $u_H \in X_H$ we define $w_{u_H,H} \in Y_H$ as the unique
  solution of the variational formulation
  \begin{equation}\label{Def wuhh}
    \int_\Omega \alpha(x) \, \nabla w_{u_H,H}(x) \cdot \nabla q_H(x) \, dx =
    \int_\Omega \mbox{\boldmath $v$}(x) \cdot \nabla u_H(x) \, q_H(x) \, dx
    \quad \mbox{for all} \; q_H \in Y_H .
  \end{equation}
  As in the continuous case we now define
  $q_{u_H,H} := u_H + w_{u_H,H} \subset Y_H$ to conclude
  \begin{eqnarray*}
    b(u_H,q_{u_H,H})
    & = & \int_\Omega \alpha(x) \,  \nabla u_H(x) \cdot
          \nabla q_{u_H,H}(x) \, dx + \int_\Omega
          \mbox{\boldmath $v$}(x) \cdot \nabla u_H(x) \, q_{u_H,H}(x) \, dx \\
    & = & \int_\Omega \alpha(x) \,  \nabla u_H(x) \cdot
          \nabla q_{u_H,H}(x) \, dx + \int_\Omega
          \alpha(x) \, \nabla w_{u_H,H}(x) \cdot \nabla q_{u_H,H}(x) \, dx \\
    & = & \int_\Omega \alpha(x) \, \nabla [u_H(x)+w_{u_H,H}(x)]
          \cdot \nabla q_{u_H,H} \, dx \\
    & = & \int_\Omega \alpha(x) \, |\nabla q_{u_H,H}(x)|^2 \, dx \, = \,
          \| q_{u_H,H} \|^2_Y \, ,
  \end{eqnarray*}
  and
  \begin{eqnarray*}
    \| q_{u_H,H} \|_Y^2
    & = & \| u_H + w_{u_H,H} \|_Y^2 \, = \,
          \int_\Omega \alpha(x) \, |\nabla_x[u_H(x)+w_{u_H,H}(x)]|^2 \, dx \\
    & = & \int_\Omega \alpha(x) \, |\nabla_x u_H(x)|^2 \, dx
          + \int_\Omega \alpha(x) \, |\nabla_x w_{u_H,H}(x)|^2 \, dx \\
    && \hspace*{4.5cm}
        + \, 2 \,
        \int_\Omega \alpha(x) \, \nabla_x w_{u_H,H}(x)
        \cdot \nabla u_H(x) \, dx \\
    & = & \| u_H \|_Y^2 + \| w_{u_H,H} \|_Y^2 + 2 \,
          \int_\Omega \mbox{\boldmath $v$}(x) \cdot
          \nabla u_H(x) \, u_H(x) \, dx \, \geq \, \| u_H \|_{X_v,H}^2,
  \end{eqnarray*}  
  i.e., the assertion follows.}

\noindent
Let $ u_H \in X_H$ be the unique Galerkin solution related to the variational
formulation \eqref{VF DBVP Advection Diffusion}, satisfying
\begin{equation}\label{FEM DBVP Advection Diffusion}
  \int_\Omega \alpha(x) \, \nabla u_H(x) \cdot \nabla q_H(x) \, dx +
  \int_\Omega \mbox{\boldmath $v$}(x) \cdot \nabla u_H(x) \,
  q_H(x) \, dx = \int_\Omega f(x) \, q_H(x) \, dx
\end{equation}
for all $q_H \in Y_H$. The Galerkin variational formulation
\eqref{FEM DBVP Advection Diffusion}
is equivalent to an algebraic system of linear equations,
$B_H \underline{u} = \underline{f}$, where the stiffness matrix is
given by its entries
\[
  B_H[j,k] = \int_\Omega \alpha(x) \, \nabla \varphi_k(x) \cdot
  \nabla \varphi_j(x) \, dx + \int_\Omega \mbox{\boldmath $v$}(x) \cdot
  \nabla \varphi_k(x) \, \varphi_j(x) \, dx
\]
for $k,j=1,\ldots,M$. Invertibility of $B_H$ follows from the discrete
inf-sup condition \eqref{discrete inf-sup}.
In fact, for any given $z \in X_v$ we can define the Galerkin projection
$z_H = G_H z \in X_H$ as the unique solution of the variational formulation
\[
b(G_H z , q_H) = b(z,q_H) \quad \mbox{for all} \; q_H \in Y_H ,
\]
and when combining the discrete inf-sup stability condition
\eqref{discrete inf-sup} and the boundedness result \eqref{boundedness b}
of the bilinear form $b(\cdot, \cdot)$ we immediately conclude
\begin{equation}\label{boundedness Gh}
  \| G_H z \|_{X_v,H} \leq \sqrt{2} \, \| z \|_{X_v} \quad \mbox{for all} \;
  z \in X_v .
\end{equation}
With \eqref{boundedness Gh} and using
\cite{Carstensen:2025, Zikatanov:2003} we then can formulate Cea's lemma,
\begin{equation}\label{Cea}
  \| u - u_H \|_{X_v,H} \leq \sqrt{2} \,
  \inf\limits_{z_H \in X_H} \| u - z_H \|_{X_v} .
\end{equation}

\begin{lemma}\label{Lemma Error Advection Diffusion}
  Let $u_H \in X_H$ be the unique solution of the Galerkin variational
  formulation \eqref{FEM DBVP Advection Diffusion}, where we assume
  $u \in H^2(\Omega)$. When assuming \eqref{Assumption v} and
  $\mbox{\rm div} \, \mbox{\boldmath $v$}(x) \leq 0 $ for almost
  all $x \in \Omega$ there holds the error estimate
  \begin{equation}\label{Advection Diffusion Error 1}
    \| u - u_H \|_Y^2 \leq c \sum\limits_{\ell=1}^N
    H_\ell^2 \, \left[ \sup\limits_{x \in \tau_\ell} \, \alpha(x) +
      \frac{c_v^2}{\underline{\alpha}}
    \right] \, |u|_{H^2(\tau_\ell)}^2 .
  \end{equation}
  Moreover, when
  $\mbox{\boldmath $v$} \in [L^\infty(\Omega)]^n$ and
  $\mbox{\rm div} \, \mbox{\boldmath $v$}(x) = 0$ for almost all
  $x \in \Omega$ is satisfied, the error estimate
  \begin{equation}\label{Advection Diffusion Error 2}
    \| u - u_H \|_Y^2 \leq c \sum\limits_{\ell=1}^N
    \left[
      H_\ell^2 \, \sup\limits_{x \in \tau_\ell} \, \alpha(x) +
      H_\ell^4 \, \sup\limits_{x \in \tau_\ell}
      \frac{|\mbox{\boldmath $v$}(x)|^2}{\alpha(x)}
    \right] \, |u|_{H^2(\tau_\ell)}^2 
  \end{equation}
  follows.
\end{lemma}
\proof{From \eqref{Cea} we first have
  \[
    \| u - u_H \|_Y^2 \, \leq \, \| u - u_H \|_{X_v,H}^2 \leq
    \| u - I_H u \|_{X_v}^2 \, = \, \| u - I_H u \|_Y^2 +  \| \mbox{\boldmath $v$} \cdot \nabla (u-I_Hu) \|_{Y^*},
  \]
  where $I_Hu$ is the piecewise linear nodal interpolation of the solution
  $u \in H^2(\Omega)$.
  Using standard local interpolation error estimates we obtain
  \begin{eqnarray*}
    \| u - I_H u \|_Y^2
    & = & \int_\Omega \alpha(x) \, |\nabla [u(x)-I_Hu(x)]|^2 \, dx \\
    & = & \sum\limits_{\ell=1}^N \int_{\tau_\ell}
          \alpha(x) \, |\nabla [u(x)-I_Hu(x)]|^2 \, dx \, \leq \,
          c \sum\limits_{\ell=1}^N \sup\limits_{x \in \tau_\ell}
          \, \alpha(x) \, H_\ell^2 \, |u|^2_{H^2(\tau_\ell)} .
  \end{eqnarray*}
  For the remaining part, and using the dual norm \eqref{dual norm},
  we have to consider
  \begin{eqnarray*}
    \| \mbox{\boldmath $v$} \cdot \nabla (u-I_Hu) \|_{Y^*}
    & = & \sup\limits_{0 \neq q \in Y}
          \frac{\langle \mbox{\boldmath $v$} \cdot \nabla (u-I_Hu),q
          \rangle_\Omega} {\| q \|_Y} = \sup\limits_{0 \neq q \in Y}
          \frac{\langle \nabla (u-I_Hu), q \, \mbox{\boldmath $v$}
          \rangle_\Omega} {\| q \|_Y} \\
    & \leq & \| \nabla (u-I_Hu) \|_{L^2(\Omega)}
             \sup\limits_{0 \neq q \in Y}
             \frac{\| q \, \mbox{\boldmath $v$} \|_{L^2(\Omega)}}{\| q \|_Y} \\
    & \leq & \frac{1}{\sqrt{\underline{\alpha}}} \,
             \| \nabla (u-I_Hu) \|_{L^2(\Omega)}
             \sup\limits_{0 \neq q \in Y}
             \frac{\| q \, \mbox{\boldmath $v$} \|_{L^2(\Omega)}}
             {\| \nabla q \|_{L^2(\Omega)}} \\
    & = & \frac{c_v}{\sqrt{\underline{\alpha}}} \,
          \| \nabla (u-I_Hu) \|_{L^2(\Omega)} \, ,
  \end{eqnarray*}
  where we have used \eqref{Assumption v}, i.e., we obtain
  \[
    \| \mbox{\boldmath $v$} \cdot \nabla (u-I_Hu) \|_{Y^*}^2 \leq
    c \, \frac{c_v^2}{\underline{\alpha}} \,
    \sum\limits_{\ell=1}^N H_\ell^2 \, |u|^2_{H^2(\tau_\ell)} \, ,
  \]
  and \eqref{Advection Diffusion Error 1} follows.

  In the case $\mbox{div} \, \mbox{\boldmath $v$}(x) = 0$ for almost
  all $x \in \Omega$ and $\mbox{\boldmath $v$} \in [L^\infty(\Omega)]^n$
  we can use integration by parts to conclude
  \begin{eqnarray*}
    \left|
    \langle \mbox{\boldmath $v$} \cdot \nabla (u-I_Hu),q \rangle_\Omega
    \right|
    & = & \left| \int_\Omega \mbox{\boldmath $v$}(x) \cdot
          \nabla [u(x)-I_Hu(x)] \, q(x) \, dx \right| \\
    && \hspace*{-2cm}
       = \, \left| - \int_\Omega [u(x)-I_Hu(x)] \,
       \mbox{div} [ q(x) \, \mbox{\boldmath $v$}(x)] \, dx \right| \\
    && \hspace*{-2cm}
       = \, \left| - \int_\Omega [u(x)-I_Hu(x)] \,
       \mbox{\boldmath $v$}(x) \cdot \nabla q(x) \, dx \right| \\
    && \hspace*{-2cm} \leq \,  \sum\limits_{\ell=1}^N \left|
       \int_{\tau_\ell} [u(x)-I_Hu(x)] \, \mbox{\boldmath $v$}(x) \cdot
       \nabla q(x) \, dx \right| \\
    && \hspace*{-2cm}  \leq \,
       \sum\limits_{\ell=1}^N \| u - I_H u \|_{L^2(\tau_\ell)}
       \left( \int_{\tau_\ell} [
       \mbox{\boldmath $v$}(x) \cdot \nabla q(x)]^2 \, dx \right)^{1/2} \\
    && \hspace*{-2cm} \leq \, c \sum\limits_{\ell=1}^N
       H_\ell^2 \, | u |_{H^2(\tau_\ell)} \sup\limits_{x \in \tau_\ell}
       \frac{|\mbox{\boldmath $v$}(x)|}{\sqrt{\alpha(x)}}
       \left( \int_{\tau_\ell} \alpha(x) \,
       |\nabla q(x)|^2 \, dx \right)^{1/2} \\
    && \hspace*{-2cm} \leq \, c \left(
       \sum\limits_{\ell=1}^N H_\ell^4 \, \sup\limits_{x \in \tau_\ell}
       \frac{|\mbox{\boldmath $v$}(x)|^2}{\alpha(x)} \, |u|^2_{H^2(\tau_\ell)}
       \right)^{1/2} \left(
       \sum\limits_{\ell=1}^N \int_{\tau_\ell} \alpha(x) \, |\nabla q(x)|^2 dx
       \right)^{1/2} \\
    && \hspace*{-2cm} = \, c \left(
       \sum\limits_{\ell=1}^N H_\ell^4 \, \sup\limits_{x \in \tau_\ell}
       \frac{|\mbox{\boldmath $v$}(x)|^2}{\alpha(x)} \,
       |u|^2_{H^2(\tau_\ell)} \right)^{1/2} \, \| q \|_Y \, .
  \end{eqnarray*}
  Hence we conclude
  \[
    \| \mbox{\boldmath $v$} \cdot \nabla (u-I_Hu) \|_{Y^*}^2
    \, \leq \, c^2
    \sum\limits_{\ell=1}^N H_\ell^4 \, \sup\limits_{x \in \tau_\ell}
    \frac{|\mbox{\boldmath $v$}(x)|^2}{\alpha(x)} \,
    |u|^2_{H^2(\tau_\ell)},
  \]
   and \eqref{Advection Diffusion Error 2} follows.}

\begin{example}
  As a first example, we consider the convection-diffusion equation
  \eqref{DBVP Advection Diffusion} in the particular case
  $\alpha(x) = 1$ and $\mbox{\boldmath $v$}(x)=v$
  for all $x \in \Omega = (0,1)$, i.e.,
  $\underline{\alpha}=\overline{\alpha}=1$.
  For a globally quasi-uniform mesh with
  $H_\ell \sim H$ for all $\ell=1,\ldots,N$ we can write the finite element
  error estimate \eqref{Advection Diffusion Error 2} as
  \begin{equation}\label{Error Example 1}
    \| \nabla (u-u_H) \|_{L^2(\Omega)}
    \leq c \, \Big( H^2 + v \, H^4 \Big)^{1/2} \, |u|_{H^2(\Omega)} ,
  \end{equation}
  when assuming $u \in H^2(\Omega)$. As in
  {\rm \cite[Example 1.2]{RoosStynesTobiska:2008}} we consider the solution
  \begin{equation}\label{eq:smooth u ex1}
    u(x) = x - \frac{\exp\left(x-1\right) - \exp\left(-1\right)}
    {1 - \exp\left(-1\right)} \quad \mbox{for} \; x\in (0,1)
  \end{equation}
  and compute the right hand side accordingly. For the discretization we
  use piecewise linear basis functions which are defined with respect to
  some uniform decompostion of $\Omega = (0,1)$ into $N$ finite elements
  of mesh size $H=1/N$. The numerical results are given in Table
  \ref{tab:Erros for advdiff1d with different velocities uniform} which
  confirm linear convergence asymptotically, as expected. On the other hand,
  for $v=1000$, we observe some initial higher order convergence
  as indicated in the error estimate \eqref{Error Example 1},
  see also the discussion on related observations in
  {\rm \cite{Pacheco:2022,Steinbach:2000}}.

  \begin{table}[h]
    \begin{center}
        \begin{adjustbox}{width=\textwidth}
          \input{tables/SimData_advdiff1d_direct_DifferentVelocities.tex}
        \end{adjustbox}
    \end{center}
    \caption{Convection-diffusion equation \eqref{DBVP Advection Diffusion} in
      $\Omega = (0,1)$, $\alpha \equiv 1$,
      $\mbox{\boldmath $v$}=v \in \{ 10, 100, 1000 \}$,
      error and estimated order of
      convergence (eoc) for the solution $u$ as given in
      \eqref{eq:smooth u ex1}.}
    \label{tab:Erros for advdiff1d with different velocities uniform}
  \end{table}

\end{example}

\begin{example}\label{ex:velocity as a function of x}
  In this example, we consider the convection-diffusion equation
  \eqref{DBVP Advection Diffusion} for $\Omega =(0,1)$, $\alpha(x) = 1$
  for $x \in (0,1)$, and a velocity field which depends on the spatial 
  variable, i.e, $\mbox{\boldmath $v$}(x)=v(x)$. As exact solution we
  consider \eqref{eq:smooth u ex1} and compute the right hand side
  accordingly to the related chosen velocity fields. Firstly, we choose  
  $v(x) = 1/\sqrt{x}$ which satisfies \eqref{Assumption v}  and
  $\divv v(x) \leq 0$ for $x \in (0,1)$. Hence, for a globally quasi-uniform 
  mesh the error estimate \eqref{Advection Diffusion Error 1} gives 
  \begin{equation*}
    \|\nabla (u - u_H )\|_{L^2(\Omega)}\leq c H \, | u |_{H^2(\Omega)}
  \end{equation*}
  when assuming $u \in H^2(\Omega)$. Secondly, we choose $v(x) = 1/x$
  with $\divv v(x) \leq 0$ for $x \in (0,1)$ but we can not use
  \eqref{Bound cv}. Although we could prove surjectivity in a different way,
  we can not apply the error estimates as given in
  Lemma \ref{Lemma Error Advection Diffusion}.
  The numerical results for both choices of the velocity field can be seen
  in Tab. \ref{tab:Erros for advdiff1d with v(x) on uniform grid}. 
  In both cases we observe a linear convergence rate which confirms the
  a priori error estimate \eqref{Advection Diffusion Error 1}.
  
  \begin{table}[h]
    \begin{center}
          \input{tables/SimData_advdiff1d_direct_VelocityAsFunctionInx.tex}
    \end{center}
    \caption{Convection-diffusion equation \eqref{DBVP Advection Diffusion} in
      $\Omega = (0,1)$, $\alpha \equiv 1$,
      $\mbox{\boldmath $v$}(x)=v(x) \in \left\{ \frac{1}{\sqrt{x}},
        \frac{1}{x} \right\}$,
      error and estimated order of
      convergence (eoc) for the solution $u$ as given in
      \eqref{eq:smooth u ex1}.}
    \label{tab:Erros for advdiff1d with v(x) on uniform grid}
  \end{table}

\end{example}

\begin{example}\label{ex:2}
  We consider $\alpha(x) = \varepsilon < \! \! < 1$,
  i.e., a singularly perturbed problem with dominating convection.
  Using again a globally quasi-uniform mesh, the error estimate
  \eqref{Advection Diffusion Error 2} now gives
  \[
  \varepsilon \, \| \nabla (u - u_H) \|^2_{L^2(\Omega)} \leq c \, \Big[
  H^2 \, \varepsilon + \| \mbox{\boldmath $v$} \|^2_{[L^\infty(\Omega)]^n} \,
  H^4 \, \varepsilon^{-1} \, \Big] \, |u|^2_{H^2(\Omega)},
  \]
  when assuming $u \in H^2(\Omega)$, i.e.,
  \[
    \| \nabla (u - u_H) \|_{L^2(\Omega)} \leq c \,
    \Big( H^2 + H^4 \, \varepsilon^{-2} \Big)^{1/2} \, |u|_{H^2(\Omega)} \, .
  \]
  As in {\rm \cite[Example 1.2]{RoosStynesTobiska:2008}} we
  consider $\Omega=(0,1)$
  and $\mbox{\boldmath $v$}(x) = v = 1$ with the exact solution to be
  \begin{equation}\label{eq:uexakt1d bndlayer}
    u(x) = x - \frac{\exp\left(-\frac{1-x}{\varepsilon}\right) -
      \exp\left(-\frac{1}{\varepsilon}\right)}
    {1 - \exp\left(-\frac{1}{\varepsilon}\right)},\quad x\in (0,1),
  \end{equation}
  which has a boundary layer at $x=1$. The convergence behaviour of the 
  numerical solutions $u_H$ for
  $\varepsilon \in \left\{10^{-2},\,10^{-4},\,10^{-5}\right\}$ 
  are depicted in Tab. \ref{tab:Erros for advdif1d uniform} and 
  Fig. \ref{abb:ex_bndlayer1d_error_plot}. We see that we obtain linear 
  convergence if the mesh size $H$ is of the order of the singular 
  perturbation parameter $\varepsilon$, i.e., $H \sim \varepsilon$.
  Before, the numerical solutions $u_H$ obtain oscillations due to an
  insufficient resolution of the boundary layer, see also
  Fig. \ref{abb:Oscillations in solution for convection diffusion}.
  However, when using a globally quasi-uniform mesh the choice 
  $H \sim \varepsilon$ is impractical in applications as 
  this leads to an unacceptably large number of grid points. This motivates 
  to consider an adaptive approach, which works for all values of the 
  singular perturbation parameter $\varepsilon$.
  
  \begin{table}
    \begin{center}
      \begin{adjustbox}{width=\textwidth}
        \input{tables/SimData_advdiff1d_direct.tex}
      \end{adjustbox}
    \end{center}
    \caption{Convection-diffusion equation \eqref{DBVP Advection Diffusion}
      in $\Omega = (0,1)$, $\alpha(x) = \varepsilon \in \{ 10^{-2},
      10^{-4}, 10^{-5} \}$, $\mbox{\boldmath $v$}=1$,
      error and estimated order of convergence (eoc) for the solution
      $u$ as given in \eqref{eq:uexakt1d bndlayer}.}
    \label{tab:Erros for advdif1d uniform}
  \end{table}

  \begin{figure}
    \centering
    \begin{tikzpicture}[scale = 0.8]
        \begin{axis}[
            xmode = log,
            ymode = log,
            xlabel= $N$,
            ylabel=errors,
            legend style={font=\tiny}, legend pos = south west]
            \addplot[mark = *,blue] table [col sep=&, y=errH1_u, x=nv]
            {tables/SimData_laplace_advdiff1d_uniform_eps1e-02.dat};
            \addlegendentry{$\varepsilon = 10^{-2}$}
        
            \addplot[mark = *,red] table [col sep=&, y=errH1_u, x=nv]
            {tables/SimData_laplace_advdiff1d_uniform_eps1e-04.dat};
            \addlegendentry{$\varepsilon = 10^{-4}$}
              
            \addplot[mark = *,green!60!black] table [col sep=&, y=errH1_u, x=nv]
            {tables/SimData_laplace_advdiff1d_uniform_eps1e-05.dat};
            \addlegendentry{$\varepsilon = 10^{-5}$}
              
            \addplot[
            domain = 100000:1000000,
            samples = 10,
            thin,
            black,
          ] {100*x^(-1)};
          \addlegendentry{$h=N^{-1}$}
        \end{axis}
    \end{tikzpicture}
    \caption{Convergence behaviour of the error
      $\|\nabla (u - u_H) \|_{L^2(0,1)}$ for a uniform
      refinement strategy in case of the function $u$
      as given in \eqref{eq:uexakt1d bndlayer}.}
      \label{abb:ex_bndlayer1d_error_plot}
    \end{figure}
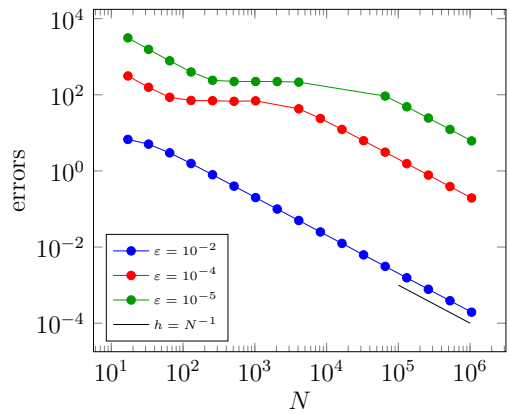
    
    \begin{figure}
      \begin{subfigure}[b]{0.32\textwidth}
          \centering
          \includegraphics[width =
          0.99\textwidth]{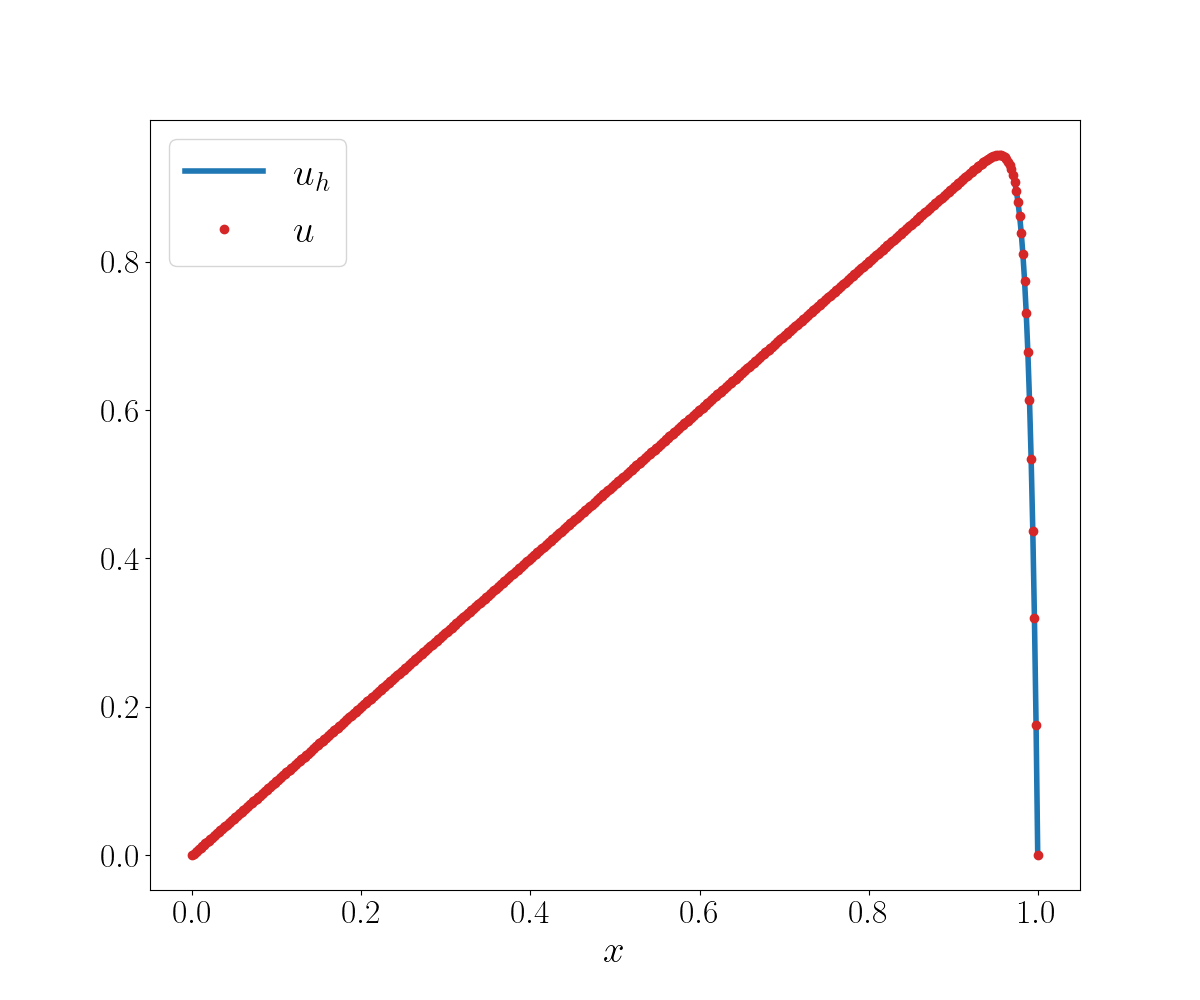} 
          \caption{$\varepsilon = 10^{-2}$}
      \end{subfigure}
      ~
      \begin{subfigure}[b]{0.32\textwidth}
          \centering
          \includegraphics[width =
          0.99\textwidth]{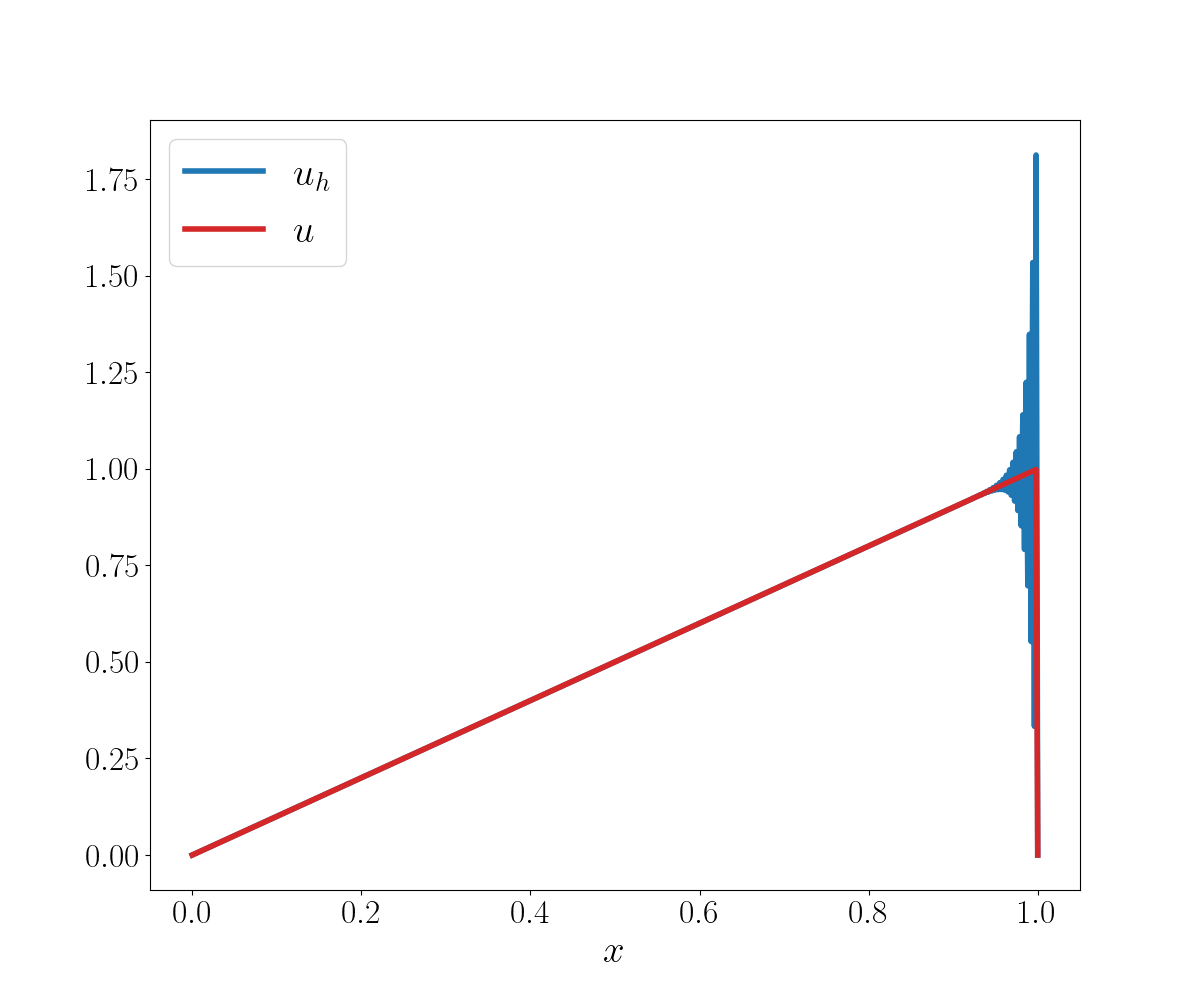}
          \caption{$\varepsilon = 10^{-4}$} 
      \end{subfigure}
      ~
      \begin{subfigure}[b]{0.32\textwidth}
          \centering
          \includegraphics[width =
          0.99\textwidth]{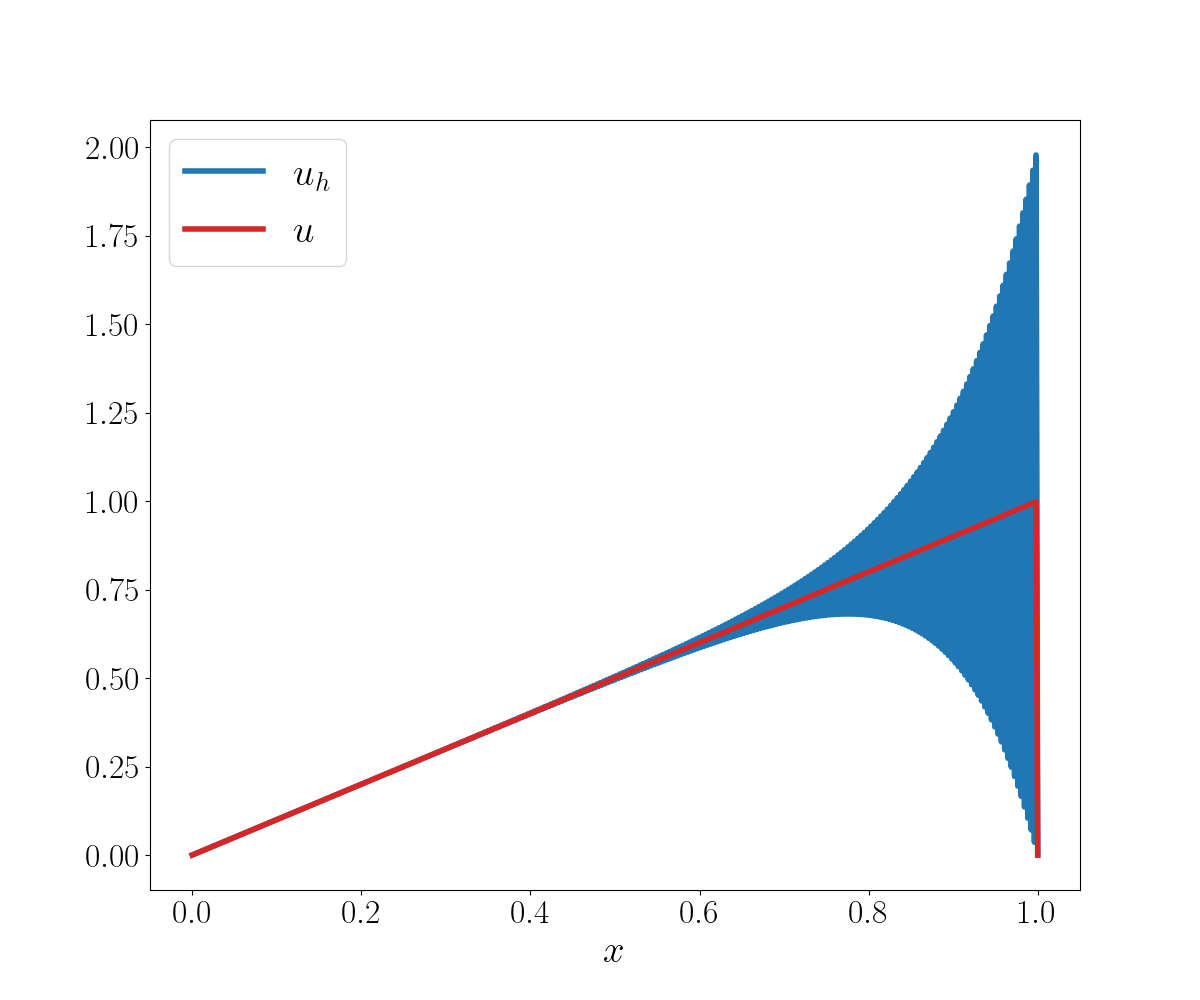}
          \caption{$\varepsilon = 10^{-5}$} 
      \end{subfigure}
      \caption{Numerical solutions $u_H$ to Example \ref{ex:2}
        in case of $N = 512$ elements.}
      \label{abb:Oscillations in solution for convection diffusion}
    \end{figure}
\end{example}

\subsection{Adaptive least squares finite element method}
The bilinear forms $a(\cdot,\cdot)$ and $b(\cdot,\cdot)$ as defined in
\eqref{Definition wu} and \eqref{VF DBVP Advection Diffusion} imply
bijective operators $A : Y \to Y^*$ and $B : X_v \to Y^*$ satisfying
\[
  \langle A p , q \rangle_\Omega = \int_\Omega \alpha(x) \,
  \nabla p(x) \cdot \nabla q(x) \, dx \quad \mbox{for all} \;
  p,q \in Y,
\]
and
\[
  \langle B u , q \rangle_\Omega = b(u,q) \quad \mbox{for all} \;
  (u,q) \in X_v \times Y .
\]
Hence, instead of the operator equation $B u = f$ in $Y^*$, i.e., of the
variational formulation \eqref{VF DBVP Advection Diffusion}, we now consider
the equivalent problem to minimize the quadratic functional
\begin{equation}\label{minimization problem}
  {\mathcal{J}}(z) = \frac{1}{2} \, \| B z - f \|_{A^{-1}}^2 \quad
  \mbox{for} \; z \in X_v .
\end{equation}
The minimizer $u \in X_v$ of \eqref{minimization problem} is given as the
unique solution of the gradient equation
\begin{equation}\label{gradient equation}
B^* A^{-1} (B u - f) \, = \, 0 ,
\end{equation}
and using the adjoint $p = A^{-1} (f-Bu) \in Y $
we have to solve a coupled system
to find $(u,p) \in X_v \times Y$ such that
\[
A p + B u = f, \quad B^* p = 0\, ,
\]
i.e.,
\begin{eqnarray}
  \label{LS1}
  \int_\Omega \alpha(x) \nabla p(x) \cdot \nabla q(x) dx +
  \int_\Omega \alpha(x) \nabla u(x) \cdot \nabla q(x) dx +
  \int_\Omega \mbox{\boldmath $v$}(x) \cdot \nabla u(x) \, q(x) dx
  \hspace*{1.5cm} && \\
  =
  \int_\Omega f(x) q(x) dx && \nonumber
\end{eqnarray}
for all $q \in Y$, and
\begin{equation}\label{LS2}
  \int_\Omega \alpha(x) \nabla z(x) \cdot \nabla p(x) dx +
  \int_\Omega \mbox{\boldmath $v$}(x) \cdot \nabla z(x) \, p(x) dx
  = 0 \quad \mbox{for all} \; z \in X_v .
\end{equation}
Unique solvability of the gradient equation \eqref{gradient equation}
and therefore of the equivalent mixed variational formulation
\eqref{LS1} and \eqref{LS2} follows as in
\cite{KLS:2023}, since $S := B^* A^{-1} B : X_v \to X_v^*$ is bounded and
elliptic. Note that in the continuous case we have $p \equiv 0$, since
$B$ and therefore $B^*$ are bijective.

For the Galerkin finite element discretization of the mixed variational
system \eqref{LS1} and \eqref{LS2} we will use again the
conforming finite element space
$X_H = \mbox{span} \{ \varphi_k \}_{k=1}^M \subset X_v$, but we introduce
a finite element space
$Y_h = \mbox{span} \{ \phi_i \}_{i=1}^{\overline{M}} \subset Y$
of piecewise linear continuous basis functions $\phi_i$ which are defined
with respect to a possibly refined decomposition of $\Omega$ into $\overline{N}$
finite elements $\overline{\tau}_\ell$ of local mesh size $h_\ell$, e.g.,
starting from $X_H$ we use at least one additional refinement to
construct $Y_h$,
i.e., $h_\ell = H_\ell/2$ when $\overline{\tau}_\ell \subset \tau_\ell$.
This definition ensures $X_H \subset Y_h$, but even the choice $Y_h = X_H$
would be allowed for a stable discretization.
In both cases, \eqref{discrete inf-sup}
implies the discrete inf-sup condition
\begin{equation}\label{discrete inf-sup 2}
\| u_H \|_{X_v,H}  
  \leq \sup\limits_{0 \neq q_h \in Y_h} \frac{b(u_H,q_h)}{\| q_h \|_Y}
  \quad \mbox{for all} \; u_H \in X_H . 
\end{equation}
The mixed finite element
discretization of \eqref{LS1} and \eqref{LS2} is to find
$(p_h,u_H) \in Y_h \times X_H$ such that
\begin{eqnarray}
  \int_\Omega \alpha(x) \nabla p_h(x) \cdot \nabla q_h(x) dx +
  \int_\Omega \alpha(x) \nabla u_H(x) \cdot \nabla q_h(x) dx +
  \int_\Omega \mbox{\boldmath $v$}(x) \cdot \nabla u_H(x) \, q_h(x) dx
  && \nonumber \\
  =
  \int_\Omega f(x) q_h(x) dx \hspace*{2cm} && \label{LS1 FEM}
\end{eqnarray}
for all $q_h \in Y_h$, and
\begin{equation}\label{LS2 FEM}
  \int_\Omega \alpha(x) \nabla z_H(x) \cdot \nabla p_h(x) dx +
  \int_\Omega \mbox{\boldmath $v$}(x) \cdot \nabla z_H(x) \, p_h(x) dx
  = 0 \quad \mbox{for all} \; z_H \in X_H .
\end{equation}
This is equivalent to a coupled system of linear algebraic equations
\begin{equation}\label{LS LGS}
  \left(
    \begin{array}{cc}
      A_h & B_h \\[1mm]
      B_h^\top &
    \end{array}
  \right)
  \left(
    \begin{array}{c}
      \underline{p} \\[1mm]
      \underline{u}
    \end{array}
  \right)
  =
  \left(
    \begin{array}{c}
      \underline{f} \\[1mm]
      \underline{0}
    \end{array}
  \right)
\end{equation}
where the entries of the stiffness matrix are now given by
\begin{eqnarray*}
  A_h[j,i]
  & = & \int_\Omega \alpha(x) \, \nabla \phi_i(x) \cdot
        \nabla \phi_j(x) \, dx, \\
  B_h[j,k]
  & = & \int_\Omega \alpha(x) \, \nabla \varphi_k(x) \cdot
        \nabla \phi_j(x) \, dx + \int_\Omega \mbox{\boldmath $v$}(x) \cdot
        \nabla \varphi_k(x) \, \phi_j(x) \, dx
\end{eqnarray*}
for $i,j=1,\ldots,\overline{M}$, $k=1,\ldots,M$. Unique solvability
of \eqref{LS LGS} follows from \eqref{discrete inf-sup 2}. Since
$A_h$ is invertibe, we can also consider the Schur complement system
\[
B_h^\top A_h^{-1} B_h \underline{u} = B_h^\top A_h^{-1} \underline{f} 
\]
which is the discrete counter part of the gradient equation
\eqref{gradient equation}. In the particular case $Y_h = X_H$ we have
$B_h = B_H$ to be invertible, and $\underline{p}=\underline{0}$
follows in this case. Hence, to have a built in error estimation we
now consider the case $X_H \subset Y_h$, $X_H \neq Y_h$ only.
As in \eqref{Def wuh} we now define $w_{u,h} \in Y_h$, and introduce
the discrete norm $\| u \|_{X_v,h} = \sqrt{\| u \|_Y^2+\| w_{u,h} \|_Y^2}$.
As in \cite[Lemma 2.4, equations (2.16), (2.31)]{KLS:2023} we then
conclude the a priori error estimate
\begin{equation}\label{LS Error}
\| u - u_H \|_{X_v,h} \leq c \, \inf\limits_{z_H \in X_H} \| u - z_H \|_{X_v},
\end{equation}
and we finally obtain the error estimates as already given in
Lemma \ref{Lemma Error Advection Diffusion}. In order to derive an
a posteriori error indicator we first have the following result:

\begin{lemma}
  Let $(p_h,u_H) \in Y_h \times X_H$ be the unique solution of the
  variational formulations \eqref{LS1 FEM} and \eqref{LS2 FEM}. Then there
  holds
  \begin{equation}
    \frac{1}{\sqrt{2}} \, \| p_h \|_Y \leq \| u - u_H \|_{X_v,h} \, .
  \end{equation}
\end{lemma}
\proof{When subtracting \eqref{LS1 FEM} from \eqref{LS1} for the test function
  $q_h = p_h \in Y_h$ this gives the Galerkin orthogonality,
  recall $p \equiv 0$,
  \begin{eqnarray*}
    \| p_h \|_Y^2
    & = & \int_\Omega \alpha(x) \, \nabla p_h(x) \cdot \nabla p_h(x) \, dx \\
    & = & \int_\Omega \alpha(x) \, \nabla [u(x)-u_H(x)] \cdot
          \nabla p_h(x) \, dx +
          \int_\Omega \mbox{\boldmath $v$}(x) \cdot \nabla [u(x)-u_H(x)] \,
          p_h(x) \, dx \\
    & = & \int_\Omega \alpha(x) \, \nabla [u(x)-u_H(x)] \cdot
          \nabla p_h(x) \, dx +
          \int_\Omega \nabla w_{u-u_H,h}(x) \cdot \nabla p_h(x) \, dx \\
    & \leq & \Big[ \| u - u_H \|_Y + \| w_{u-u_H,h} \|_Y \Big] \,
             \| p_h \|_Y \, \leq \,
             \sqrt{2} \, \Big( \| u - u_H \|_Y^2 +
             \| w_{u-u_H,h} \|^2_Y \Big)^{1/2} \,
             \| p_h \|_Y \\
    & = & \sqrt{2} \, \| u - u_H \|_{X_v,h} \| p_h \|_Y,
  \end{eqnarray*}
  i.e., the assertion follows.}

\noindent
It remains to prove that the error estimator $\| p_h \|_Y$ is also
reliable. For this we consider the variational formulations \eqref{LS1 FEM} and
\eqref{LS2 FEM}, using
$X_h=Y_h$, to find $(\overline{p}_h,u_h) \in Y_h \times X_h$ such that
\begin{equation}\label{Def uh}
  \langle A \overline{p}_h , q_h \rangle_\Omega +
  \langle B u_h , q_h \rangle_\Omega = \langle f , q_h \rangle_\Omega,
  \quad \langle B z_h , \overline{p}_h \rangle_\Omega = 0
\end{equation}
is satisfied for all $(z_h,q_h) \in Y_h \times X_h$. As in
\eqref{discrete inf-sup} we have the discrete inf-sup condition
\begin{equation}\label{discrete inf-sup h}
  \| u_h \|_{X_v,h} \leq \sup\limits_{0 \neq q_h \in Y_h}
  \frac{b(u_h,q_h)}{\| q_h \|_Y} \quad \mbox{for all} \; u_h \in X_h,
\end{equation}
from which unique solvability of \eqref{Def uh} follows.

\begin{lemma}
  Let $(p_h,u_H) \in Y_h \times X_H$ and
  $(\overline{p}_h,u_h) \in X_h \times Y_h$ be the unique solutions of the
  variational formulations \eqref{LS1 FEM}--\eqref{LS2 FEM} and
  \eqref{Def uh}, respectively. Assume the saturation condition
  \begin{equation}\label{saturation}
    \| u - u_h \|_{X_v,h} \leq \eta \, \| u-u_H \|_{X_v,h} \quad
    \mbox{for some} \; \eta \in (0,1) .
  \end{equation}
  Then there holds the reliability estimate
  \begin{equation}\label{reliability}
    \| u - u_H \|_{X_v,h} \leq \frac{1}{1-\eta} \, \| p_h \|_Y \, .
  \end{equation}
\end{lemma}
\proof{Since the variational formulation \eqref{Def uh} is considered
  for $X_h=Y_h$, $\overline{p}_h \equiv 0$ follows. Hence we obtain,
  using \eqref{LS1 FEM},
  \[
    \langle B u_h , q_h \rangle_\Omega =
    \langle f , q_h \rangle_\Omega =
    \langle A p_h , q_h \rangle_\Omega + \langle B u_H , q_h \rangle_\Omega
    \quad \mbox{for all} \; q_h \in Y_h ,
  \]
  i.e.,
  \[
    \langle B(u_h - u_H) , q_h \rangle_\Omega  =
    \langle A p_h , q_h \rangle_\Omega \quad
    \mbox{for all} \; q_h \in Y_h .
  \]
  Now, using the discrete inf-sup condition \eqref{discrete inf-sup h}
  for $u_h - u_H \in X_h$ this gives
  \[
    \| u_h - u_H \|_{X_v,h} \leq
    \sup\limits_{0 \neq q_h \in Y_h} \frac{\langle B(u_h-u_H),q_h\rangle_\Omega}
    {\| q_h \|_Y} = \sup\limits_{0 \neq q_h \in Y_h}
    \frac{\langle A p_h , q_h \rangle_\Omega}{\| q_h \|_Y} \leq
    \| p_h \|_Y \, .
  \]
  Then, using the triangle inequality and the saturation assumption
  \eqref{saturation}, we obtain
  \[
    \| u - u_H \|_{X_v,h}
    \leq \| u - u_h \|_{X_v,h} + \| u_h - u_H \|_{X_v,h} 
    \leq \eta \, \| u - u_H \|_{X_v,h} + \| p_h \|_Y \, ,
  \]
  i.e., the assertion follows.}
  
\begin{remark}
  The numerical results indicate that the saturation condition
  \eqref{saturation} is satisfied for the choice $h = H/2$, i.e.,
  one additional refinement to define $Y_h$ when starting from $X_H$.
  Otherwise, one may use some more additional refinement steps when
  required.
\end{remark}

\begin{example}\label{Example 2.3}
  In this example we apply the least-squares approach
  \eqref{minimization problem} to the singularly perturbed problem already
  considered in Example \ref{ex:2}. We use the finite element spaces
  $X_H = S_H^1(\mathcal{T}_H) \cap X_v$ and $Y_h = S_H^2(\mathcal{T}_H) \cap Y$,
  which are defined with respect to an admissible and locally
  quasi-uniform decomposition $\mathcal{T}_H = \{ \tau_\ell \}_{\ell=1}^N$
  of the interval $\Omega = (0,1)$ into finite elements $\tau_\ell$ with
  $\widetilde{M} = N+1$ nodes. We use the global error estimator 
  \begin{equation}\label{global error estimator}
    \eta_H^2 = \|p_h\|_Y^2 =
    \varepsilon \int_0^1 \nabla p_h(x) \cdot \nabla p_h(x) \, dx 
    = \sum_{\ell=1}^N \eta_\ell^2
  \end{equation}
  with the local error indicators 
  \begin{equation}\label{local error indicator}
    \eta_\ell^2 = \varepsilon \int_{\tau_\ell}
    \nabla p_h(x) \cdot \nabla p_h(x) \, dx
  \end{equation}
  to drive an adaptive refinement scheme with a Dörfler marking strategy
  \cite{Doerfler:1996}. The sparse direct solver \texttt{Pardiso} is used
  to solve the resulting linear systems. The convergence behaviour of the
  error and the estimator are given in Fig. \ref{abb:ex3 error vs estimator}. 
  As expected, we observe a linear rate. The numerical solutions obtained
  on different refinement levels are provided in
  Fig. \ref{abb:ex3 numerical solutions uH}. We see that the sequence
  of iterates from the adaptive refinement process converges to the 
  physical true solution. However, the first few iterates obtain a
  constant shift from the true solution and some minor oscillations at
  $x=0$ and at $x=1$, respectively. Note that this behaviour was also 
  observed in {\rm \cite{BacutaHayes:2024,CohenDahmenWelper:2012}}. Further,
  we see that the numerical solutions on the first few refinement levels have
  negative values even though the true solution is not negative. Therefore,
  the discrete maximum principle is not satisfied for all the solutions
  from the adaptive refinement process. However, the inbuilt error estimator
  detects how many degrees of freedom need to be added in order to obtain
  a numerical solution $u_H$, which is in good accordance with the physical
  correct reference solution $u$. Finally, we want to mention that in the case
  of $\varepsilon = 10^{-5}$ the numerical solution $u_H$ on $L=20$ was
  computed on a mesh with $60$ vertices (indicated in black) and $179$
  degrees of freedom (dof) for the corresponding saddle point system. This
  already lead to a satisfactory result. In comparison, the direct approach
  needed about $100000$ dofs to
  give a satisfactory approximation to the solution, 
  see Fig. \ref{abb:ex_bndlayer1d_error_plot}
  and Tab. \ref{tab:Erros for advdif1d uniform}.
  
  \begin{figure}
    \centering
    \begin{subfigure}{0.48\textwidth}
      \includegraphics[width=0.98\textwidth]
      {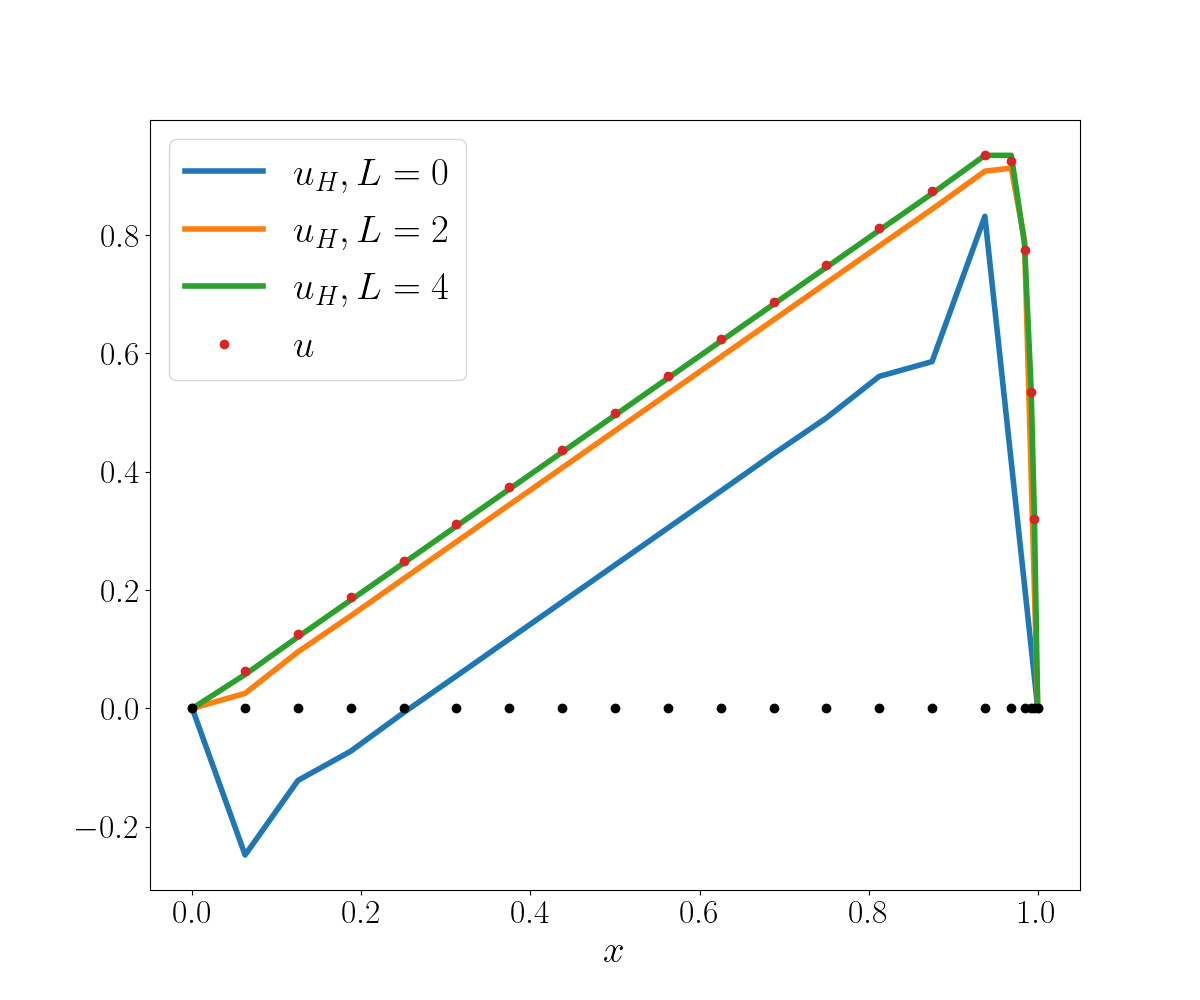}
      \caption{$\varepsilon = 10^{-2}$}
    \end{subfigure}
    \begin{subfigure}{0.48\textwidth}
      \includegraphics[width = 0.98\textwidth]
      {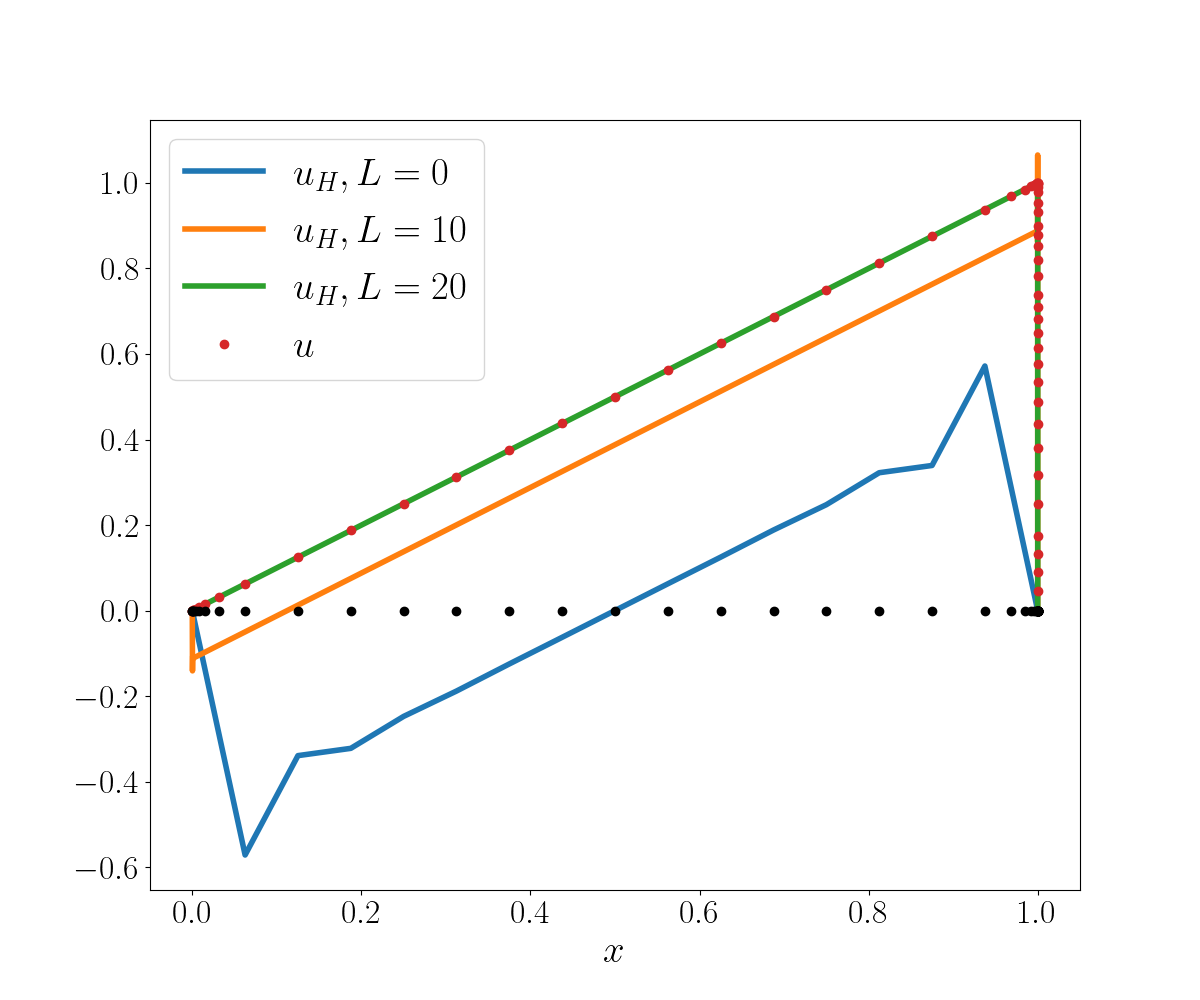}
      \caption{$\varepsilon = 10^{-5}$}
    \end{subfigure}
    \caption{Numerical solutions $u_H$ to Example \ref{Example 2.3} in the
      case of the adaptive refinement scheme.}
    \label{abb:ex3 numerical solutions uH}
  \end{figure}

  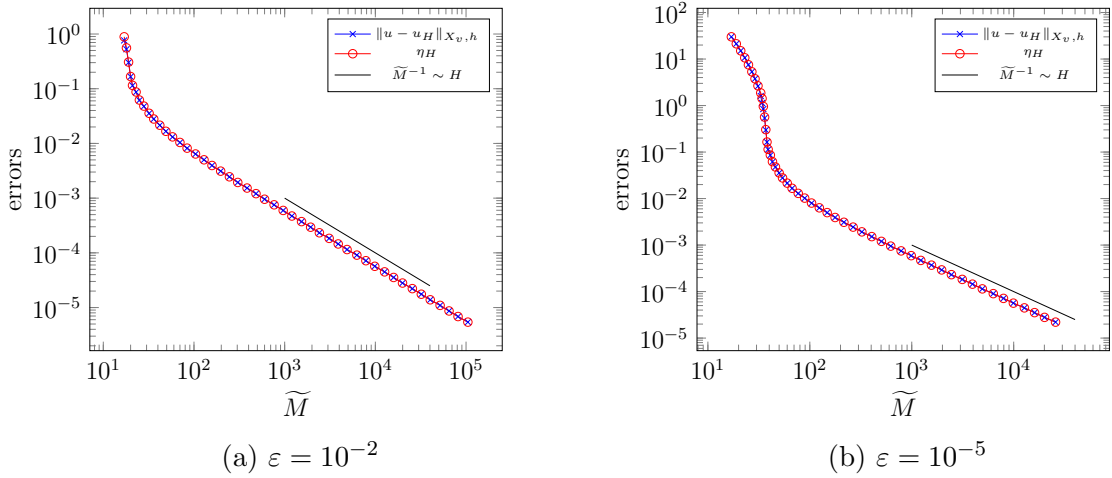
\begin{figure}
    \centering
    \begin{subfigure}{0.49\textwidth}
      \begin{tikzpicture}[scale = 0.8]
        \begin{axis}[
            xmode = log,
            ymode = log,
            xlabel= $\widetilde{M}$,
            ylabel=errors,
            legend style={font=\tiny}, legend pos = north east]
            \addplot[mark = x,blue] table [col sep=&, y=errXh_u, x=nv]
            {tables/SimData_LS_laplace_advdiff1d_eps1e-02_adaptive.dat};
            \addlegendentry{$\|u - u_H\|_{X_v,h}$}
        
            \addplot[mark = o,red] table [col sep=&, y=L2_dxph, x=nv]
            {tables/SimData_LS_laplace_advdiff1d_eps1e-02_adaptive.dat};
            \addlegendentry{$\eta_H$}
              
            \addplot[
            domain = 1000:40000,
            samples = 10,
            thin,
            black,
          ] {1*x^(-1)};
          \addlegendentry{$\widetilde{M}^{-1}\sim H$}
        \end{axis}
      \end{tikzpicture}
      \caption{$\varepsilon = 10^{-2}$}
    \end{subfigure}
    \begin{subfigure}{0.49\textwidth}
      \begin{tikzpicture}[scale = 0.8]
        \begin{axis}[
            xmode = log,
            ymode = log,
            xlabel= $\widetilde{M}$,
            ylabel=errors,
            legend style={font=\tiny}, legend pos = north east]
            \addplot[mark = x,blue] table [col sep=&, y=errXh_u, x=nv]
            {tables/SimData_LS_laplace_advdiff1d_eps1e-05_adaptive.dat};
            \addlegendentry{$\|u - u_H\|_{X_v,h}$}
        
            \addplot[mark = o,red] table [col sep=&, y=L2_dxph, x=nv]
            {tables/SimData_LS_laplace_advdiff1d_eps1e-05_adaptive.dat};
            \addlegendentry{$\eta_H$}
              
            \addplot[
            domain = 1000:40000,
            samples = 10,
            thin,
            black,
          ] {1*x^(-1)};
          \addlegendentry{$\widetilde{M}^{-1}\sim H$}
        \end{axis}
      \end{tikzpicture}
      \caption{$\varepsilon = 10^{-5}$}
    \end{subfigure}
    
    \caption{Convergence behaviour of the errors and the estimators.}
    \label{abb:ex3 error vs estimator}
  \end{figure}
\end{example}

\begin{example}[{\cite[Example 3.3]{JohnKnobloch:2023}}]
  This problem deals with a non-constant convection field {\boldmath{$v$}}.
  In particular, we consider $\Omega = (0,1)^2$ with
  {\boldmath{$v$}}$=(-y,x)^\top$, $\alpha(x) = \varepsilon = 10^{-5}$, and 
  $f=0$. We prescribe homogeneous Dirichlet boundary conditions 
  on the boundaries $\{1\}\times \left[0,1\right]$ and
  $\left[0,1\right]\times \{1\}$, i.e., on the right and top boundary.
  At the inlet boundary $\left[0,1\right]\times \{0\}$ we consider the 
  inhomogeneous boundary condition given by,  with $\xi = 10^{-3}$,
  \begin{equation}\label{eq:inhomogeneous bnd condition}
    u(x,0) = 
    \begin{cases}
      1 - \frac{1}{4}\left( 1 - \cos\left(\frac{1/3 + \xi  -x }{2 \xi}
          \pi \right)\right)^2
      & \text{for }x \in \left[\frac{1}{3} - \xi, \frac{1}{3} + \xi \right], \\
      1 & \text{for }x \in \left(\frac{1}{3} + \xi, \frac{2}{3} - \xi
      \right), \\
      1 - \frac{1}{4}\left( 1 - \cos\left(\frac{x-2/3 + \xi }{2 \xi}
          \pi \right)\right)^2 
      & \text{for }x \in \left[\frac{2}{3} - \xi , \frac{2}{3} + \xi\right], \\
      0 & \text{else}.
    \end{cases}
  \end{equation}
  On the remaining outlet boundary $\{0\} \times \left(0,1\right)$ we
  prescribe homogeneous Neumann boundary conditions. In our implementation
  we use piecewise linear and piecewise quadratic basis functions to define
  the finite element spaces $X_H$ and $Y_h$, respectively. Furthermore, we
  use the inbuilt error estimator to drive an adaptive refinement scheme
  with a Dörfler marking strategy {\rm \cite{Doerfler:1996}}. All resulting
  linear systems are solved with the sparse direct solver \texttt{Pardiso}.
  In order to study the satisfaction of the global
  discrete maximum principle (DMP)
  we evaluate as in {\rm \cite{JhaJohnKnobloch:2023}} the quantity 
  \begin{equation}\label{eq:global osc}
    \osc_{max}(u_H) = \max _{(x,y) \in \overline{\Omega}}u_H(x,y) -1 - 
    \min_{(x,y)\in \overline{\Omega}} u_H(x,y).
  \end{equation}
  In order to assess the accuracy of the numerical solution three
  characteristic values of the solution at the outflow boundary are
  provided in {\rm \cite{JohnKnobloch:2023}}. The reference values read:
  \begin{itemize}
    \item width of the lower layer: $0.01439869$,
    \item width of the upper layer:  $0.01439637$,
    \item outflow profile width: $0.3482541$.
  \end{itemize}
  
  \begin{figure}
    \centering
    \begin{subfigure}{0.48\textwidth}
      \includegraphics[width=0.98\textwidth]
      {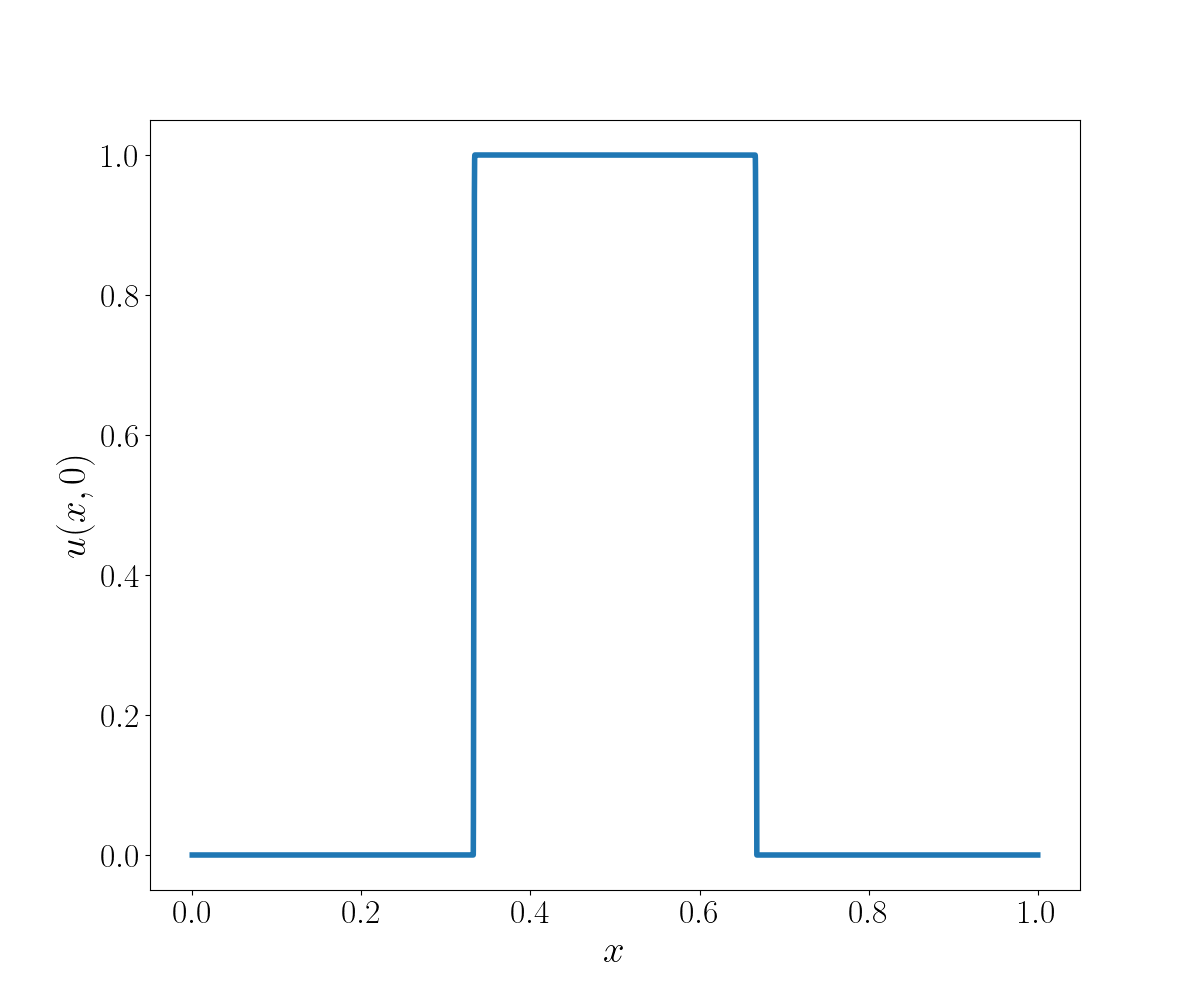}
    \end{subfigure}
    \begin{subfigure}{0.48\textwidth}
      \includegraphics[width = 0.98\textwidth]
      {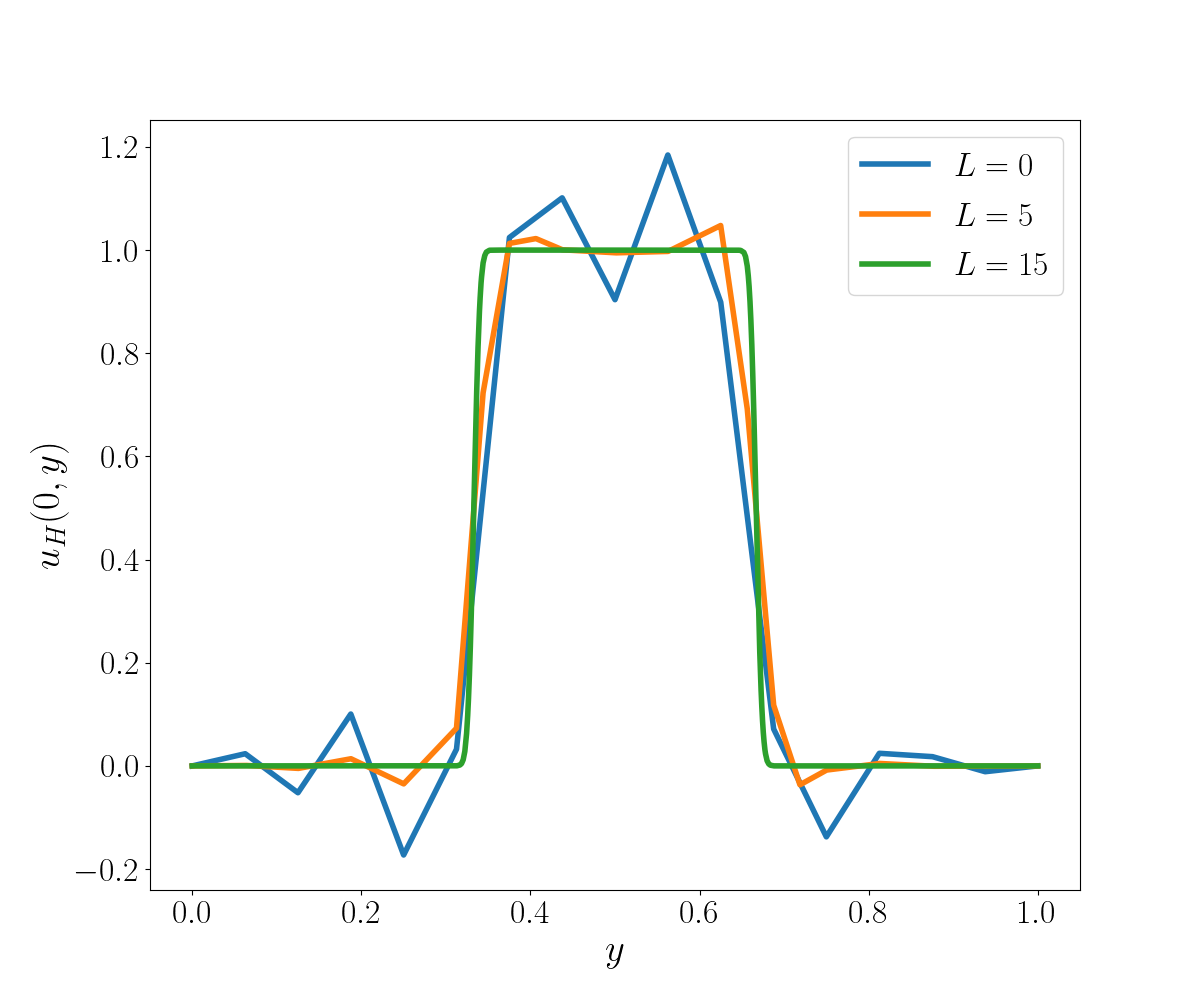}
    \end{subfigure}
    \caption{Left: Inhomogeneous boundary condition
      \eqref{eq:inhomogeneous bnd condition} prescribed on 
      inlet boundary. Right: Numerical solution $u_H$ at the outlet
      boundary for different refinement levels.}
    \label{abb:example laplace nonconstant convection inlet outlet bnd}
  \end{figure}

  \begin{figure}
    \centering
    \begin{subfigure}{0.48\textwidth}
      \includegraphics[trim = 550px 300px 550px 100px,clip,
      width=0.97\textwidth]
      {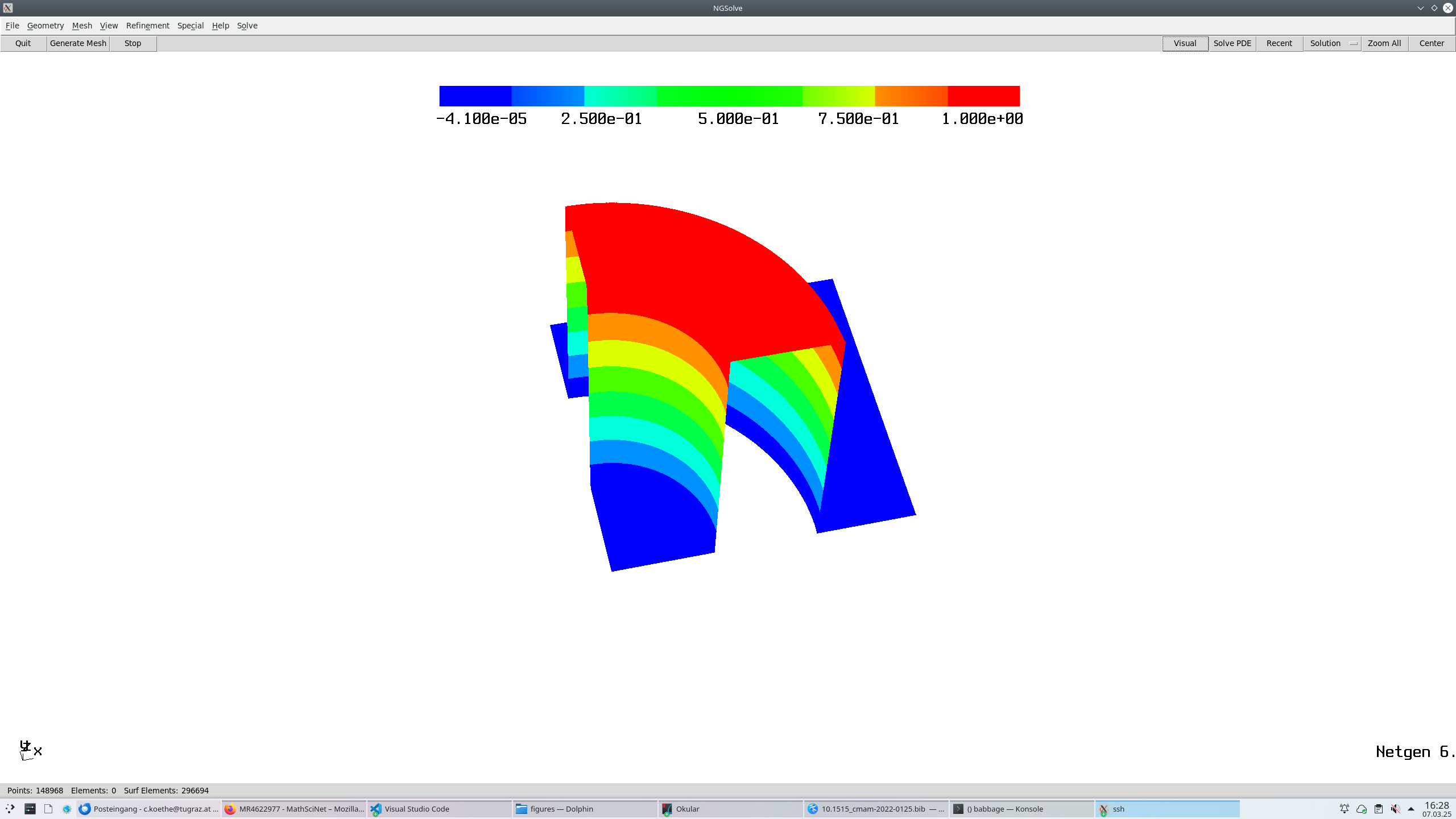}
    \end{subfigure}
    \begin{subfigure}{0.48\textwidth}
      \includegraphics[trim = 450px 200px 450px 100px,clip,
      width = 0.97\textwidth]
      {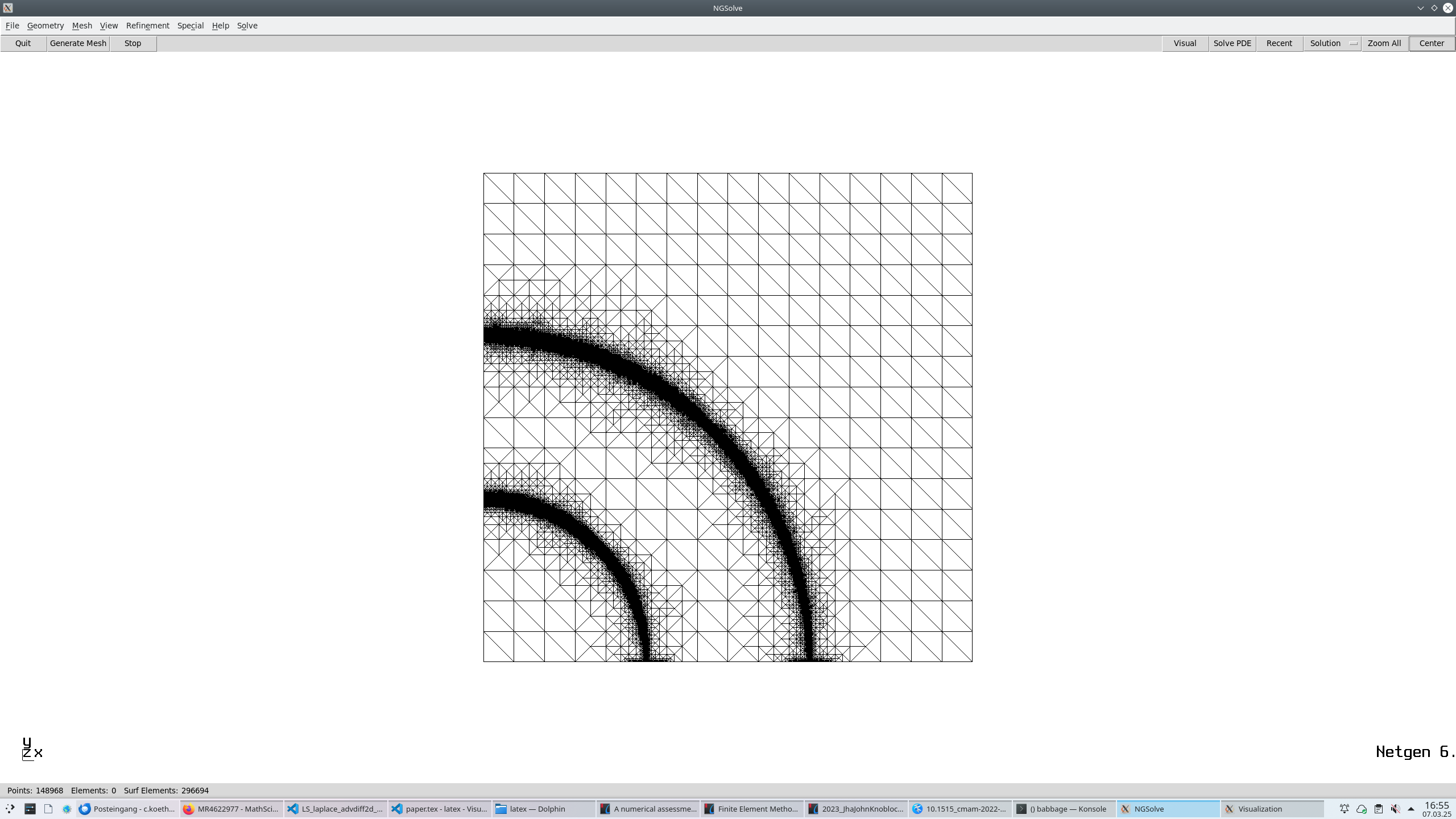}
    \end{subfigure}
    \caption{Left: Numerical solution $u_H$ on $L=15$. Right: Adaptive
      mesh on $L=15$ with 148968 dofs.}
    \label{abb:example laplace nonconstant convection uH and mesh}
  \end{figure}

  \begin{figure}
    \centering
    \begin{subfigure}{0.48\textwidth}
      \begin{tikzpicture}[scale = 0.8]
        \begin{axis}[
            xmode = log,
            ymode = log,
            xlabel= $\widetilde{M}$,
            ylabel=difference to reference values,
            legend style={font=\tiny}, legend pos = north east]
            \addplot[mark = x,blue] table [col sep=&, y=diff_profile, x=nv]
            {tables/SimData_LSLaplace_advdiff2d_nonconstantConvection.dat};
            \addlegendentry{profile width}
        
            \addplot[mark = x,red] table [col sep=&, y=diff_lower, x=nv]
            {tables/SimData_LSLaplace_advdiff2d_nonconstantConvection.dat};
            \addlegendentry{lower width}
            
            \addplot[mark = x,green!60!black] table [col sep=&,
            y=diff_upper, x=nv]
            {tables/SimData_LSLaplace_advdiff2d_nonconstantConvection.dat};
            \addlegendentry{upper width}
        \end{axis}
      \end{tikzpicture}
    \end{subfigure}
    \begin{subfigure}{0.48\textwidth}
      \begin{tikzpicture}[scale = 0.8]
        \begin{axis}[
            xmode = log,
            xlabel= $\widetilde{M}$,
            ylabel=$\osc_{max}(u_H)$,
            legend style={font=\tiny}, legend pos = north east]
            \addplot[mark = x,blue] table [col sep=&, y=oscmax, x=nv]
            {tables/SimData_LSLaplace_advdiff2d_nonconstantConvection.dat};
        \end{axis}
      \end{tikzpicture}
    \end{subfigure}
    \caption{Left: Comparison of the computed characteristic values to
      the reference values. Right: Satisfaction 
    of the global DMP on different refinement levels.}
  \label{abb:Comparison of characteristic values and evaluation of
    osc for nonconstant convection}
  \end{figure}

  \begin{table}
    \begin{center}
      \input{tables/SimData_LS_laplace_advdiff2d_nonconstantconvection.tex}
    \end{center}
    \caption{Computed characteristic values and evaluation of
      \eqref{eq:global osc}.}
    \label{tab:Computed characteristic values and osc for nonconstant
      convection}
  \end{table}
  
  \noindent
  In Fig. \ref{abb:example laplace nonconstant convection inlet outlet bnd}
  we provide a plot of the inhomogeneous boundary condition
  \eqref{eq:inhomogeneous bnd condition} on the inlet boundary as well as
  the numerical solution $u_H$ at the outlet boundary for different
  refinement levels.
  In Fig. \ref{abb:example laplace nonconstant convection uH and mesh} 
  the numerical solution $u_H$ as well as the adaptive mesh generated on 
  refinement level $L=15$ are depicted. The adaptive mesh was generated from 
  a structured initial mesh with $16\times 16$ elements. Further, in Tab. 
  \ref{tab:Computed characteristic values and osc for nonconstant convection}
  and  Fig. \ref{abb:Comparison of characteristic values and evaluation of osc
    for nonconstant convection} we provide a comparison 
  of the computed characteristic values to the reference values and the
  satisfaction of the global DMP. We see that for
  earlier refinement levels over and under shoots are 
  visible, but they are almost vanishing on higher refinement levels.
  This can be also seen from the quantity \eqref{eq:global osc} in Table 
  \ref{tab:Computed characteristic values and osc for nonconstant convection}
  and Fig. 
  \ref{abb:Comparison of characteristic values and evaluation of osc for
    nonconstant convection}, which show a good 
  satisfaction of the global DMP for higher refinements. In addition to that 
  Fig. \ref{abb:Comparison of characteristic values and evaluation of osc for
    nonconstant convection} shows that the difference of the computed
  characteristic quantities to the reference values is in the order of
  floating point precision for higher refinement levels, i.e., 
  we have a good agreement of the computed reference values with the
  characteristic values.
\end{example}

\begin{example}[{Hemker problem}]
  In this example we consider the Hemker problem {\rm \cite{Hemker:1996}},
  which is a standard benchmark problem for steady state convection-diffusion
  problems. The domain is given by
  $\Omega = (-3,9) \times (-3,3) \setminus \{(x,y):\, x^2 +y^2 \leq 1 \}$.
  The velocity is given by {\boldmath{$v$}}$=(1,0)^\top$, and the right-hand
  side $f$ is set equal to zero. Further we have the boundary conditions 
  \begin{equation*}
    u(x,y) = 
    \begin{cases}
      0 & x=-3,\\
      1 & x^2 +y^2 = 1, \\
      \alpha(x)\nabla u \cdot \mathbf{n} = 0
      & x=9 \, \lor \, y=-3 \, \lor \, y=3   
    \end{cases}.
  \end{equation*}
  As in {\rm \cite{AugustinJohn:2011,JhaJohnKnobloch:2023,JohnKnobloch:2023}}
  we consider the diffusion coefficient to be
  $\alpha(x) = \varepsilon = 10^{-4}$.
  In order to assess the accuracy of the numerical solution a value for the
  width of the interior layer at $x=4$ was provided in
  {\rm \cite{AugustinJohn:2011}}. This width is defined to be the length of the
  interval, where $u(4,y) \in \left[0.1,0.9\right]$.
  In {\rm \cite{AugustinJohn:2011}} the reference value $0.0723$ is provided
  for the upper layer, i.e., where $y\geq 0$. Furthermore, we evaluate the
  quantity \eqref{eq:global osc} to measure the satisfaction of the global
  DMP. In our numerical experiments we use piecewise linear finite elements
  for the trial space $X_H$ and piecewise quadratic finite elements for the
  test space $Y_h$. Further we use the local error indicators
  \eqref{local error indicator} to drive an adaptive refinement scheme. As
  a marking strategy we use the Dörfler criterion with parameter
  $\theta = 0.5$. The marked elements are then refined using newest
  vertex bisection. The resulting linear systems are solved with the sparse
  direct solver \texttt{Pardiso}.

  \noindent 
  In Fig. \ref{abb:Hemker meshes} the initial mesh and the adaptive mesh
  obtained on level $L=17$ are depicted. We see stronger refinements at the
  boundary layer around the circle, which at the top and bottom of the circle
  passes into an interior layer that spread into the direction of the
  convection. The refinements at the left boundray may be explained in terms
  of the over- and undershoots that occur for lower refinements.
  In Fig. \ref{abb:Hemker uH 676072 dofs} we depict the numerical solution
  $u_H$ obtained on refinement level $L=17$.
  In Fig. \ref{abb:Hemker cut lines at x4 and y1} we provide plots of the
  cut lines $u_H(x,1)$ for $-3\leq x \leq 9 $ and $u_H(4,y)$ for
  $-3\leq y \leq 3$. The cut lines possess some oscillations for lower
  refinement levels which get significantly reduced for higher refinement
  levels. In Fig. \ref{abb:Comparison of characteristic values and evaluation
    of osc for Hemker problem} and
  Tab. \ref{tab:Computed characteristic values and osc for Hemker problem} 
  we provide a comparison of the computed characteristic value with the
  reference value as well as a plot of the satisfaction of the global DMP.
  We see that the computed value converges to the reference value given for
  the upper layer width. Furthermore, we see a reduction of the quantity 
  \eqref{eq:global osc} for higher refinements.

  \begin{figure}
    \centering
    \begin{subfigure}{0.48\textwidth}
      \includegraphics[trim = 490px 300px 550px 300px,clip,
      width=0.98\textwidth]
      {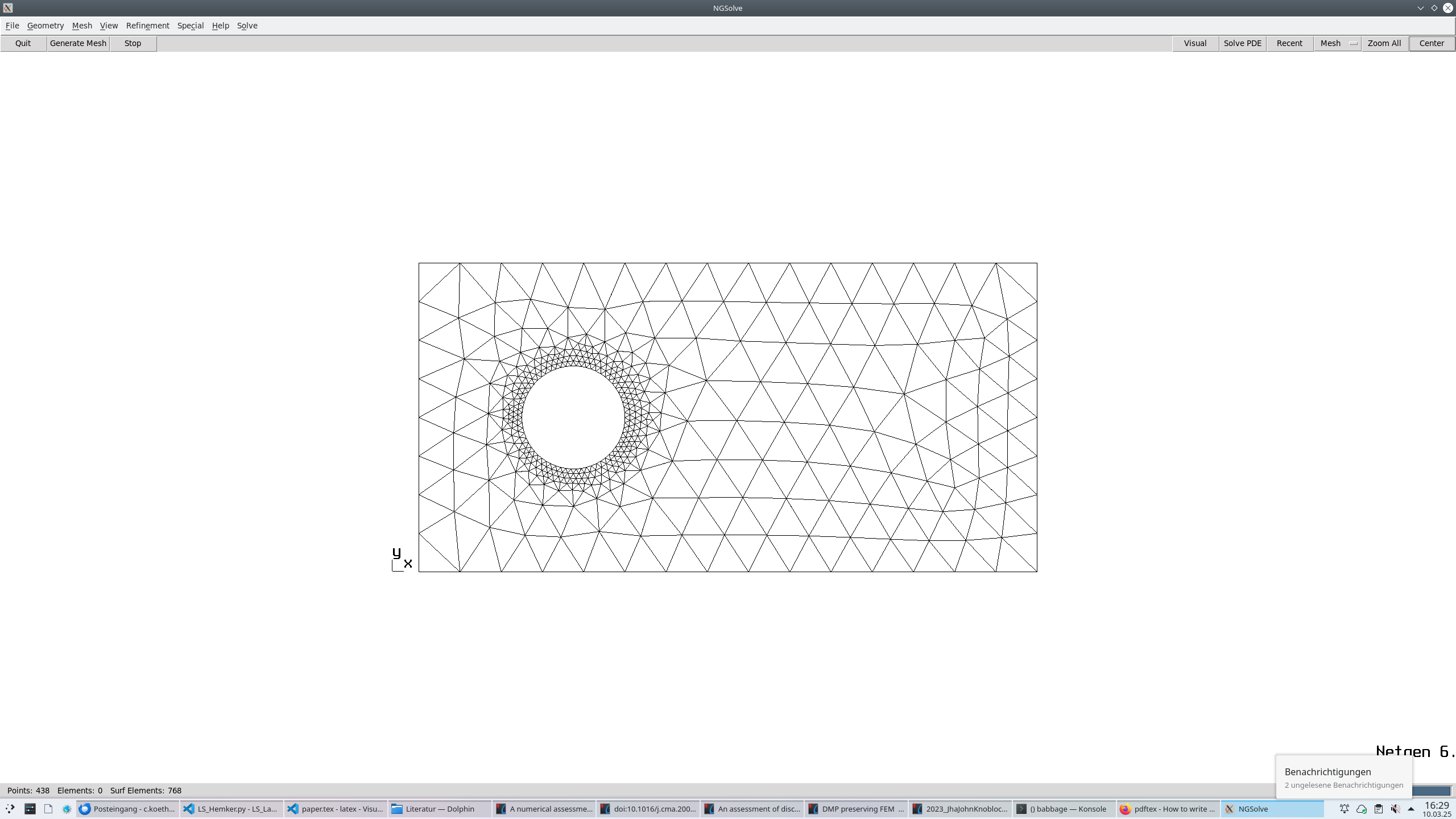}
    \end{subfigure}
    \begin{subfigure}{0.48\textwidth}
      \includegraphics[trim = 490px 300px 550px 300px,clip,
      width = 0.98\textwidth]
      {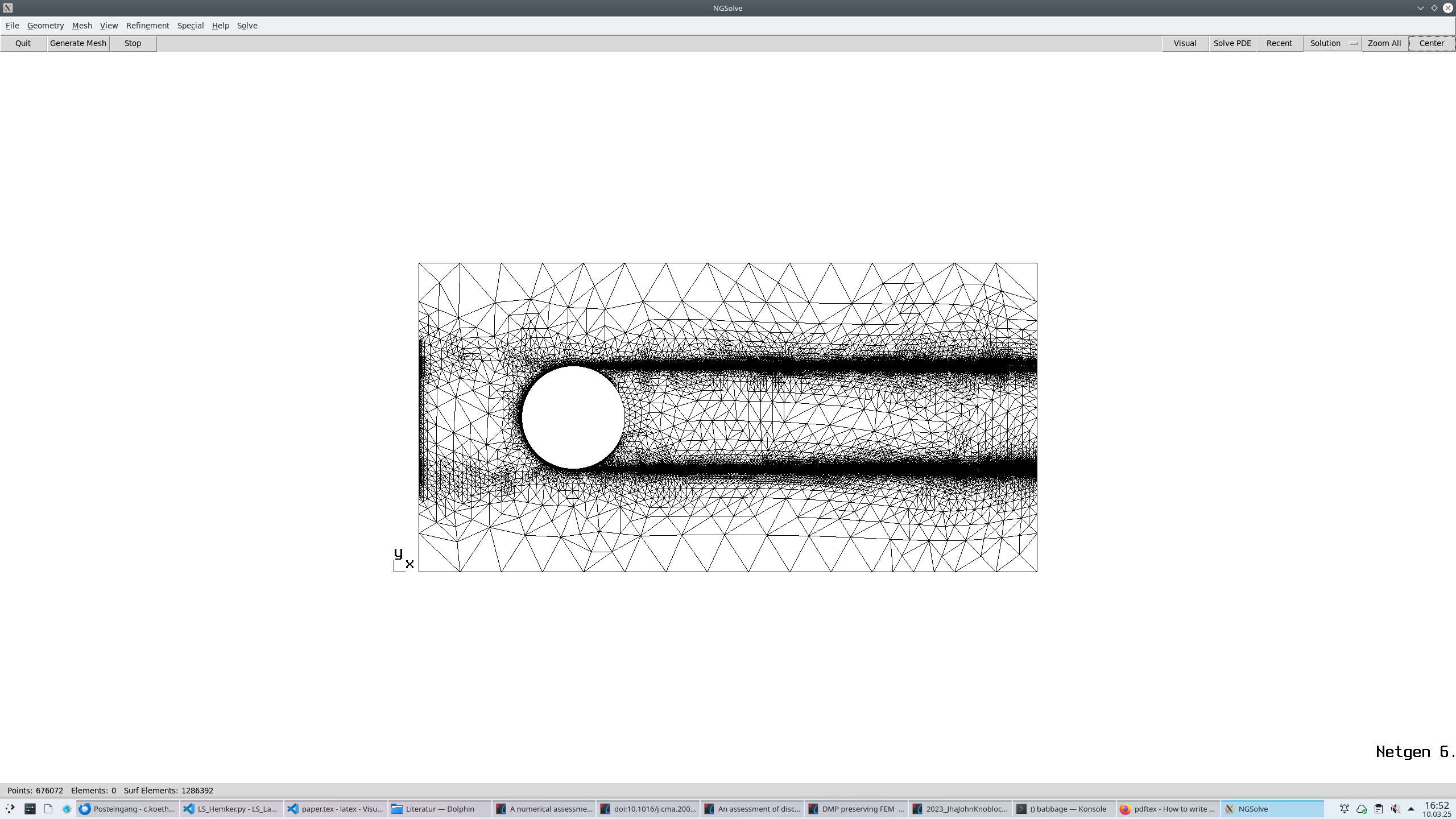}
    \end{subfigure}
    \caption{Left: Initial mesh with 438 dofs. Right: Adaptive mesh
      on $L=17$ with 676072 dofs.}
    \label{abb:Hemker meshes}
  \end{figure}

  \begin{figure}
    \centering
    \includegraphics[trim = 450px 300px 450px 300px,clip,
    width=0.8\textwidth]{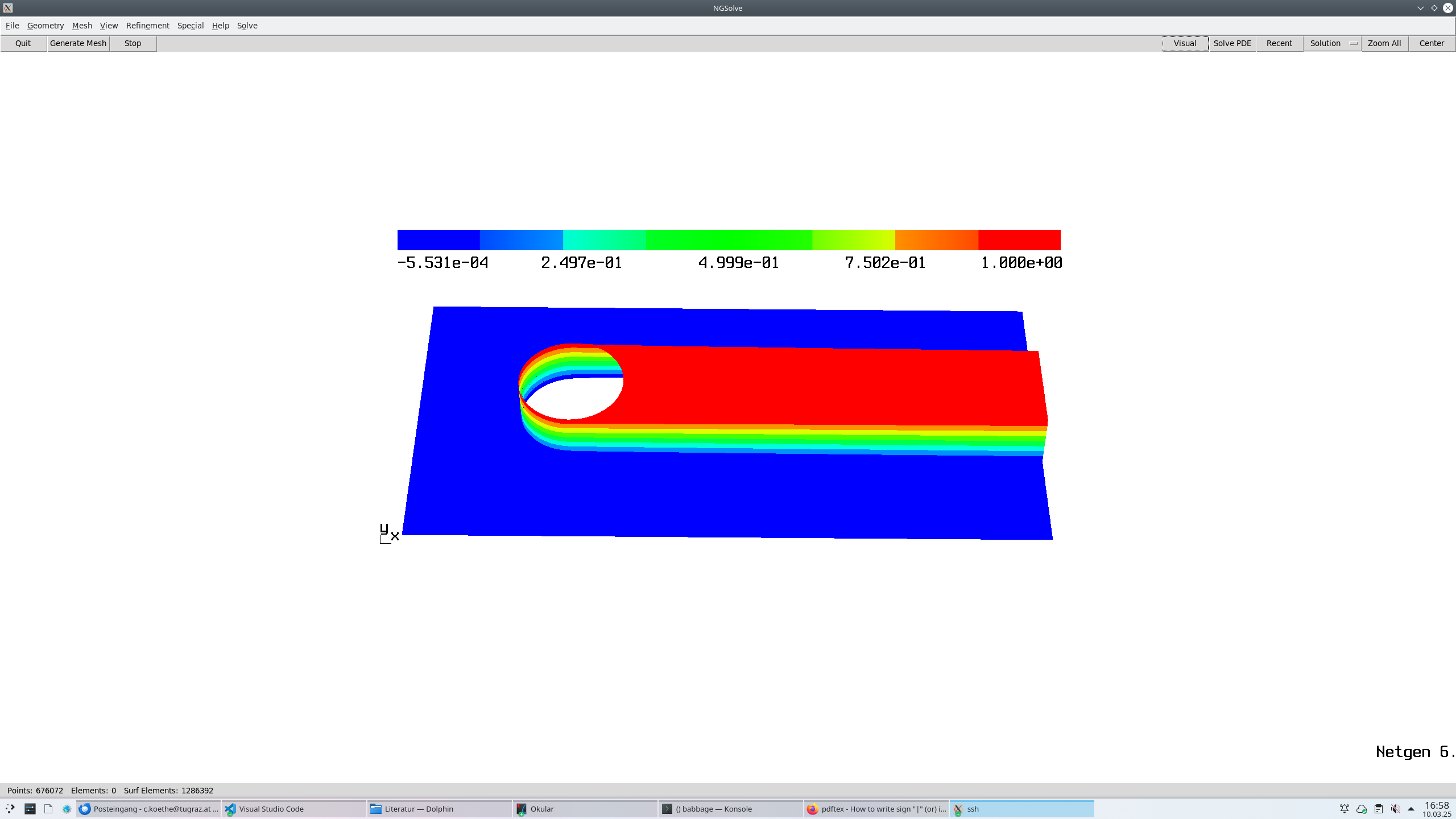}
    \caption{Numerical solution $u_H$ on $L=17$.}
    \label{abb:Hemker uH 676072 dofs}
  \end{figure}

  \begin{figure}
    \centering
    \begin{subfigure}{0.48\textwidth}
      \includegraphics[width=0.98\textwidth]{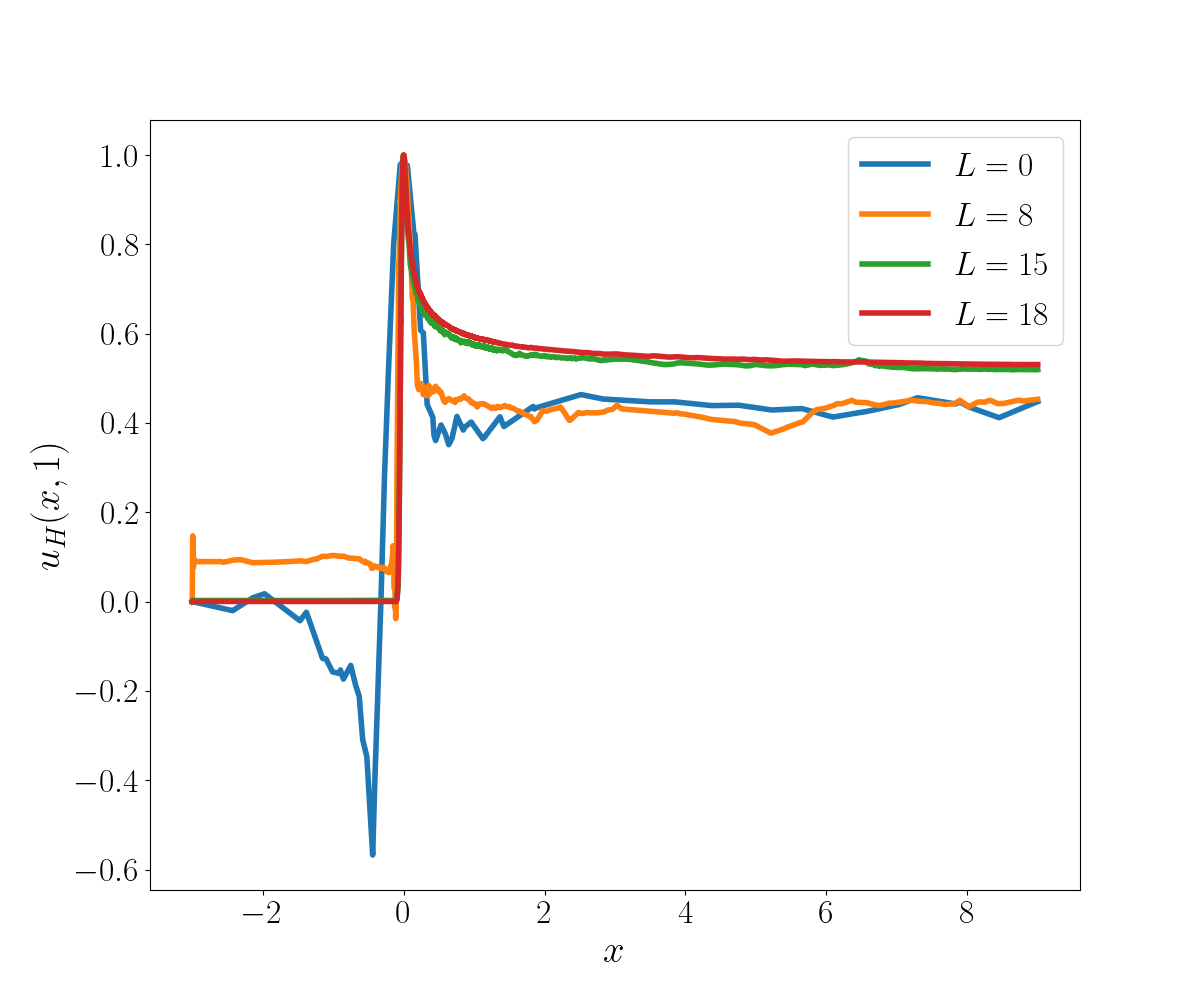}
    \end{subfigure}
    \begin{subfigure}{0.48\textwidth}
      \includegraphics[width = 0.98\textwidth]{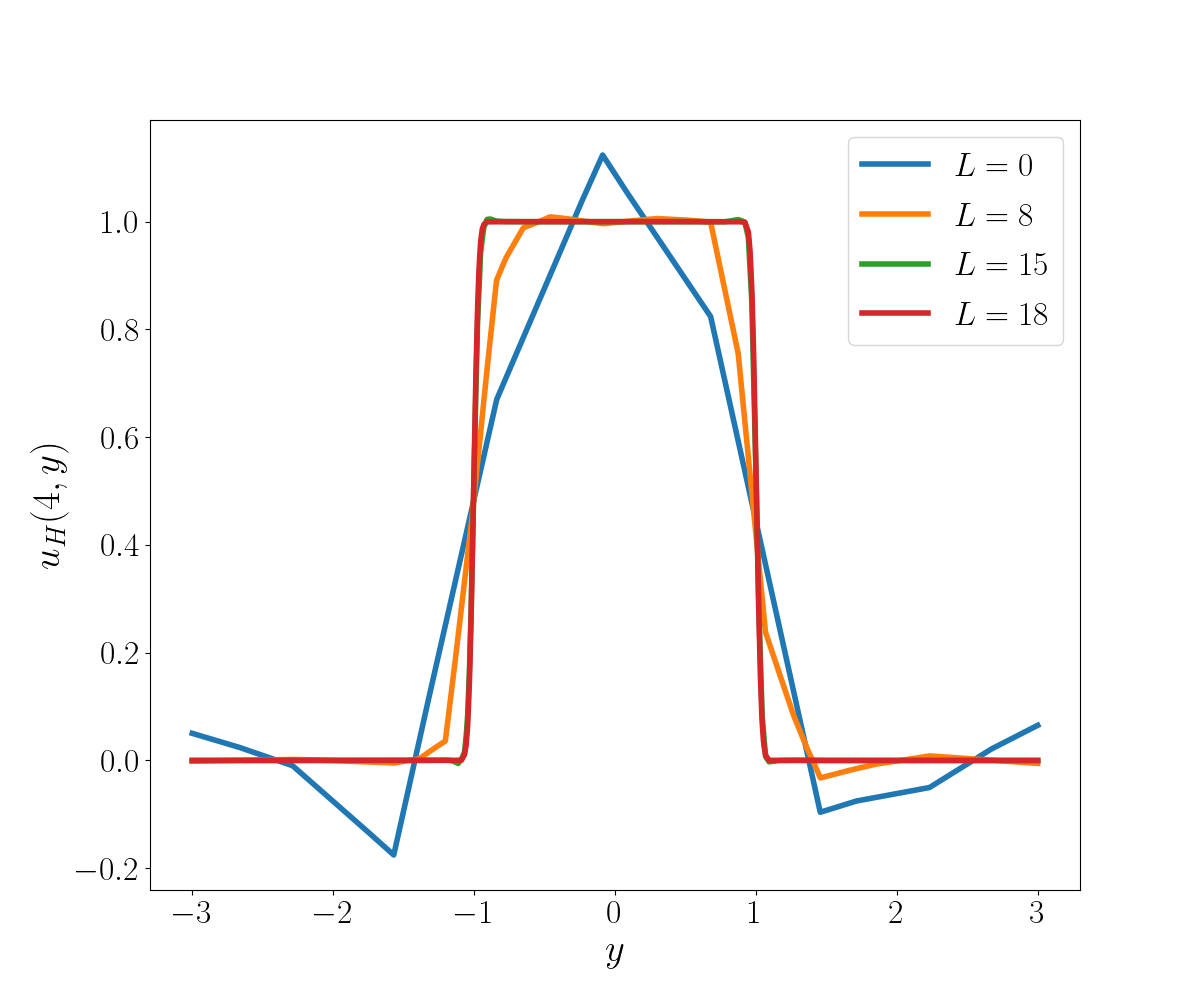}
    \end{subfigure}
    \caption{Cut lines of the numerical solution on different refinement
      levels. Left: $u_H(x,1)$. Right: $u_H(4,y)$.}
    \label{abb:Hemker cut lines at x4 and y1}
  \end{figure}

  \begin{figure}
    \centering
    \begin{subfigure}{0.48\textwidth}
      \begin{tikzpicture}[scale = 0.8]
        \begin{axis}[
            xmode = log,
            ymode = log,
            xlabel= $\widetilde{M}$,
            ylabel=difference to reference value,
            legend style={font=\tiny}, legend pos = north east]
            \addplot[mark = x,blue] table [col sep=&, y=diff_upper, x=nv]
            {tables/SimData_LSHemker.dat};
            \addlegendentry{upper width}
        \end{axis}
      \end{tikzpicture}
    \end{subfigure}
    \begin{subfigure}{0.48\textwidth}
      \begin{tikzpicture}[scale = 0.8]
        \begin{axis}[
            xmode = log,
            xlabel= $\widetilde{M}$,
            ylabel=$\osc_{max}(u_H)$,
            legend style={font=\tiny}, legend pos = north east]
            \addplot[mark = x,blue] table [col sep=&, y=oscmax, x=nv]
            {tables/SimData_LSHemker.dat};
        \end{axis}
      \end{tikzpicture}
    \end{subfigure}
    \caption{Left: Comparison of the computed characteristic values
      to the reference values. Right: Satisfaction 
    of the global DMP on different refinement levels.}
  \label{abb:Comparison of characteristic values and evaluation
    of osc for Hemker problem}
  \end{figure}
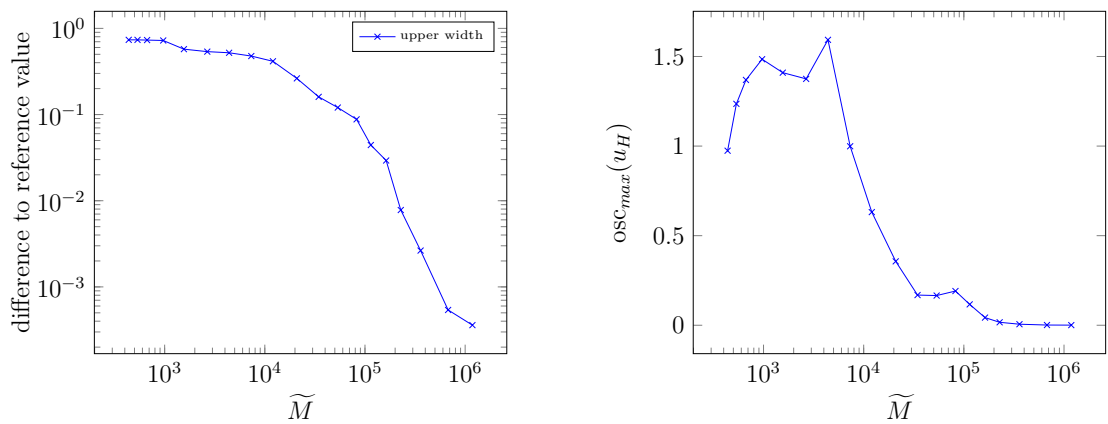

  \begin{table}
    \begin{center}
      \input{tables/SimData_Hemker.tex}
    \end{center}
    \caption{Computed characteristic value and evaluation of
      \eqref{eq:global osc} for the Hemker problem.}
    \label{tab:Computed characteristic values and osc for Hemker problem}
  \end{table}
\end{example}

\section{Instationary convection-diffusion problems}
Instead of \eqref{DBVP Advection Diffusion} we now consider the time-dependent
Dirichlet boundary value problem for the convection-diffusion equation,
\begin{equation}\label{IBVP Advection Diffusion}
  \begin{array}{rclcl}
    \partial_t u(y,t) +
    \mbox{\boldmath $v$}(y,t) \cdot \nabla_y u(y,t)
    - \mbox{div}_y [\alpha(y) \, \nabla_y u(y,t) ] 
    & = & f(y,t)
    & & \mbox{for} \; (y,t) \in Q, \\[1mm]
    u(y,t)
    & = & 0
    & & \mbox{for} \; (y,t) \in \Sigma, \\[1mm]
    u(x,0) & = & 0 & & \mbox{for} \; x \in \Omega ,
  \end{array}
\end{equation}
where the space-time domain $Q$ is given by
\[
  Q := \Big \{ (y,t) \in {\mathbb{R}}^{n+1} : y = \varphi(t,x),
  x \in \Omega \subset {\mathbb{R}}^n, t \in (0,T) \Big \},
\]
with the lateral boundary
\[
  \Sigma := \Big \{ (y,t) \in {\mathbb{R}}^{n+1} : y = \varphi(t,x),
  x \in \partial \Omega, t \in (0,T) \Big \}.
\]
Here, $T >0$ is a given time horizon, and
$\mbox{\boldmath $v$}(y,t) = \frac{d}{dt}y(t)$ is the velocity along the
trajectory $y(t) = \varphi(t,x) \in {\mathbb{R}}^n$ for a reference point
$x \in \Omega$. We assume that the deformation $\varphi$ is bijective and
sufficient regular for all $t \in (0,T)$, satisfying $\varphi(0,x)=x$
for all $x \in \Omega$, and $\mbox{div}_y \mbox{\boldmath $v$}(y,t)=0$.
Then we can write Reynold's transport theorem as
\begin{equation}\label{Reynolds}
  \frac{d}{dt} \int_{\Omega(t)} u(y,t) \, dy \, = \,
  \int_{\Omega(t)} \Big[ \partial_t u(y,t) + \mbox{\boldmath $v$}(y,t) \cdot
  \nabla_y u(y,t) \Big] \, dy \, = \,
  \int_{\Omega(t)} \frac{d}{dt}u(y,t) \, dy.
\end{equation}
As in the stationary case we assume that the diffusion coefficient $\alpha(y)$
is bounded and strictly positive, see \eqref{Assumption alpha}. 

\subsection{Variational formulation}
Similar as in the stationary case we now define the Bochner spaces
$Y := L^2(0,T;H^1_0(\Omega(t)))$ with the norm
\[
  \| q \|_Y^2 := \int_0^T \int_{\Omega(t)} \alpha(y) \, |\nabla_y q(y,t)|^2 \,
  dy \, dt,
\]
and
\[
  X_v := \Big \{ u \in Y : \frac{d}{dt} u =
  \partial_t u + \mbox{\boldmath $v$} \cdot \nabla_y u \in Y^*,
  u(x,0) = 0 \; \mbox{for} \; x \in \Omega \Big \} .
\]
A norm in $X_v$ is given by the graph norm
\[
  \| u \|_{X_v} := \sqrt{\| u \|_Y^2 +
    \| \partial_t u + \mbox{\boldmath $v$} \cdot \nabla_y u \|_{Y^*}^2} =
  \sqrt{\| u \|_Y^2 + \| w_u \|_Y^2},
\]
where $w_u \in Y$ is the unique solution of the variational problem
such that
\begin{eqnarray}\label{Definition wu ST}
  && a(w_u,q) := \int_0^T \int_{\Omega(t)} \alpha(y) \, \nabla_y w_u(y,t) \cdot
     \nabla_y q(y,t) \, dy \, dt \\
  && \hspace*{4cm} = \,
     \int_0^T \int_{\Omega(t)} \Big[
     \partial_t u(y,t) + \mbox{\boldmath $v$}(y,t) \cdot \nabla_y u(y,t)
     \Big] \, q(y,t) \, dy \, dt \nonumber
\end{eqnarray}
is satisfied for all $q \in Y$. The space-time variational formulation of
\eqref{IBVP Advection Diffusion} then reads to find $u \in X_v$ such that
\begin{eqnarray}\label{IBVP VF}
  b(u,q) & := & \int_0^T \int_{\Omega(t)} \Big[
  \partial_t u(y,t) + \mbox{\boldmath $v$}(y,t) \cdot \nabla_y u(y,t)
                \Big] q(y,t) \, dy \, dt \\
         && \nonumber
            + \int_0^T \int_{\Omega(t)} \alpha(y) \, \nabla_y u(y,t) \cdot
            \nabla_y q(y,t) \, dy \, dt \, = \,
            \int_0^T \int_{\Omega(t)} f(y,t) q(y,t) \, dy \, dt
\end{eqnarray}
is satisfied for all $q \in Y$. The boundedness of the bilinear form
$b(\cdot,\cdot) : X \times Y \to {\mathbb{R}}$ follows as in the proof
of Lemma \ref{Lemma solvability}, i.e.,
\[
  |b(u,q)| \leq \sqrt{2} \, \| u \|_{X_v} \| q \|_Y \quad
  \mbox{for all} \; (u,q) \in X_v \times Y .
\]
Moreover, choosing $q_u := u + w_u \in Y$, we conclude
\begin{eqnarray*}
  b(u,q_u)
  & = & \langle \partial_t u + \mbox{\boldmath $v$} \cdot \nabla_y u,
        u + w_u \rangle_Q \, + \, \langle \nabla_y u , \nabla_y (u + w_u)
        \rangle_{L^2(Q)} \\
  & = & 2 \, \langle \nabla_y w_u,
        \nabla_y u \rangle_Q \, + \, \langle \nabla_y w_u , \nabla_y w_u
        \rangle_{L^2(Q)} \, + \, \langle \nabla_y u , \nabla_y u
        \rangle_{L^2(\Omega)} \\
  & = & \langle \nabla_y (u+w_u) , \nabla_y (u+w_u) \rangle_{L^2(Q)}
        \, = \, \| q_u \|_Y^2 \, . 
\end{eqnarray*}
Using Reynold's transport theorem \eqref{Reynolds}
and $u(x,0)=0$ for $x \in \Omega$, this gives
\begin{eqnarray*}
  && \hspace*{-5mm} 2 \, \langle \partial_t u +
  \mbox{\boldmath $u$} \cdot \nabla_y u,u \rangle_Q
  = 2 \int_0^T \int_{\Omega(t)} \frac{d}{dt}u(y,t) \, u(y,t) \, dy \, dt \\
  && \hspace*{-3mm}
     = \int_0^T \int_{\Omega(t)} \frac{d}{dt}[u(y,t)]^2 \, dy \, dt \, = \,
     \int_0^T \frac{d}{dt} \int_{\Omega(t)} [u(y,t)]^2 \, dy \, dt \, = \,
     \int_{\Omega(T)} [u(y,T)]^2 dy \, \geq \, 0.
\end{eqnarray*}
Hence we obtain
\begin{eqnarray*}
  \| q_u \|^2_Y
  & = & \| u + w_u \|_Y^2 \, = \, \| u \|_Y^2 + \| w_u \|_Y^2
        + 2 \, \langle \nabla_y w_u , \nabla_y u \rangle_{L^2(Q)} \\
  & = & \| u \|_Y^2 + \| w_u \|^2_Y + 2 \, \langle \partial_t u +
        \mbox{\boldmath $u$} \cdot \nabla_y u,u \rangle_Q \\
  & \geq & \| u \|_Y^2 + \| w_u \|_Y^2 \, = \, \| u \|_{X_v}^2,
\end{eqnarray*}
and as in \eqref{inf-sup advection-diffusion} we therefore conclude
the inf-sup stability condition
\begin{equation}\label{inf-sup time}
  \| u \|_ {X_v} \leq \sup\limits_{0 \neq q \in Y}
  \frac{b(u,q)}{\| q \|_Y} \quad \mbox{for all} \; u \in X_v .
\end{equation}
In order to apply the Babu\v{s}ka--Ne\v{c}as theory it remains to prove
surjectivity.

\begin{lemma}\label{Lemma surjectivity}
  Consider the gradient of the velocity $\mbox{\boldmath $v$}(y,t)$,
  \[
    A_v = \left(
      \begin{array}{ccc}
        \displaystyle \frac{\partial}{\partial y_1} v_1
        & \displaystyle \frac{\partial}{\partial y_1} v_2
        & \displaystyle \frac{\partial}{\partial y_1} v_3 \\[4mm]
           \displaystyle \frac{\partial}{\partial y_2} v_1
        & \displaystyle \frac{\partial}{\partial y_2} v_2
        & \displaystyle \frac{\partial}{\partial y_2} v_3 \\[4mm]
           \displaystyle \frac{\partial}{\partial y_3} v_1
        & \displaystyle \frac{\partial}{\partial y_3} v_2
        & \displaystyle \frac{\partial}{\partial y_3} v_3 
      \end{array}
    \right)
  \]
  which is assumed to be positive semi-definite. Then, for all
  $q \in Y \backslash \{ 0 \}$ there exists a $u_q \in X_v$ such that
  $b(u_q,q) > 0$.
\end{lemma}

\proof{For any given $q \in Y\backslash \{ 0 \}$ we define $u_q \in X_v$,
  \[
    u_q(y,t) = u_q(\varphi(t,x),t) = \int_0^t q(\varphi(s,x),s) \, ds, \quad
    \frac{d}{dt} u_q(y,t) = q(y,t) .
  \]
  Hence we have
  \[
    b(u_q,q) = \langle \frac{d}{dt} u_q , q \rangle_Q +
    \langle \nabla_y u_q , \nabla_y q \rangle_{L^2(Q)} =
    \| q \|^2_{L^2(Q)} +
    \langle \nabla_y u_q , \nabla_y \frac{d}{dt} u_q \rangle_{L^2(Q)},
  \]
  and it is sufficient to prove that the second summand is non-negative.
  A direct computation gives
  \[
    \nabla_y u_q \cdot \nabla_y \frac{d}{dt} u_q \, = \,
    \frac{1}{2} \, \frac{d}{dt} \, |\nabla_y u_q|^2 +
    (A_v \nabla_y u_q , \nabla_y u_q ),
  \]
  and therefore,
  \begin{eqnarray*}
    \langle \nabla_y u_q , \nabla_y \frac{d}{dt} u_q \rangle_{L^2(Q)}
    & = & \int_0^T \int_{\Omega(t)} \left[
          \frac{1}{2} \, \frac{d}{dt} \, |\nabla_y u_q|^2 +
          (A_v \nabla_y u_q , \nabla_y u_q )
          \right] dy dt \\
    & & \hspace*{-2cm}
        = \,\frac{1}{2} \int_{\Omega(T)} |\nabla_y u_q(y,T)|^2 \, dy +
          \int_0^T \int_{\Omega(t)} (A_v \nabla_y u_q , \nabla_y u_q ) \,
          dy \, dt \, \geq \, 0
  \end{eqnarray*}
  follows. This concludes the proof.}

\begin{example}
  Consider the deformation $y = \varphi(t,x)$ to be a planar rotation,
  i.e., using polar coordinates we have
  \[
    y(t) = r \left(
      \begin{array}{c}
        \cos (\phi + \alpha t) \\[1mm]
        \sin (\phi + \alpha t)
      \end{array}
    \right),
    \quad
    \mbox{\boldmath $v$}(y,t) = \frac{d}{dt}y(t) =
    \left( \begin{array}{c} -y_2(t) \\[1mm] y_1(t)
             \end{array} \right).
  \]
  In this case we conclude
  \[
    A_v = \left(
      \begin{array}{cc}
        0 & 1 \\[1mm] -1 & 0
      \end{array}
    \right), \quad
    (A_v \nabla_y u , \nabla_y u) = 0 \, .
  \]
  This particular case was already considered in
  {\rm \cite{GanglGobrialSteinbach:2025}}.
\end{example}

\begin{example}
  For any velocity field
  $\mbox{\boldmath $v$}(y,t) = \mbox{\boldmath $v$}(t)$ we obviously
  have $A_v = 0$, and the assumption of Lemma \ref{Lemma surjectivity}
  is trivially satisfied.
\end{example}

\noindent
As in Lemma \ref{Lemma solvability} we now conclude unique
solvability of the variational formulation
\eqref{IBVP VF}, since all assumptions of the
Babu\v{s}ka--Ne\v{c}as theory are satisfied.

\subsection{Adaptive least-squares space-time finite element method}
Instead of the bilinear form $b(u,q)$ as defined in
\eqref{VF DBVP Advection Diffusion} we now consider the definition of
$b(u,q)$ as used in \eqref{IBVP VF} for $(u,q) \in X_v \times Y$,
and instead of the bilinear form $a(w_u,q)$ as defined in
\eqref{Definition wu} we now use the definition as given in
\eqref{Definition wu ST}. With these definitions we consider the minimization
of the quadratic functional \eqref{minimization problem} whose
minimizer $u \in X_v$ is given as the unique solution of the
gradient equation \eqref{gradient equation}. Instead of
\eqref{LS1} and \eqref{LS2} we therefore have to solve the coupled
system to find $(u,p) \in X_v \times Y$ such that
\begin{eqnarray} \nonumber
  && \int_0^T \int_{\Omega(t)} \alpha(y) \nabla_y p(y,t) \cdot
     \nabla_y q(y,t) \, dy \, dt \\ \label{LS ST 1}
  && \hspace*{5mm} + \int_0^T \int_{\Omega(t)} \Big[ \partial_t u(y,t) +
     \mbox{\boldmath $v$}(y,t) \cdot \nabla_y u(y,t) \Big] q(y,t) \, dy \, dt
  \\
  && \hspace*{1cm} \nonumber
     + \int_0^T \int_{\Omega(t)} \alpha(y) \nabla_y u(y,t) \cdot
     \nabla_y q(y,t) \, dy \, dt \, = \,
     \int_0^T \int_{\Omega(t)} f(y,t) q(y,t) \, dy \, dt
\end{eqnarray}
is satisfied for all $q \in Y$, and
\begin{eqnarray} \label{LS ST 2}
  && \int_0^T \int_{\Omega(t)} \Big[ \partial_t z(y,t) +
     \mbox{\boldmath $v$}(y,t) \cdot \nabla_y z(y,t) \Big] p(y,t) \, dy \, dt
  \\
  && \hspace*{4cm} \nonumber
     + \int_0^T \int_{\Omega(t)} \alpha(y) \nabla_y z(y,t) \cdot
     \nabla_y p(y,t) \, dy \, dt \, = \, 0
\end{eqnarray}
is satisfied for all $z \in X_v$. Unique solvability of the mixed
variational formulation \eqref{LS ST 1} and \eqref{LS ST 2} follows
again as in \cite{KLS:2023}, in particular when using the
inf-sup stability condition \eqref{inf-sup time}, and the ellipticity
of the bilinear form $a(q,q) = \| q \|_Y^2$ for all $q \in Y$.
Since $B : X_v \to Y^*$ is surjective, $p \equiv 0$ follows.

For the space-time finite element discretization of the mixed variational
formulation \eqref{LS ST 1} and \eqref{LS ST 2} we use the conforming
finite element space $X_H = \mbox{span} \{ \varphi_k \}_{k=1}^M \subset X_v$ 
of piecewise linear continuous basis functions $\varphi_k$ which are
defined with respect to some admissible locally quasi-uniform decomposition
of the space-time domain $Q$ into simplicial space-time finite elements
$\tau_\ell$ of mesh size $H_\ell$. In addition, we introduce a space-time
finite element space
$Y_h = \mbox{span} \{ \phi_i \}_{i=1}^{\overline{M}} \subset Y$
of piecewise linear continuous basis functions $\phi_i$ which are defined
with respect to a possibly refined decomposition of $Q$ into simplicial
finite elements $\overline{\tau}_\ell$ of local mesh size $h_\ell$. As in the
stationary case we may use one additional refinement, i.e.,
$h_\ell = H_\ell/2$ when $\overline{\tau}_\ell \subset \tau_\ell$.

According to \eqref{Definition wu ST} we define $w_{u,h} \in Y_h$ as
unique solution of the variational formulation such that
\begin{eqnarray}\label{Definition wuh ST}
  && \int_0^T \int_{\Omega(t)} \alpha(y) \, \nabla_y w_{u,h}(y,t) \cdot
     \nabla_y q_h(y,t) \, dy \, dt \\
  && \hspace*{4cm} = \,
     \int_0^T \int_{\Omega(t)} \Big[
     \partial_t u(y,t) + \mbox{\boldmath $v$}(y,t) \cdot \nabla_y u(y,t)
     \Big] \, q_h(y,t) \, dy \, dt \nonumber
\end{eqnarray}
is satisfied for all $q_h \in Y_h$. Hence we can define the discrete norm
\[
  \| u \|_{X_v,h} := \sqrt{\| u \|_Y^2 + \| w_{u,h} \|_Y^2} \leq \| u \|_{X_v}
  \quad \mbox{for all} \; u \in X_v,
\]
and as in \eqref{inf-sup time} we can establish the discrete
inf-sup stability condition
\begin{equation}\label{inf-sup time discrete}
  \| u_H \|_{X_v,h} \leq \sup\limits_{0 \neq q_h \in Y_h}
  \frac{b(u_H,q_h)}{\|q_h \|_Y} \quad \mbox{for all} \; u_H \in X_H ,
\end{equation}
which remains true for all $X_H \subseteq Y_h$. Hence we conclude
unique solvability of the space-time finite element discretization of
\eqref{LS ST 1} and \eqref{LS ST 2} to find $(u_H,p_h) \in X_H \times Y_h$
such that
\begin{eqnarray} \nonumber
  && \int_0^T \int_{\Omega(t)} \alpha(y) \nabla_y p_h(y,t) \cdot
     \nabla_y q_h(y,t) \, dy \, dt \\ \label{LS ST 1 FEM}
  && \hspace*{5mm} + \int_0^T \int_{\Omega(t)} \Big[ \partial_t u_H(y,t) +
     \mbox{\boldmath $v$}(y,t) \cdot \nabla_y u_H(y,t) \Big]
     q_h(y,t) \, dy \, dt \\
  && \hspace*{1cm} \nonumber
     + \int_0^T \int_{\Omega(t)} \alpha(y) \nabla_y u_H(y,t) \cdot
     \nabla_y q_h(y,t) \, dy \, dt \, = \,
     \int_0^T \int_{\Omega(t)} f(y,t) q_h(y,t) \, dy \, dt
\end{eqnarray}
is satisfied for all $q_h \in Y_h$, and
\begin{eqnarray} \label{LS ST 2 FEM}
  && \int_0^T \int_{\Omega(t)} \Big[ \partial_t z_H(y,t) +
     \mbox{\boldmath $v$}(y,t) \cdot \nabla_y z_H(y,t) \Big]
     p_h(y,t) \, dy \, dt \\
  && \hspace*{4cm} \nonumber
     + \int_0^T \int_{\Omega(t)} \alpha(y) \nabla_y z_H(y,t) \cdot
     \nabla_y p_h(y,t) \, dy \, dt \, = \, 0
\end{eqnarray}
is satisfied for all $z_H \in X_H$. Moreover, and as in
\cite{GanglGobrialSteinbach:2025,KLS:2023,Steinbach:2015}, we can derive
estimates for the space-time finite element error
$\| u - u_H \|_Y \leq c \, H \, |u|_{H^2(Q)}$ when assuming $u \in H^2(Q)$.
While for $Y_h = X_H$ we have $p_h \equiv 0$, in the more general case
$X_h \subset Y_h$ but $Y_h \neq X_H$ and assuming an saturation condition
such as in  \eqref{saturation}, we can use $\| p_h \|_Y$ as a posteriori
error indicator for $\| u - u_H \|_{X_{v,h}}$ to drive an adaptive scheme.

\begin{example}
  As a first example, we consider the one dimensional domain 
  $\Omega = (0,1)$, and the time horizon $T = 1$, i.e., $Q=(0,1)^2$. Further,
  we consider $\alpha(x)=\varepsilon = 10^{-2}$ and $\mbox{\boldmath $v$} = 1$.
  As exact solution, and similar as in {\rm \cite{CangianiGeorgoulis:2014}},
  we choose the smooth function
  \begin{equation}
    u(x,t):=\left(1-e^{-t/\varepsilon}\right)
      \left( \frac{e^{(x-1)/\varepsilon}-1}{e^{-1/\varepsilon}-1}
    + x -1 \right),
  \label{eq:uexact_bndlayer}
  \end{equation}
  and we compute $f = \partial_t u - \varepsilon \Delta_x u
  + \mbox{\boldmath $v$} \cdot \nabla_x u$ accordingly. The smooth function
  \eqref{eq:uexact_bndlayer} exhibits a spatial (at $x=1$) and a temporal
  (at $t=0$) boundary layer. The numerical results for both a uniform and
  an adaptive refinement strategy are shown in
  Fig.~\ref{abb:ex_bndlayer1d1d_error_plot}.
  We observe a rate of $\mathcal{O}(H)$ for the error in the energy norm
  and $\mathcal{O}(H^2)$ for the $L^2$ error, as expected.
  In Fig.~\ref{abb:ex_bndlayer1d1d_plot_ph_u-uh} we present a comparison
  between the errors $\| u - u_H \|_Y$, $\| u - u_H \|_{X_v,h}$ and the 
  error estimator $\eta_H = \| p_h \|_Y$.
  One can see that the error indicator is effective and that the error 
  in the norm $\|\cdot\|_{X_v,h}$ is mainly driven by the spatial part of
  the norm. Finally, in Fig.~\ref{abb:ex_bndlayer_resultsadaptive} we
  present the related finite element mesh and the numerical solution $u_H$.
  
  \begin{figure}[!htbp]
    \centering
    \begin{subfigure}[b]{0.44\textwidth}
      \begin{tikzpicture}[scale = 0.8]
        \begin{axis}[
          xmode = log,
          ymode = log,
          xlabel= $\widetilde{M}$,
          ylabel= errors,
          legend style={font=\tiny}, legend pos = south west]
          \addplot[mark = square*,blue] table [col sep=&, y=errL2_u, x=nv]
          {tables/SimData_LSHeat_advdiff1d1d_bndlayer_eps1e-02_beta1_adaptive_EpsInNorm.dat};
          \addlegendentry{$\norm{u - u_H}_{L^2(Q)}$,a}
        
          \addplot[mark = square*,red] table [col sep=&, y=errH1_u, x=nv]
          {tables/SimData_LSHeat_advdiff1d1d_bndlayer_eps1e-02_beta1_adaptive_EpsInNorm.dat};
          \addlegendentry{$\norm{u-u_H}_{Y}$,a}
          
          \addplot[mark = *,blue, dashed] table [col sep=&, y=errL2_u, x=nv]
          {tables/SimData_LSHeat_advdiff1d1d_bndlayer_eps1e-02_beta1_uniform_EpsInNorm.dat};
          \addlegendentry{$\norm{u - u_H}_{L^2(Q)}$,u}
          
          \addplot[mark = *,red, dashed] table [col sep=&, y=errH1_u, x=nv]
          {tables/SimData_LSHeat_advdiff1d1d_bndlayer_eps1e-02_beta1_uniform_EpsInNorm.dat};
          \addlegendentry{$\norm{u-u_H}_{Y}$,u}
          
          \addplot[
          domain = 10000:200000,
          samples = 10,
          thin,
          black,
          ] {50*x^(-1)};
          \addlegendentry{$\widetilde{M}^{-1}\sim H^{2}$}
          
          \addplot[
          domain = 10000:200000,
          samples = 10,
          thin,
          black,
          dashed,
        ] {10*x^(-1/2)};
        \addlegendentry{$\widetilde{M}^{-1/2}\sim H$}
        \end{axis}
      \end{tikzpicture}
      \caption{Errors $\| u - u_H \|_{Y}$ and
            $\| u - u_H \|_{L^2(Q)}$ for uniform and adaptive
            refinement strategies.}
      \label{abb:ex_bndlayer1d1d_error_plot}
    \end{subfigure}
    ~
    \begin{subfigure}[b]{0.44\textwidth}
      \begin{tikzpicture}[scale = 0.8]
        \begin{axis}[
          xmode = log,
          ymode = log,
          xlabel= $\widetilde{M}$,
          ylabel= errors,
          legend style={font=\tiny}, legend pos = north east]
          \addplot[mark = x,blue] table [col sep=&, y=L2_dxph, x=nv]
          {tables/SimData_LSHeat_advdiff1d1d_bndlayer_eps1e-02_beta1_adaptive_EpsInNorm.dat};
          \addlegendentry{$\eta_H$}
          
          \addplot[mark = o,green!60!black] table [col sep=&, y=errXH_u_hat,
          x=nv]{tables/SimData_LSHeat_advdiff1d1d_bndlayer_eps1e-02_beta1_adaptive_EpsInNorm.dat};
          \addlegendentry{$\| u-u_H \|_{X_v,h}$}
        
          \addplot[mark = x,red] table [col sep=&, y=errH1_u, x=nv]
          {tables/SimData_LSHeat_advdiff1d1d_bndlayer_eps1e-02_beta1_adaptive_EpsInNorm.dat};
          \addlegendentry{$\norm{u-u_H}_{Y}$}
        \end{axis}
      \end{tikzpicture}
      \caption{Comparison between estimator and true error for an adaptive refinement strategy.}
      \label{abb:ex_bndlayer1d1d_plot_ph_u-uh}
    \end{subfigure}
    
    \caption{Convergence results in the case of a smooth solution for a nonstationary
          convection-diffusion equation.}
  \end{figure}
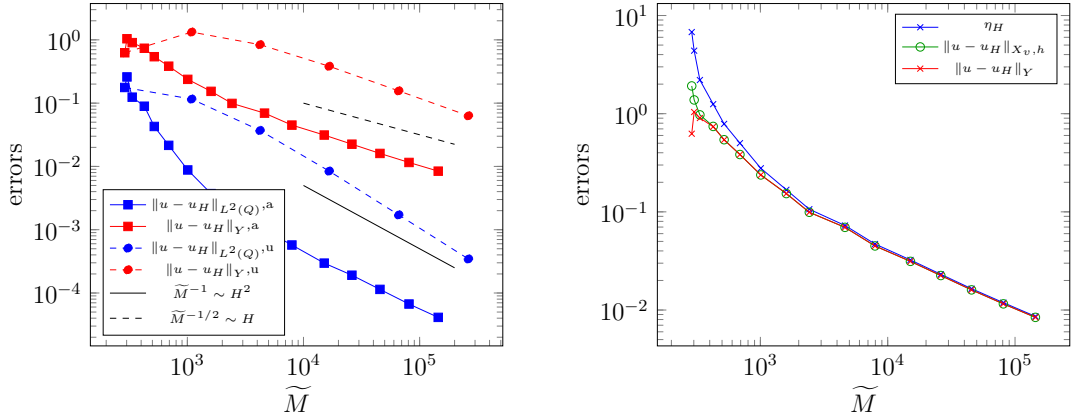

  \begin{figure}[!htbp]
    \centering
    \begin{subfigure}[b]{0.47\textwidth}
      \includegraphics[trim = 350px 150px 350px 150px,clip, width=\textwidth]
          {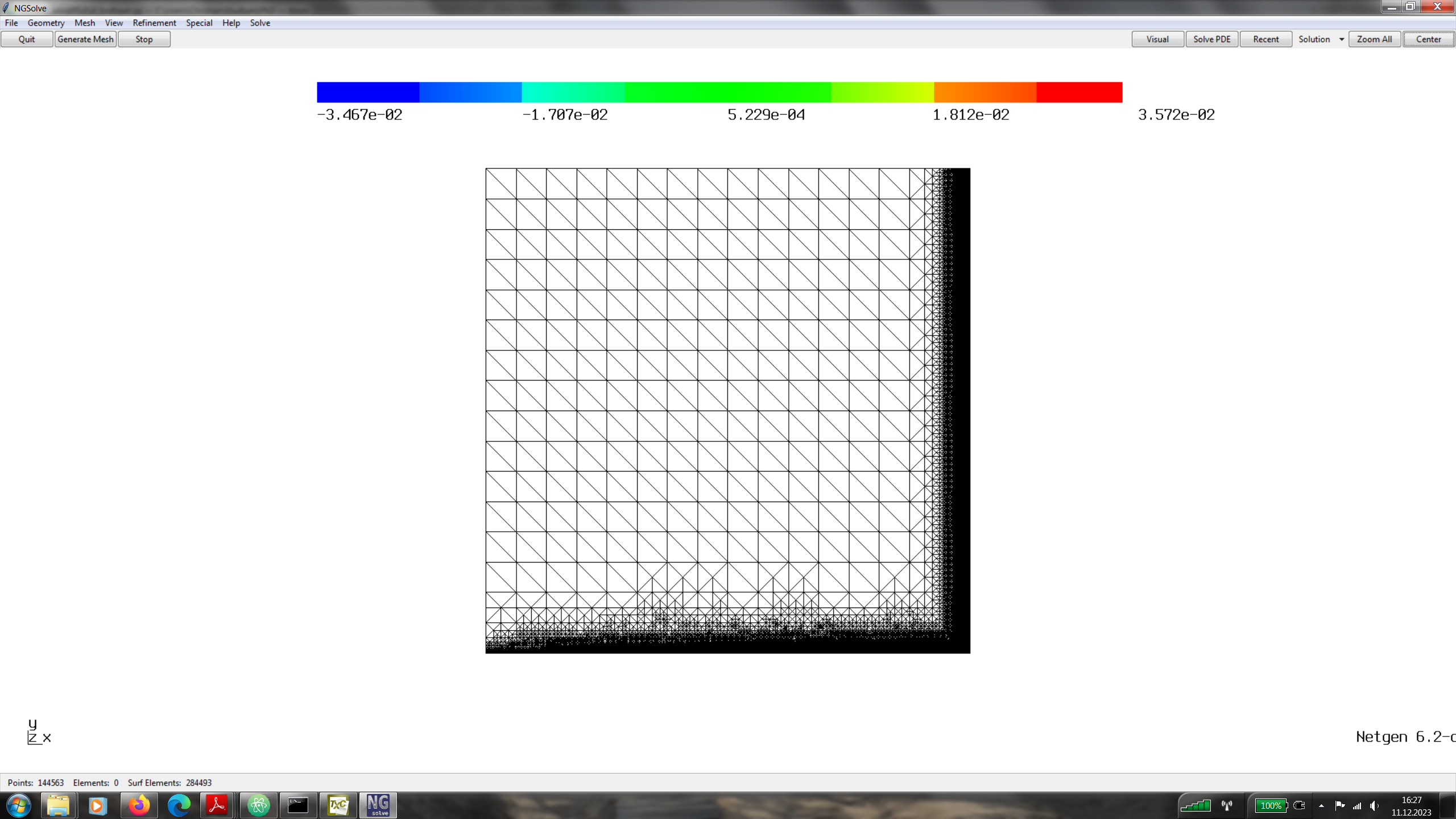}
      \caption{Adaptive mesh on $L=15$, $144563$ dofs}
    \end{subfigure}
    ~
    \begin{subfigure}[b]{0.47\textwidth}
      \includegraphics[trim = 350px 150px 350px 150px,clip, width=\textwidth]
          {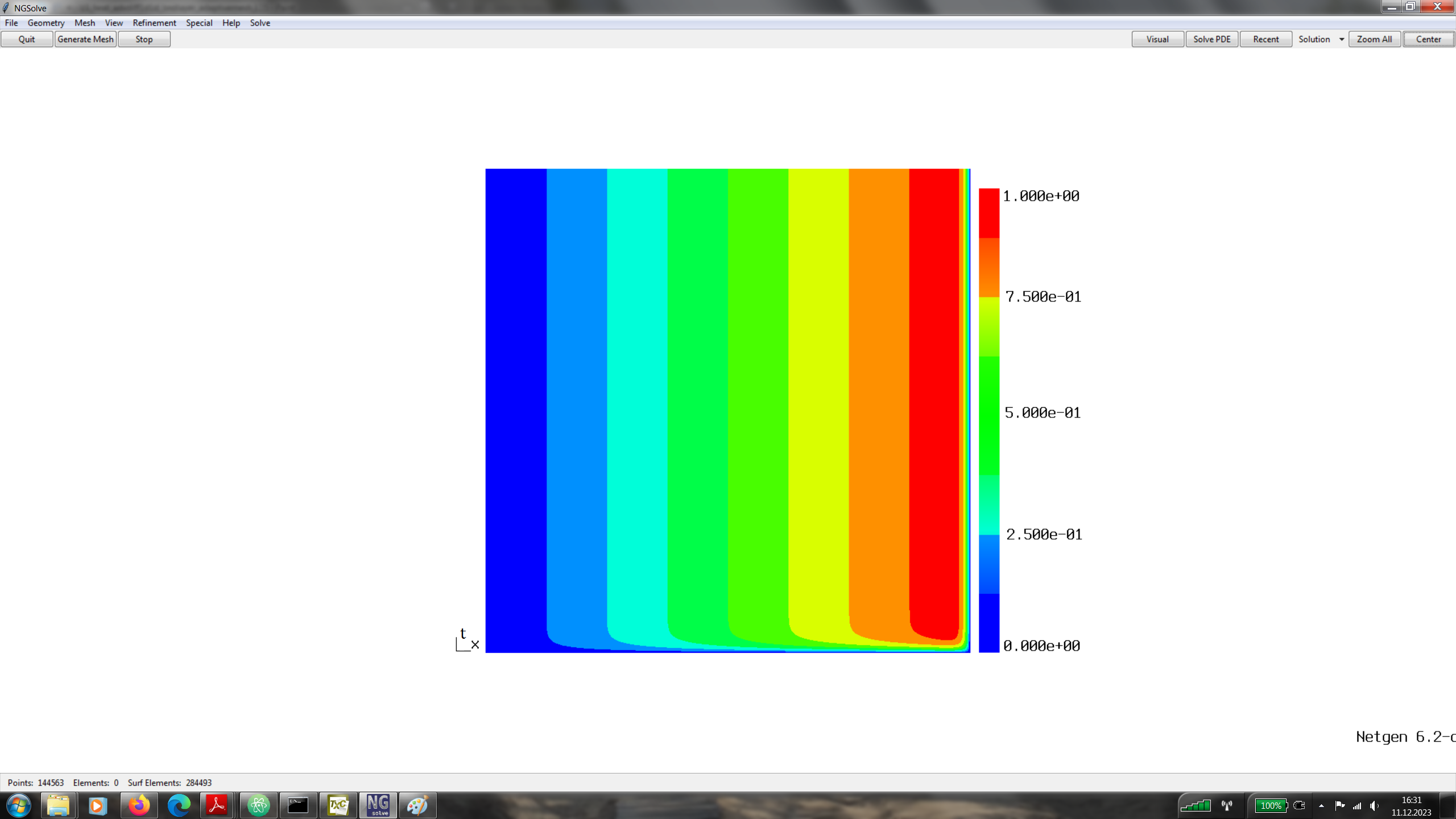}
      \caption{Solution $u_H$ on adaptive mesh}
    \end{subfigure}
    \caption{Simulation results for the adaptive refinement process.}%
    \label{abb:ex_bndlayer_resultsadaptive}%
  \end{figure}

\end{example}

\begin{example}
  As a second example we consider the two-dimensional domain
  $\Omega = (0,1)^2$ and the time horizon $T=1$, i.e., $Q=(0,1)^3$. As
  initial state $u(x,0)=u_0(x)$, $x \in \Omega$, we consider similar as
  in {\rm \cite{LOSSepulveda:2023}} the function 
  \begin{equation*}
    u_0(x):=\psi(10 \, \|x-x_0\|_2),\quad 
    \psi(r):=\begin{cases}
              (1-r^2)^2, & \mbox{for} ~ r\leq 1, \\[1mm]
              0, & \mbox{for} ~ r>1,
            \end{cases} \quad
            x_0 = \left( \begin{array}{c} 0.5 \\ 0.5 \end{array} \right).
  \end{equation*}
  We compute numerical solutions to \eqref{IBVP Advection Diffusion} for the
  velocity field $\mbox{\boldmath $v$} = (0,1)^\top$ and without a source term,
  i.e. $f \equiv 0$. Furthermore, 
  $\alpha(x) = \varepsilon \in \left\{10^{-3}, 10^{-5}, 10^{-6}\right\}$.
  The results for a mesh with $32 \times 32 \times 32$ elements ($35937$ dofs)
  can be seen in Fig. \ref{abb:LS_Heat_advdiff_ex2_ComparisonDirectvsLS}.
  In the top row one can see the numerical solution $u_H$ computed by solving
  \eqref{IBVP VF} with the space-time finite element method described in
  {\rm \cite{Steinbach:2015}}. This leads to oscillations in the solution as
  the mesh size is not sufficiently small. In the bottom row one can see the
  solution using the developed least-squares formulation with piecewise linear
  trial and piecewise quadratic test functions. This formulation leads to
  stable results.
  \begin{figure}
    \centering
    \begin{subfigure}[b]{0.3\textwidth}
      \includegraphics[width=\textwidth]{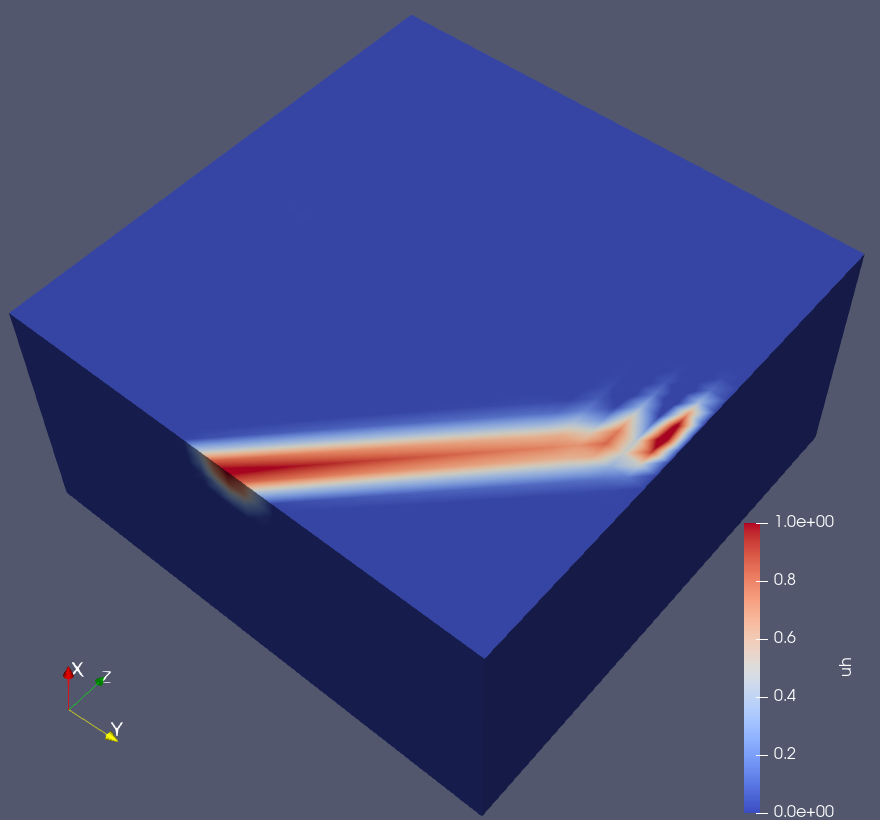}
      \caption{$\varepsilon = 10^{-3}$}
    \end{subfigure}
    ~
    \begin{subfigure}[b]{0.3\textwidth}
      \includegraphics[width=\textwidth]{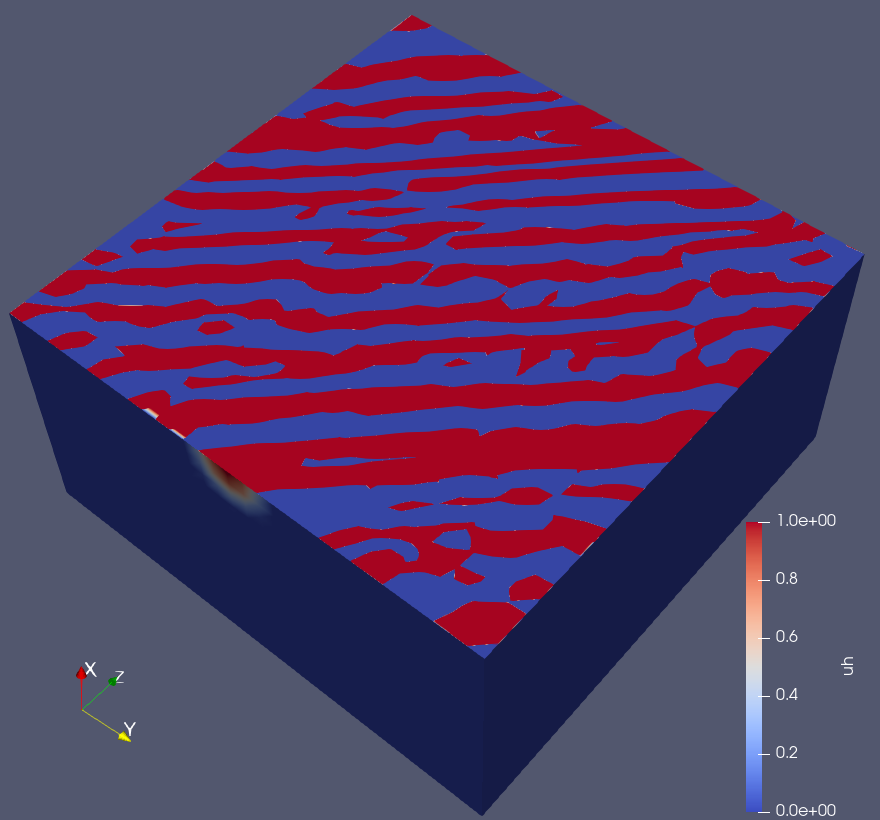}
      \caption{$\varepsilon = 10^{-5}$}
    \end{subfigure}
    ~
    \begin{subfigure}[b]{0.3\textwidth}
      \includegraphics[width=\textwidth]{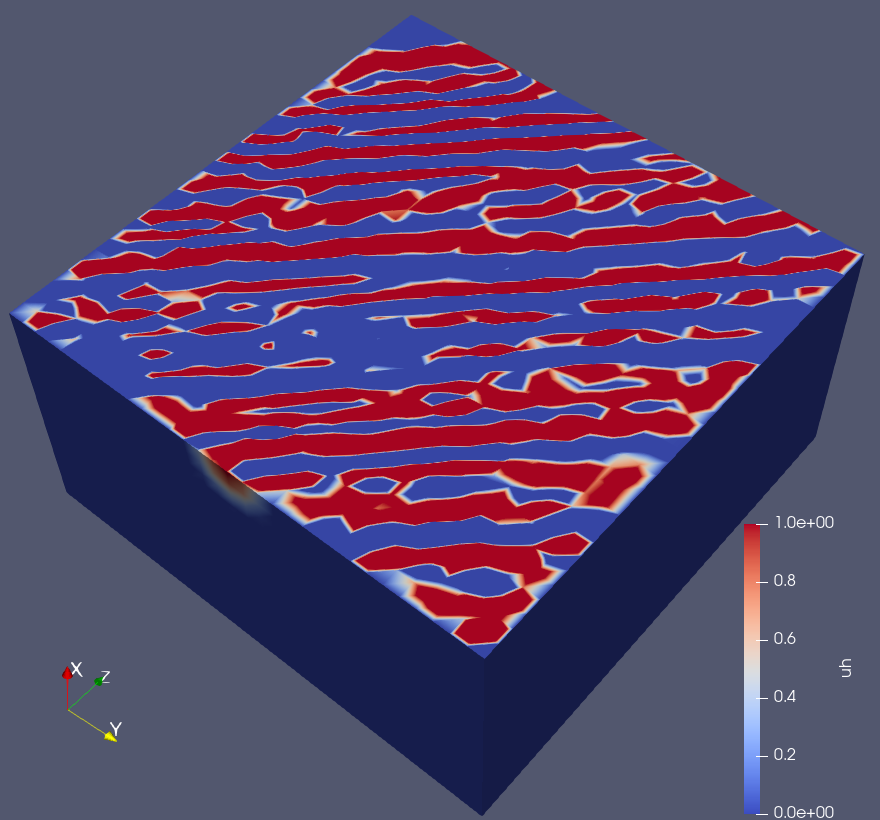}
      \caption{$\varepsilon = 10^{-6}$}
    \end{subfigure}
    \\
    \begin{subfigure}[b]{0.3\textwidth}
      \includegraphics[width=\textwidth]{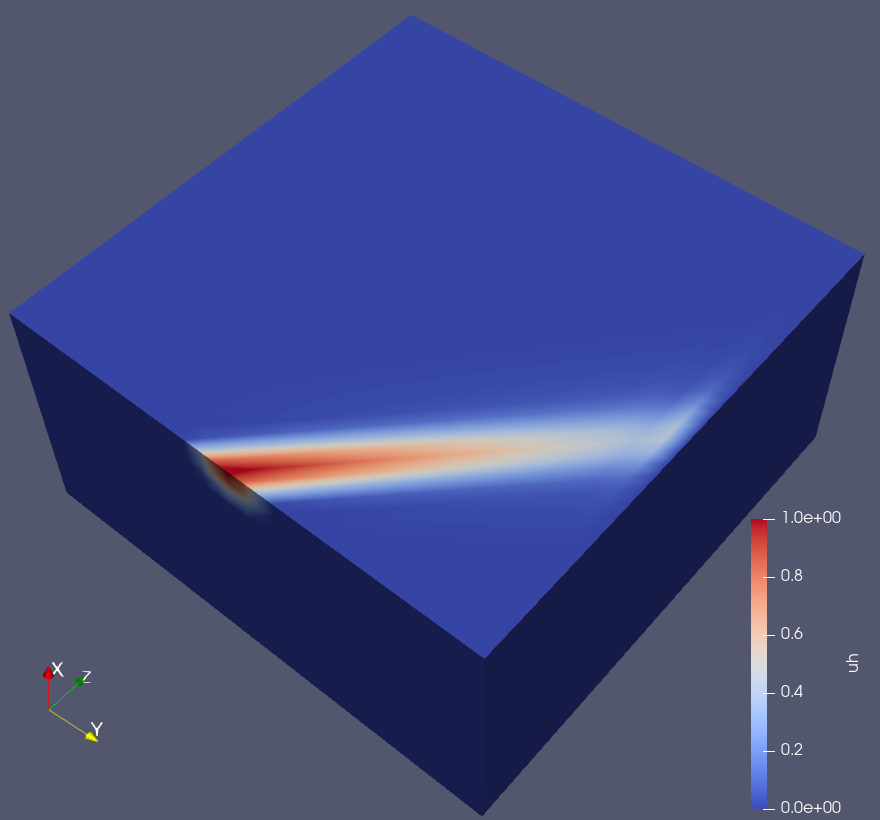}
      \caption{$\varepsilon = 10^{-3}$}
    \end{subfigure}
    ~
    \begin{subfigure}[b]{0.3\textwidth}
      \includegraphics[width=\textwidth]{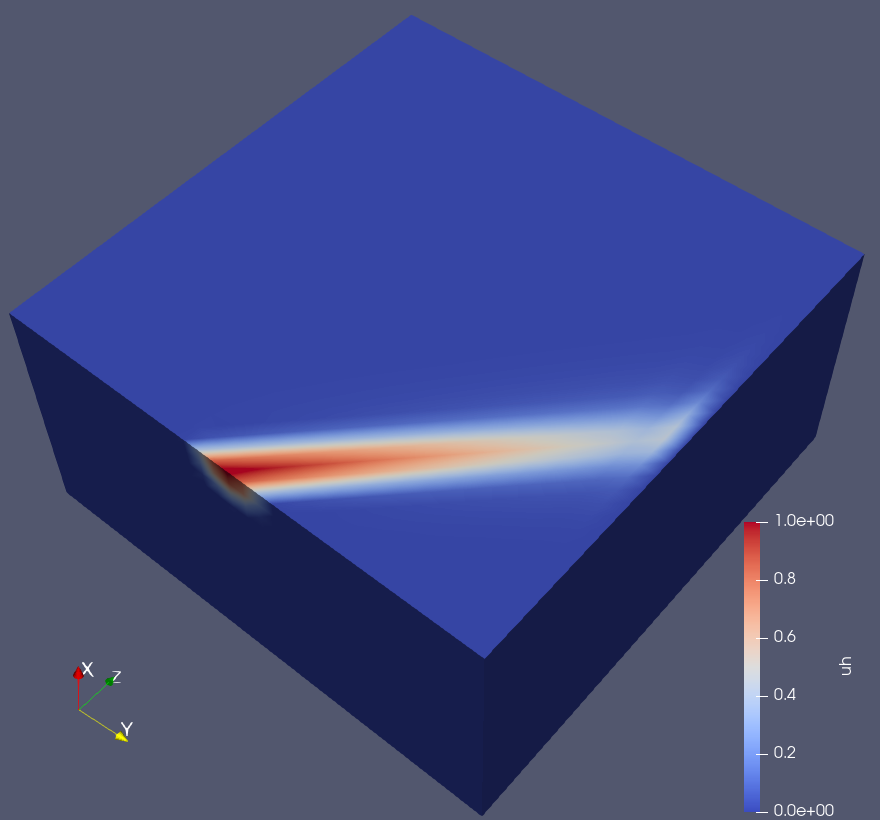}
      \caption{$\varepsilon = 10^{-5}$}
    \end{subfigure}
    ~
    \begin{subfigure}[b]{0.3\textwidth}
      \includegraphics[width=\textwidth]{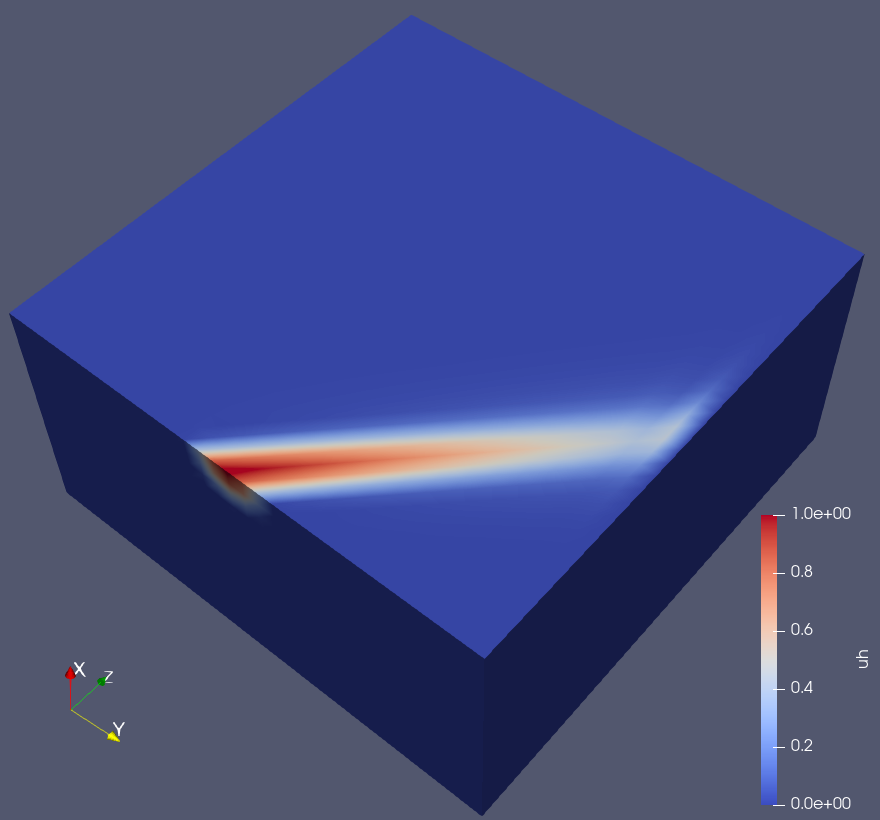}
      \caption{$\varepsilon = 10^{-6}$}
    \end{subfigure}
    
  \caption{Numerical results for {\boldmath{$v$}} $= (0,1)^T$ on a mesh with $32\times 32 \times 32$
      elements. Top: no stabilization via direct formulation \cite{Steinbach:2015}, 
      bottom: stabilization via developed least-squares 
      formulation.}
  \label{abb:LS_Heat_advdiff_ex2_ComparisonDirectvsLS}
  \end{figure}
  In a further step we use the inbuilt error estimator to drive an adaptive
  refinement scheme for the parameters $\varepsilon = 10^{-3}$,
  $\mbox{\boldmath $v$} = (0,0.3)^\top$, and $\varepsilon = 10^{-6}$, 
  $\mbox{\boldmath $v$} = (0,1)^\top$. In Fig. \ref{abb:advdiff_ex2_estimator}
  the convergence rate of the error estimator in case of a uniform and an
  adaptive refinement strategy for both sets of parameters is depicted.
  In the case $\varepsilon = 10^{-3}$, $\mbox{\boldmath $v$} = (0,0.3)^\top$,
  we observe a linear rate $\mathcal{O}(H)$ for both refinement strategies.
  In the case $\varepsilon = 10^{-6}$, $\mbox{\boldmath $v$} = (0,1)^\top$
  we observe a reduced rate of $\mathcal{O}(H^{0.4})$ in the uniform case. 
  However, we can recover the full rate of $\mathcal{O}(H)$ in the adaptive
  case. The  obtained adaptive meshes as well as the corresponding numerical 
  solutions are depicted in
  Fig. \ref{abb:advdiff_ex2_ResultsAfterAdaptiveRefinement_eps1e-06} and 
  \ref{abb:advdiff_ex2_ResultsAfterAdaptiveRefinement_eps1e-03}. Note that
  we obtain a mesh which is fully unstructured in space and time.
  \begin{figure}[!htbp]
    \centering
    \begin{subfigure}[b]{0.45\textwidth}
    \begin{tikzpicture}[scale=0.7]
          \begin{axis}[
            xmode = log,
            ymode = log,
            xlabel= $\widetilde{M}$,
            ylabel= errors,
            xmin = 1000,
            legend style={font=\tiny}, legend pos = south west]
            \addplot[mark = *,red] table [col sep=&, y=L2_dxph, x=nv]
            {tables/SimData_LSHeat_advdiff2d1d_eps1e-03_beta3e-01_uniform_MinresWithJacboi_epsInNorm.dat};
            \addlegendentry{$\eta_H$ uniform}
            \addplot[mark = *,blue] table [col sep=&, y=L2_dxph, x=nv]
            {tables/SimData_LSHeat_advdiff2d1d_eps1e-03_beta3e-01_adaptive_pardiso_epsInNorm.dat};
            \addlegendentry{$\eta_H$ adaptive}			
            \addplot[
              domain = 80000:900000,
              samples = 10,
              thin,
              black,
            ] {0.3*x^(-1/3)}; 
            \addlegendentry{$\widetilde{M}^{-1/3}\sim H$}
          \end{axis}
        \end{tikzpicture}
      \caption{$\varepsilon = 10^{-3}$, {\boldmath{$v$}}$ = (0,0.3)^\top$}
    \end{subfigure}
    ~
    \begin{subfigure}[b]{0.45\textwidth}
    \begin{tikzpicture}[scale=0.7]
          \begin{axis}[
            xmode = log,
            ymode = log,
            xlabel= $\widetilde{M}$,
            ylabel= errors,
            xmin = 1000,
            legend style={font=\tiny}, legend pos = south west]
            \addplot[mark = *,red] table [col sep=&, y=L2_dxph, x=nv]
            {tables/SimData_LSHeat_advdiff2d1d_eps1e-06_beta1_uniform_MinresWithJacboi_epsInNorm.dat};
            \addlegendentry{$\eta_H$ uniform}
            \addplot[mark = *,blue] table [col sep=&, y=L2_dxph, x=nv]
            {tables/SimData_LSHeat_advdiff2d1d_eps1e-06_beta1_adaptive_pardiso_epsInNorm.dat};
            \addlegendentry{$\eta_H$ adaptive}			
            \addplot[
              domain = 20000:900000,
              samples = 10,
              thin,
              black,
            ] {43*x^(-1/3)}; 
            \addlegendentry{$\widetilde{M}^{-1/3}\sim H$}
            \addplot[
              domain = 50000:2000000,
              samples = 10,
              thin,
              black,
              dashed,
            ] {17*x^(-0.4/3)}; 
            \addlegendentry{$\widetilde{M}^{-0.4/3}\sim H^{0.4}$}
          \end{axis}
        \end{tikzpicture}
        \caption{$\varepsilon = 10^{-6}$, {\boldmath{$v$}}$ = (0,1)^\top$}
      \end{subfigure}
    \caption{Error estimator $\eta_H = \| p_h \|_{Y}$ in case of an adaptive
      and uniform refinement strategy for the second example.}%
    \label{abb:advdiff_ex2_estimator}%
  \end{figure}
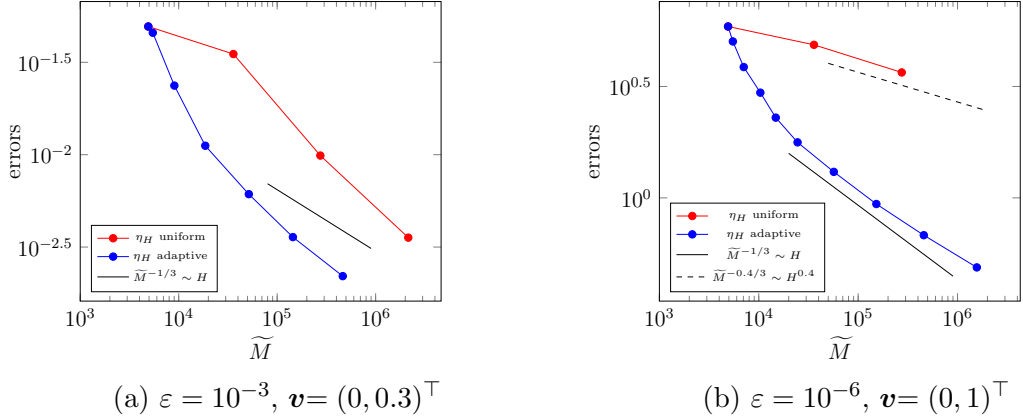

  \begin{figure}
    \centering
    \begin{subfigure}[b]{0.47\textwidth}
      \includegraphics[width=0.95\textwidth]
          {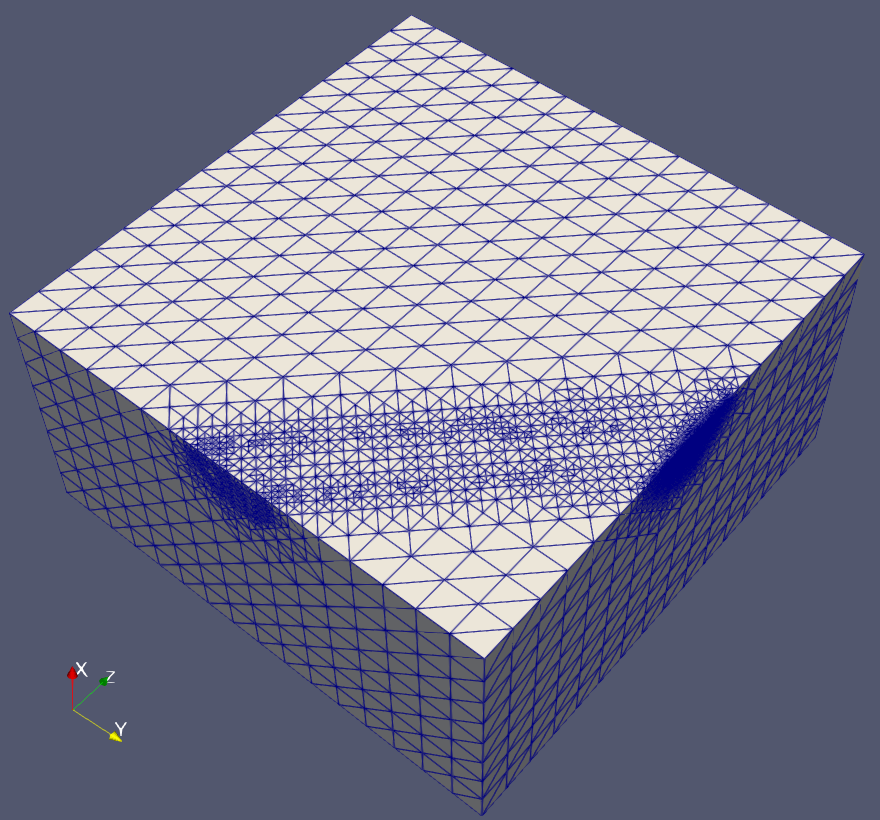}
      \caption{Adaptive mesh on $L=7$, $152513$ dofs}
    \end{subfigure}
    ~
    \begin{subfigure}[b]{0.47\textwidth}
      \includegraphics[width=0.95\textwidth]
          {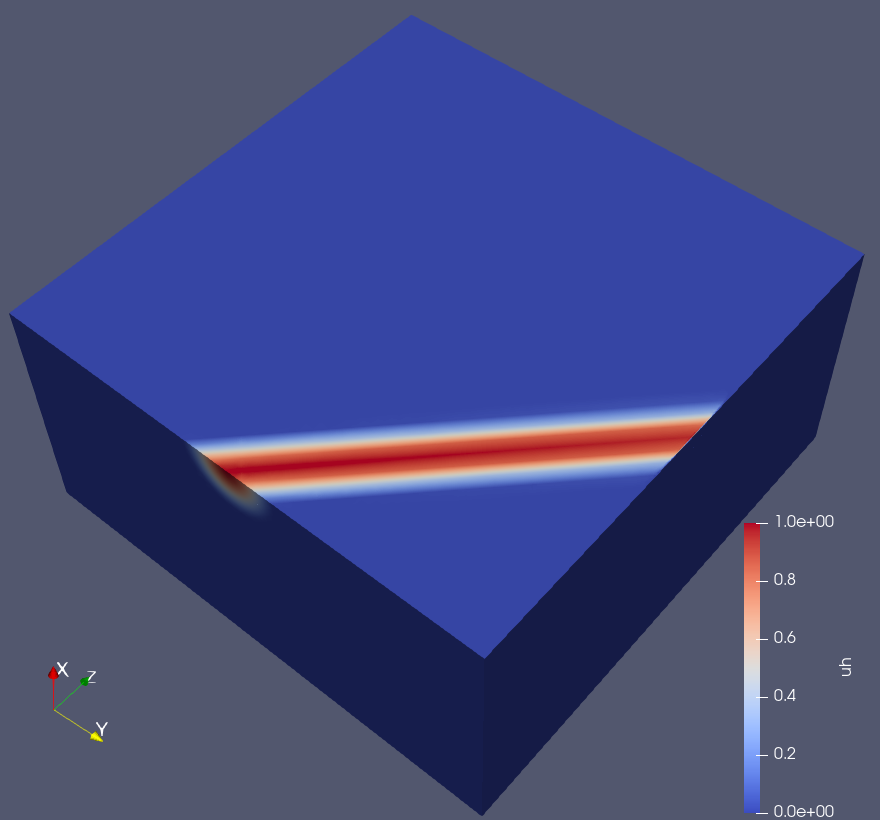}
      \caption{Solution $u_H$ on the adaptive mesh}
    \end{subfigure}
    \caption{Obtained results for $\varepsilon = 10^{-6}$ and {\boldmath{$v$}}$ = (0,1)^\top$ after 
      the adaptive refinement process.}%
    \label{abb:advdiff_ex2_ResultsAfterAdaptiveRefinement_eps1e-06}%
  \end{figure}

  \begin{figure}
    \centering
    \begin{subfigure}[b]{0.47\textwidth}
      \includegraphics[width=0.95\textwidth]
          {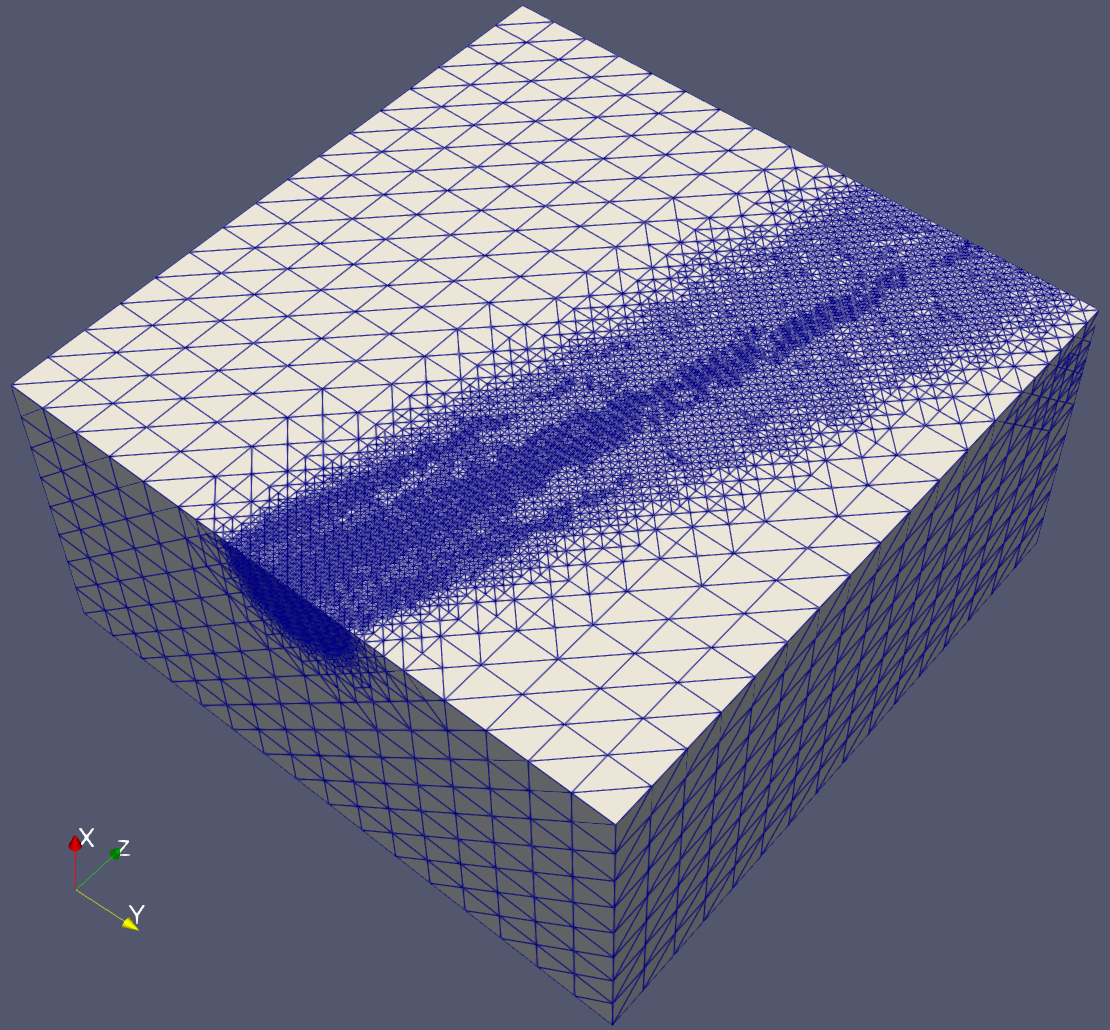}
      \caption{Adaptive mesh on $L=6$, $463696$ dofs}
    \end{subfigure}
    ~
    \begin{subfigure}[b]{0.47\textwidth}
      \includegraphics[width=0.95\textwidth]
          {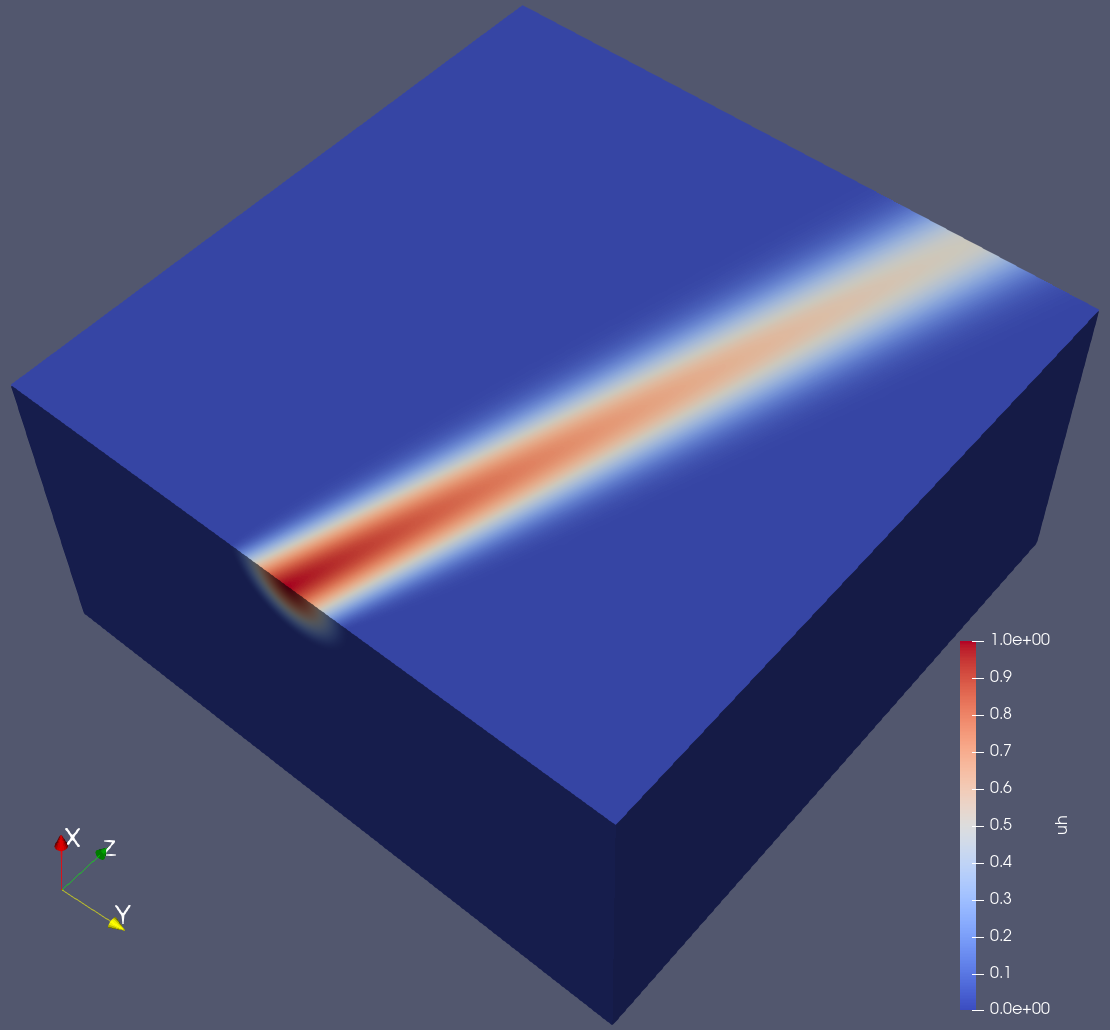}
      \caption{Solution $u_H$ on adaptive the mesh}
    \end{subfigure}
    \caption{Obtained results for $\varepsilon = 10^{-3}$ and {\boldmath{$v$}}$ = (0,0.3)^\top$ after 
      the adaptive refinement process.}%
    \label{abb:advdiff_ex2_ResultsAfterAdaptiveRefinement_eps1e-03}%
  \end{figure}

\end{example}

\begin{example}
  As a third example, we consider again the unit cube in the space-time domain,
  i.e., $Q=(0,1)^3$. We choose $u_0=0$, $\alpha(x) = \varepsilon = 10^{-2}$
  and the source term to be $f = 1$. The velocity field is 
  a time dependent function with $\mbox{\boldmath $v$}(y,t) =
  (\sin(2\pi t), \cos(2\pi t))^\top$ for $(y,t) \in Q$. Thus, the solution $u$
  has a boundary layer whose location depends on time. Note that a similar
  example is considered in {\rm \cite[Example~4]{CangianiGeorgoulis:2014}}.
  In Fig. \ref{abb:advdiff_movingbndlayer_estimator} a comparison of the
  error estimator in case of an adaptive and uniform refinement strategy is
  depicted. We observe a convergence rate of $\mathcal{O}(H)$ in both cases.
  The generated grids in the adaptive case at different fixed times can be
  seen in Fig. \ref{abb:ResultsMovingBndlayerAdaptive}. The circular movement
  of the boundary layer in time is visible.
  \begin{figure}
    \centering
    \begin{tikzpicture}[scale=0.8]
      \begin{axis}[
        xmode = log,
        ymode = log,
        xlabel= $\widetilde{M}$,
        ylabel= errors,
        xmin = 1000,
        legend style={font=\tiny}, legend pos = south west]
        \addplot[mark = *,red] table [col sep=&, y=L2_dxph, x=nv]
        {tables/SimData_LSHeat_advdiff2d1d_movingbndlayer_uniform_pardiso_epsInNorm.dat};
        \addlegendentry{$\eta_H$ uniform}
        \addplot[mark = *,blue] table [col sep=&, y=L2_dxph, x=nv]
        {tables/SimData_LSHeat_advdiff2d1d_movingbndlayer_adaptive_pardiso_epsInNorm.dat};
        \addlegendentry{$\eta_H$ adaptive}			
        \addplot[
          domain = 80000:900000,
          samples = 10,
          thin,
          black,
          ] {5*x^(-1/3)};
        \addlegendentry{$\widetilde{M}^{-1/3}\sim H$}
      \end{axis}
    \end{tikzpicture}
  \caption{Error estimator $\eta_H = \| p_h \|_Y$ in case of an adaptive
      and uniform refinement strategy for the third example.}%
  \label{abb:advdiff_movingbndlayer_estimator}%
  \end{figure}
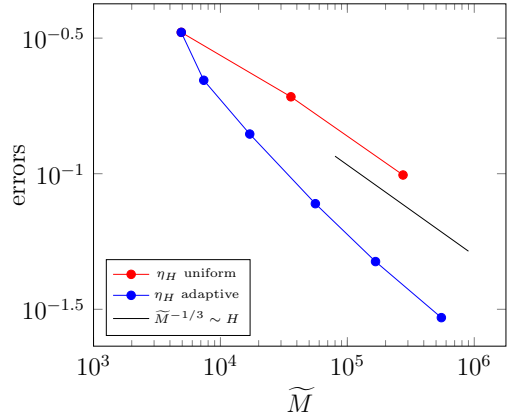
  \begin{figure}
    \centering
    \begin{subfigure}[b]{0.47\textwidth}
          \centering
      \includegraphics[width=0.95\textwidth]
          {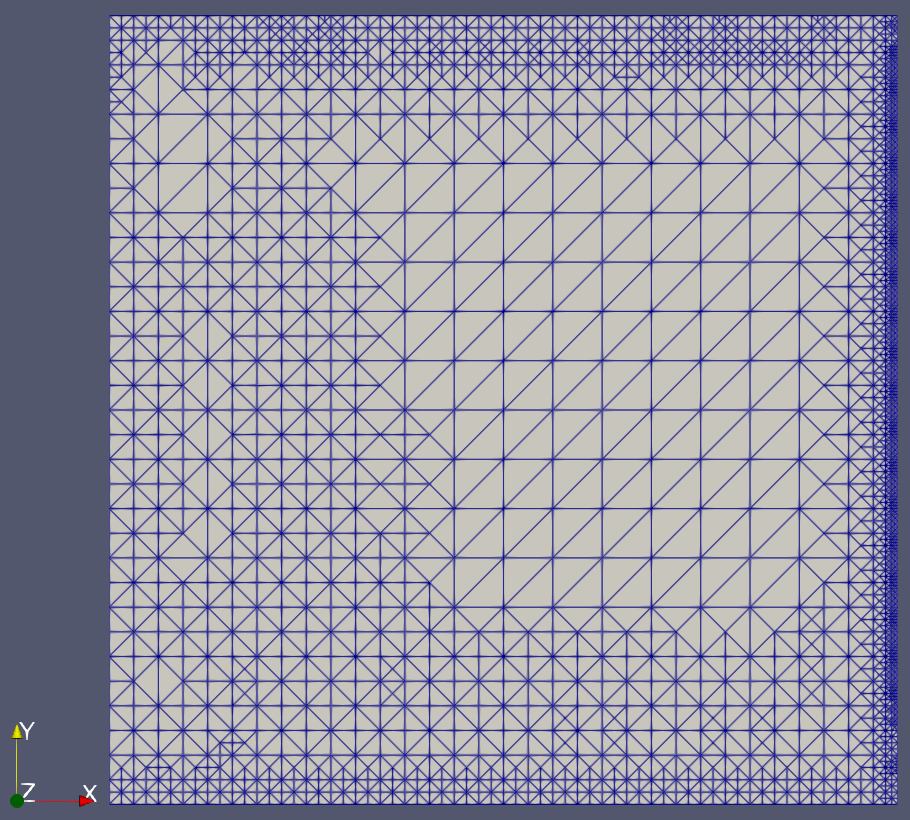}
      \caption{$t=0.25$}
    \end{subfigure}
    ~
    \begin{subfigure}[b]{0.47\textwidth}
          \centering
      \includegraphics[width=0.95\textwidth]
          {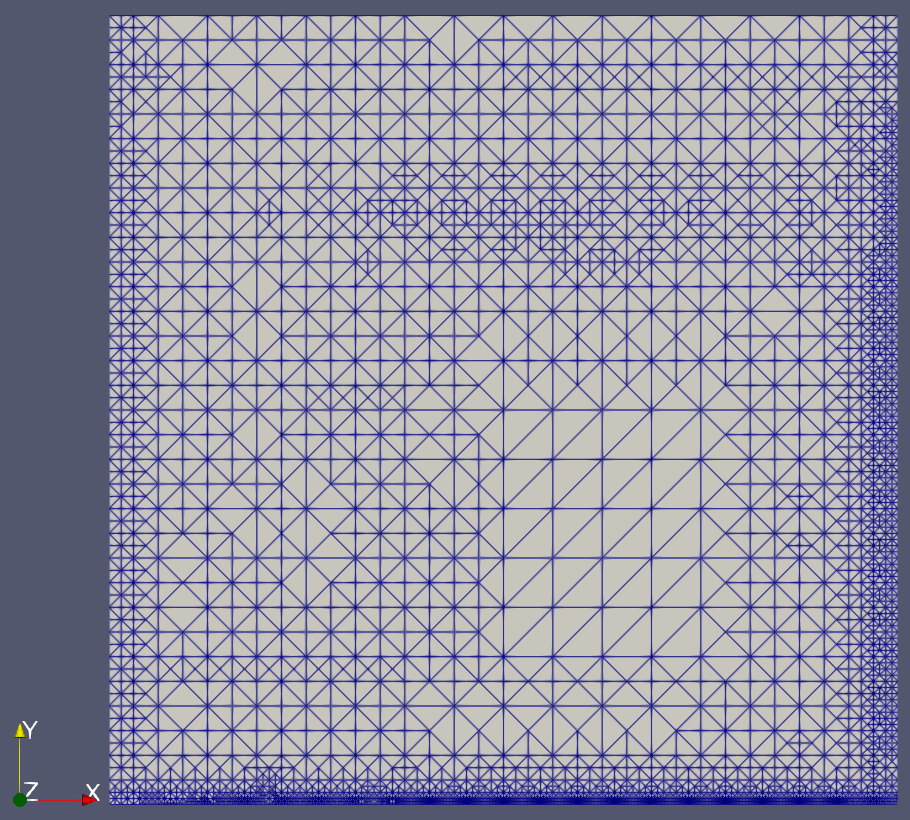}
      \caption{$t=0.5$}
    \end{subfigure}
    \\
    \begin{subfigure}[b]{0.47\textwidth}
          \centering
      \includegraphics[width=0.95\textwidth]
          {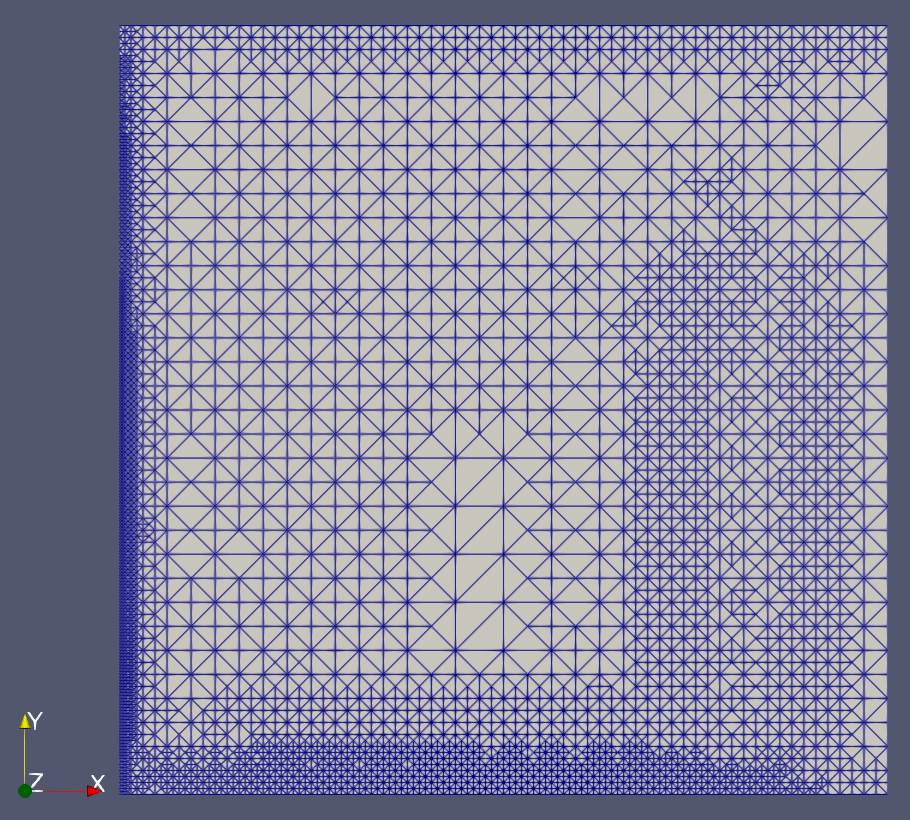}
      \caption{$t=0.75$}
    \end{subfigure}
    ~
    \begin{subfigure}[b]{0.47\textwidth}
          \centering
      \includegraphics[width=0.95\textwidth]
          {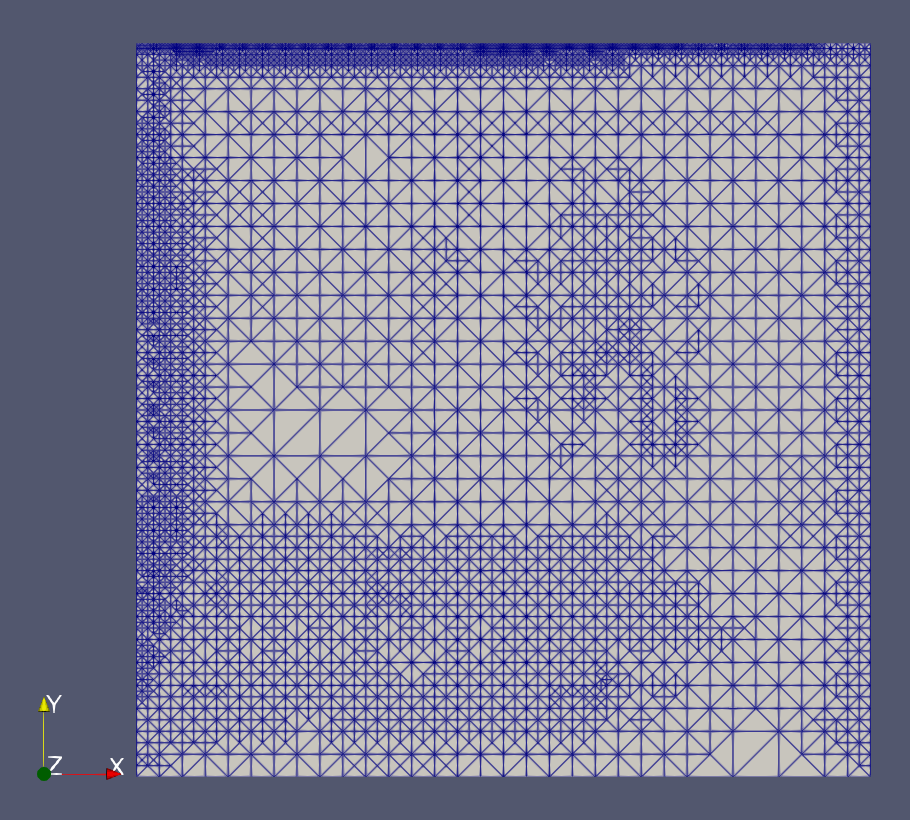}
      \caption{$t=1$}
    \end{subfigure}
    \caption{Generated mesh on refinement level $L = 5$ at different times $t$
          in case of the time dependent velocity field 
          {\boldmath{$v$}}$(x,t) = \left(\sin(2\pi t), \cos(2\pi t)\right)^\top$.}%
    \label{abb:ResultsMovingBndlayerAdaptive}%
  \end{figure}
\end{example}

\pagebreak

\noindent
\textbf{Acknowledgement:}
This work is supported by the joint
DFG/FWF Collaborative Research Centre CREATOR
(DFG: Project-ID 492661287/TRR361; FWF: 10.55776/F90) at
TU Darmstadt, TU Graz and JKU Linz.

\end{document}

%% file: tables/SimData_advdiff1d_direct_DifferentVelocities.tex
\begin{tabular}{rcrrcrrcrr} \toprule
    && \multicolumn{2}{c}{$v= 10$} && \multicolumn{2}{c}{$v=100$} && \multicolumn{2}{c}{$v=1000$} \\
    \cmidrule{3-4} \cmidrule{6-7} \cmidrule{9-10}
    $N$ && $\| \nabla(u-u_H) \|_{L^2(0,1)}$ & eoc  
        && $\| \nabla(u-u_H) \|_{L^2(0,1)}$ & eoc  
        && $\| \nabla(u-u_H) \|_{L^2(0,1)}$ & eoc \\
    \midrule  
    2 && 1.567e-01 &  && 5.283e-01 &  && 5.073e+00 &  \\
    4 && 7.554e-02 & 1.053 && 9.943e-02 & 2.410 && 6.511e-01 & 2.962 \\
    8 && 3.759e-02 & 1.007 && 3.877e-02 & 1.359 && 8.967e-02 & 2.860 \\ 
    16 && 1.877e-02 & 1.002 && 1.890e-02 & 1.036 && 2.155e-02 & 2.057 \\
    32 && 9.385e-03 & 1.000 && 9.400e-03 & 1.008 && 9.563e-03 & 1.172 \\ 
    64 && 4.692e-03 & 1.000 && 4.694e-03 & 1.002 && 4.714e-03 & 1.021 \\ 
    128 && 2.346e-03 & 1.000 && 2.346e-03 & 1.000 && 2.349e-03 & 1.005 \\ 
    256 && 1.173e-03 & 1.000 && 1.173e-03 & 1.000 && 1.173e-03 & 1.001 \\ 
    \bottomrule
\end{tabular}

%% file: tables/SimData_advdiff1d_direct_VelocityAsFunctionInx.tex
\begin{tabular}{rcrrcrr} \toprule
    && \multicolumn{2}{c}{$v(x) = \frac{1}{\sqrt{x}}$} && \multicolumn{2}{c}{$v(x) = \frac{1}{x}$}  \\
    \cmidrule{3-4} \cmidrule{6-7}
    $N$ && $\norm{\nabla(u-u_H)}_{L^2(0,1)}$ & eoc 
        && $\norm{\nabla(u-u_H)}_{L^2(0,1)}$ & eoc \\
    \midrule  
    2 && 1.485e-01 & 0.000 && 1.494e-01 & 0.000 \\ 
    4 && 7.487e-02 & 0.988 && 7.503e-02 & 0.994 \\
    8 && 3.751e-02 & 0.997 && 3.754e-02 & 0.999 \\
    16 && 1.876e-02 & 0.999 && 1.877e-02 & 1.000 \\
    32 && 9.383e-03 & 1.000 && 9.384e-03 & 1.000 \\
    64 && 4.692e-03 & 1.000 && 4.692e-03 & 1.000 \\
    128 && 2.346e-03 & 1.000 && 2.346e-03 & 1.000 \\
    256 && 1.173e-03 & 1.000 && 1.173e-03 & 1.000 \\
    \bottomrule
\end{tabular}
 

%% file: tables/SimData_advdiff1d_direct.tex
\begin{tabular}{rrrrrrrrr} \toprule
    \multicolumn{3}{c}{$\varepsilon = 10^{-2}$} & \multicolumn{3}{c}{$\varepsilon=10^{-4}$} & \multicolumn{3}{c}{$\varepsilon=10^{-5}$} \\
    \cmidrule{2-3} \cmidrule{5-6} \cmidrule{8-9}
    $N$ & $\norm{\nabla(u-u_H)}_{L^2(0,1)}$ & eoc & 
    $N$ & $\norm{\nabla(u-u_H)}_{L^2(0,1)}$ & eoc & 
    $N$ & $\norm{\nabla(u-u_H)}_{L^2(0,1)}$ & eoc \\
    \midrule  
    128 & 1.567e+00 &  & 4096 & 4.280e+01 & & 65536 & 9.237e+01 &  \\ 
    256 & 7.938e-01 & 0.981 & 8192 & 2.389e+01 & 0.841 & 131072 & 4.843e+01 & 0.931 \\ 
    512 & 3.982e-01 & 0.995 & 16384 & 1.233e+01 & 0.955 & 262144 & 2.452e+01 & 0.982 \\ 
    1024 & 1.993e-01 & 0.999 & 32768 & 6.213e+00 & 0.988 & 524288 & 1.230e+01 & 0.995 \\ 
    2048 & 9.966e-02 & 1.000 & 65536 & 3.113e+00 & 0.997 & 1048576 & 6.154e+00 & 0.999 \\
    \bottomrule
    \end{tabular}

%% file: tables/SimData_LS_laplace_advdiff2d_nonconstantconvection.tex
\begin{tabular}{rrrrrr}
    \toprule 
    $L$ & $\widetilde{M}$ & lower layer & upper layer & profile & $\osc_{max}(u_H)$ \\
    \midrule 
    0 & 289 & 0.17989000 & 0.06040999 & 0.49811000 & 0.43594821 \\
    5 & 1666 & 0.04909000 & 0.05312000 & 0.37743000 & 0.28357896 \\ 
    10 & 15160 & 0.01785999 & 0.01824999 & 0.35253999 & 0.07359097 \\ 
    13 & 51480 & 0.01422000 & 0.01436000 & 0.34822000  & 1.523e-03 \\
    15 & 148968 & 0.01440000 & 0.01440000 & 0.34825000 & 7.373e-04 \\
    17 & 473779 & 0.01439000 & 0.01439000 & 0.34825000 & 4.033e-04 \\
    \bottomrule
\end{tabular}

%% file: tables/SimData_Hemker.tex
\begin{tabular}{rrrrrr}
    \toprule 
    $L$ & $\widetilde{M}$ & upper layer & $\osc_{max}(u_H)$ \\
    \midrule 
    0 & 438 & 0.8087 & 0.9742 \\
    8 & 12048 & 0.4873 & 0.6315 \\
    12 & 82453 & 0.1603 & 0.1909 \\
    15 & 227728  & 0.0801 & 0.0165 \\
    17 & 676072  & 0.0728 & 1.187e-03 \\
    18 & 1181853 & 0.0727 & 6.823e-04 \\ 
    \bottomrule
\end{tabular}